\documentclass[11pt]{article} 
\usepackage[utf8]{inputenc}
\usepackage[T1]{fontenc}
\usepackage[english]{babel}

\usepackage{amsmath,amsfonts,amssymb,amsthm}
\usepackage{enumitem}
\usepackage{mathtools}
\usepackage{booktabs}
\usepackage[normalem]{ulem}
\usepackage{authblk}

\DeclareSymbolFont{symbolsC}{U}{txsyc}{m}{n}
\DeclareMathSymbol{\Searrow}{\mathrel}{symbolsC}{117}

\usepackage[mono=false]{libertine}
\usepackage[cmintegrals,libertine]{newtxmath}
\usepackage[cal=euler, scr=boondoxo]{mathalfa}

\useosf
\usepackage[a4paper,vmargin={3.5cm,3.5cm},hmargin={2.5cm,2.5cm}]{geometry}
\usepackage[font={small,sf}, labelfont={sf,bf}, margin=1cm]{caption}
\usepackage[pdftex,colorlinks=true]{hyperref}
\usepackage[pdftex]{color,graphicx}

\usepackage{stackrel}
\usepackage{soul} 

\theoremstyle{plain}
\newtheorem{theorem}{Theorem}

\newtheorem{proposition}[theorem]{Proposition}
\newtheorem{lemma}[theorem]{Lemma}
\theoremstyle{definition}

\newcommand{\Hnr}{{ \bar{\mathcal{H}} }}

\newcommand{\Onr}{{ \bar{\mathcal{O}} }}

	\title{\bf A family of triangulated $3$-spheres \\ constructed from trees}
\author[1]{Timothy Budd}
\author[1,2]{Luca Lionni}
\affil[1]{IMAPP, Radboud University, Nijmegen, The Netherlands.}
\affil[2]{Heidelberg University, Institut für Theoretische Physik, Philosophenweg 19, 69120 Heidelberg, Germany. }
\affil[ ]{\textnormal{\texttt{\href{mailto:t.budd@science.ru.nl}{t.budd@science.ru.nl}, \href{mailto:lionni@thphys.uni-heidelberg.de}{lionni@thphys.uni-heidelberg.de}}}}
\date{\today}
	
\begin{document}
	\maketitle
\begin{abstract}
The search for universality in random triangulations of manifolds, like those featuring in (Euclidean) Dynamical Triangulations, is central to the random geometry approach to quantum gravity.
In case of the 3-sphere, or any other manifold of dimension greater than two for that matter, the pursuit is held back by serious challenges, including the wide open problem of enumerating triangulations.
In an attempt to bypass the toughest challenges we identify a restricted family of triangulations, of which the enumeration appears less daunting.
In a nutshell, the family consists of triangulated 3\nobreakdash-spheres decorated with a pair of trees, one spanning its tetrahedra and the other its vertices, with the requirement that after removal of both trees one is left with a tree-like 2-complex. 
We prove that these are in bijection with a combinatorial family of triples of plane trees, satisfying restrictions that can be succinctly formulated at the level of planar maps. 
An important ingredient in the bijection is a step-by-step reconstruction of the triangulations from triples of trees, that results in a natural subset of the so-called locally constructible triangulations, for which spherical topology is guaranteed, through a restriction of the allowed moves.
We also provide an alternative characterization of the family in the framework of discrete Morse gradients. 
Finally, several exponential enumerative bounds are deduced from the triples of trees and some simulation results are presented.
\end{abstract}

	\tableofcontents

\section{Introduction}

\subsection{Universality classes of random geometry} 
The geometry of uniform random triangulations of the 2-sphere, and of more general models of random planar maps, has been the focus of intensive research at the interface of theoretical physics, combinatorics and probability theory. 
A driving question behind these developments has been whether natural notions of continuous random metrics exist on the 2-sphere, in analogy to the Wiener measure on paths on the real line. 
This question can be approached by starting from random discrete geometries, encoded by random planar maps, and seeking scaling limits of the geometry, in which the size of the building blocks decreases while their number increases. 
It has been answered in the affirmative in the form of the \emph{Brownian sphere}, a random continuous metric on the 2-sphere with interesting fractal properties, like a Hausdorff dimension of 4 almost everywhere. 
It is not just the scaling limit of uniform triangulations of the 2-sphere \cite{LeGall_Uniqueness_2013}, but of many more models of random planar maps, see e.g. \cite{LeGall_Uniqueness_2013,Miermont_Brownian_2013,Bettinelli_scaling_2014,Marzouk_Scaling_2018}. 
As all these models share the same scaling limit, they are said to belong to the same \emph{universality class}, terminology derived from the phenomenon of universality in statistical physics which states that the macroscopic laws are often largely independent of the microscopic details of a model.

Besides striving to better understand the properties of the Brownian universality class, effort has gone into finding other classes. Insights from physics tell us that a natural way to change the universality class is by including critical matter systems that interact with the geometry. 
This can be achieved by random sampling of a planar map together with a statistical system supported on the map, like an Ising model or a uniform spanning tree.
It should open up the possibility of scaling limits belonging to a 1-parameter family of universality classes, that are related in the continuum to Liouville Quantum Gravity at different values of its coupling constant $\gamma$ \cite{Gwynne_Existence_2021,Ding_Introduction_2021}. 
The Brownian universality class, often referred to as "pure gravity" in the physics literature, corresponds to the special value $\gamma = \sqrt{8/3}$.
Although full scaling limits (in the Gromov-Hausdorff sense) are still to be established, there is a long list of models that are expected to belong to non-Brownian universality classes (i.e.\ universality classes different from the Brownian one) based on weaker convergence results.
For example, triangulations coupled to a critical Ising model and triangulations decorated by a uniform spanning tree are expected to belong to the universality classes of Liouville Quantum Gravity with $\gamma = \sqrt{3}$ and $\gamma = \sqrt{2}$ respectively, see e.g. \cite{Albenque_Geometric_2022} and \cite{Gwynne_Mating_2019}.
Several other families of universality classes (not corresponding to Liouville Quantum Gravity) have been identified in planar map models, like the stable maps \cite{LeGall_scaling_2010}, stable shredded spheres \cite{Bjoernberg_Stable_2019}, stable looptrees \cite{Curien_Random_2014}, but although they are expected to be planar (homeomorphic to a subset of the 2-sphere) none of these random metric spaces have the topology of the $2$-sphere, essentially because they contain macroscopic holes. A conjectural family of universality classes that may have the desired $2$-sphere topology was proposed in \cite{Budd_Geometry_2017,Bertoin_Martingales_2018} and is in some sense dual to the family of stable maps of \cite{LeGall_scaling_2010}. 
Finally, a
sequence of random metric spaces obtained as limits of random planar maps with unconventional distributions, with Hausdorff dimensions $2^D$ for any $D\ge 3$, has been put forward in \cite{Lionni_Iterated_2019} that may display spherical topology.

Although many scaling limits are yet to be established rigorously and many of their properties to be determined, one may say that the arena of random geometries on the 2-sphere is well populated.
This cannot be said of its higher-dimensional counterpart.
To our knowledge, \emph{not a single non-trivial universality class with the topology of the 3-sphere, or any other manifold of dimension three or larger, is currently known}, despite decades of research spurred by the relevance of such a class to the problem of quantum gravity in physics (see Section~\ref{sec:qg}).
In analogy with the two-dimensional case, a natural approach is to focus on triangulations of the 3-sphere, i.e. piecewise linear geometries with the topology of $S^3$ assembled from identical tetrahedra that are glued along their boundary triangles. 
One could choose such a triangulation with a fixed number $n$ of tetrahedra uniformly and consider the large-$n$ limit. 
Contrary to the planar case, the number of vertices of the triangulation is not determined by $n$, so it is natural to consider the more general distribution in which probabilities are proportional to $x^{N}$ where $N$ is the number of vertices of a triangulation.
Modulo some choices on the class of triangulations allowed, e.g. whether tetrahedra are allowed to be glued to themselves, this one-parameter family of models is known in the physics literature as three-dimensional Dynamical Triangulations \cite{Ambjoern_Three_1991}.

\subsection{Challenges on the road to higher dimensions}
However, investigation of this model has not yet led to new universality classes, neither rigorously nor conjecturally based on numerical evidence. 
The challenges in this direction can be summarized as:
\begin{enumerate}[label = {(C\arabic*)}]
	\item\label{item:topology} \textbf{Topology of manifolds in three dimensions is a lot more complicated than in two.} Although algorithms exist that recognize 3-sphere triangulations \cite{Rubinstein_algorithm_1995,Thompson_Thin_1994}, no such polynomial-time algorithm is known. This is in stark contrast with the two-dimensional case, where spherical topology can be recognized by a simple computation of the Euler characteristic. The state of the art is that 3-sphere recognition is in NP, meaning that for every triangulation of the 3-sphere there exists a certificate that can be checked in polynomial time \cite{Schleimer_Sphere_2011}.
	\item\label{item:enumeration} \textbf{Enumerating triangulations of the 3-sphere is hard.} This difficulty is illustrated by the still wide open problem whether there exists an exponential bound in $n$ on the number of 3-sphere triangulations with $n$ tetrahedra. The formulation of this problem goes back at least to Durhuus and Jonsson \cite{Durhuus-Jonsson} and was popularized by Gromov \cite{Gromov_Spaces_2000}. Various subclasses of triangulated spheres with exponential bounds have been identified \cite{Benedetti-Ziegler}, notably the locally constructible class \cite{Durhuus-Jonsson}, that we will discuss extensively in this paper.\\
	The enumeration problem of 3-dimensional triangulations has been extensively studied in the context of generalizations of matrix models to rank-3 tensor models (see e.g.~\cite{Gurau2017, Lionni:17}).
	These are related through formal expansions to generating functions of $3$-dimensional pseudo-manifolds organized by their number of tetrahedra and certain non-negative invariants. 
	In the case of matrix models the expansion is organized by the genus, and in the large-matrix limit the matrix integrals provide access to the enumeration of maps on the 2-sphere.
	In comparison, for tensor models the invariants are non-topological and therefore disentangling the enumeration of spherical triangulations remains a challenge. 
	The limit of large tensors only gives access to a small subset of triangulated 3-spheres that are treelike. Moreover, analytical tools to compute tensor integrals are still lacking.
	\item\label{item:critical} \textbf{Simulations have not uncovered promising critical phenomena.} In lack of analytical results, one may still build intuition based on numerics. Dynamical triangulations have been the subject of extensive Monte Carlo simulations in the physics literature, mainly in the nineties \cite{Boulatov_phase_1991,Ambjoern_Three_1992,Ambjoern_vacuum_1992,Catterall_Entropy_1995,Hagura_Phases_1998,Hotta_Multicanonical_1998,Thorleifsson_Three_1999}. First of all, the numerics lend support to the existence of an exponential bound on the number of 3-sphere triangulations. Secondly, the phase diagram, parametrized by the weight $x$ per vertex, has been explored and two phases have been identified separated by a phase transition at $x=x_*$: a \emph{branched polymer phase} ($x>x_*$) and a \emph{crumpled phase} ($x<x_*$). In the branched polymer phase, large random triangulations appear to be structured in a treelike fashion, and it has been conjectured that for any $x>x_*$, the model lives in the universality class of the Continuum Random Tree (CRT) of Aldous \cite{Aldous_continuum_1991}, which arises as the universal scaling limit of random discrete plane trees.\footnote{In the case of colored triangulations, this has been shown to hold for $x=\infty$ \cite{Gurau_Melons_2014}, where only the triangulations with a maximal number of vertices contribute. These triangulations arise naturally in the context of random tensor models in the limit of large tensors.} In the crumpled phase, the edge graph of random triangulations is highly connected and graph distances appear to grow only slowly with the number of tetrahedra $n$, in other words portraying some ``small world phenomena''. This could indicate that no scaling limit exists or, if it does exist, that it has a high Hausdorff dimension. Based on numerical data it is hard to distinguish these cases, because one needs very large size $n$ to reach sizeable graph distances.
	Finally, the phase transition at $x=x_*$ appears to be first order, so no new critical phenomena are expected to be found there.  
\end{enumerate}

\subsection{Guiding principles}
To keep hopes of formulating an analytically tractable model of random triangulations of the 3-sphere belonging to a new universality class, without fully resolving these challenges, it is necessary to explicitly circumvent them. 
This naturally leads to the following guidelines:
\begin{enumerate}[label = {(G\arabic*)}]
	\item\label{item:exponentialsubset} In view of the challenge \ref{item:enumeration}, the model's probability distribution should have support only on a proper subset of all triangulated $3$-spheres, for which the existence of an exponential bound is within reach.\\
	There are two reasons to opt for this, the first of which is purely pragmatic. As mentioned above there is good numerical support that triangulated 3-spheres are actually exponentially bounded, so the lack of a rigorous exponential bound on some subset of all triangulated $3$-spheres is quite likely a consequence of our ignorance about enumerating that subset.
	It goes without saying that we wish to avoid such ignorance in our quest for an analytically tractable model. 
	The second reason is that a probability distribution without size constraints on a super-exponential set of triangulations cannot arise from simple local Boltzmann weights\footnote{A local Boltzmann weight means a product over all tetrahedra of a function that only depends on the local geometry around each tetrahedron.} like the weight $z^nx^N$ used in Dynamical Triangulations for a triangulation with $n$ tetrahedra and $N$ vertices. 
	\item\label{item:certifiedtopology} For this subset of triangulations, the 3-sphere topology should be efficiently verifiable. In view of the first challenge \ref{item:topology}, this suggests that a certificate of $3$-sphere topology should be easily deduced from the triangulation or incorporated as an additional piece of information in the configuration.
	\item\label{item:departDT} The conclusions one could draw from \ref{item:critical} are rather subjective, but we settle on the following: preferably the distribution of the model should be sufficiently distinct from that of Dynamical Triangulations, that qualitatively different phenomena appear in simulations. In particular, one should be on the lookout for qualitatively different phases and higher-order phase transitions.
\end{enumerate}
A final guiding principle underlying this work is that it would be beneficial if the model would feature (combinatorial) trees as building blocks. 
There are multiple reasons for this. 
First of all, random trees and their scaling limits have been studied extensively. 
The Continuum Random Tree (CRT) of Aldous \cite{Aldous_continuum_1991} is perhaps the best understood universality class of random geometries around.
Second, in many cases in combinatorics, the existence of an exact enumeration result of a combinatorial family can be explained bijectively by uncovering an encoding of the family in terms of trees or collections of trees (perhaps with decorations).
For many families of planar maps, such encodings are known, with famous examples being the Cori-Vauqelin-Schaeffer bijection between quadrangulations and labeled plane trees \cite{Cori_Planar_1981,Schaeffer_Conjugaison_1998} and its generalization by Bouttier, Di Francesco and Guitter to arbitrary planar maps with control on face degrees \cite{Bouttier_Planar_2004}. 
Beyond enumeration, these trees provide insight into statistics of graph distances in the maps.
The combination of these properties and that the trees live in the universality class of the CRT has played a crucial role in establishing scaling limits of random planar map models and the construction of the Brownian sphere out of the CRT \cite{Marckert_Limit_2006,LeGall_Uniqueness_2013}. 
More recently, it has been realized \cite {Duplantier_Liouville_2014,Gwynne_Mating_2019} that all other universality classes in the family of Liouville Quantum Gravity can also be assembled from CRTs via the so-called mating of trees approach.
As a consequence one may understand the CRT as the building block of essentially all known universality classes of random geometries on the 2-sphere. 
With this background in mind it seems natural to focus the search for new universality classes at geometries that can be encoded in treelike structures:
\begin{enumerate}[label = {(G\arabic*)}, resume]
	\item\label{item:treeencoding} Preferably the configurations of the model, meaning the triangulations together with decorations, can be encoded in trees.
\end{enumerate}

In \cite{Lionni_Iterated_2019} the second author and Marckert introduced a model of non-planar discrete random graphs taking the last principle seriously: by iterating the Cori-Vauqelin-Schaeffer bijection these graphs were encoded by three or more trees.
Several partial results hint at the existence of new scaling limits, even though convergence in the Gromov-Hausdorff sense is still out of reach. 
However, it is not yet clear whether the featured graphs have the structure of a triangulated manifold (or more general cell-complexes) while the topologies of the potential scaling limits remain unknown.

\

\emph{In this paper, we propose a new model that certainly satisfies \ref{item:exponentialsubset}, \ref{item:certifiedtopology} and \ref{item:treeencoding} and for which preliminary simulations give indications of \ref{item:departDT}.}

\subsection{Relevance for Quantum Gravity}\label{sec:qg}
Going by the significant challenges faced in the search of universality classes on higher-dimensional manifolds, one may be inclined to label the search as premature. 
However, in light of the great importance of the problem to fundamental physics, any progress in this direction, even at a conjectural level, is of value. 
Let us take a minute to introduce the high-level motivation coming from the search for quantum gravity, especially for those that are not so familiar with the challenges in theoretical high-energy physics. 

Our best understanding of the elementary particles and fundamental forces of physics, excluding gravity, is encapsulated in the Standard Model of particle physics, a quantum field theory living on 4-dimensional Minkowski space. 
Of course, this theory is far from a rigorous mathematical definition in a constructive sense (e.g. as a measure on an appropriate function space). 
But when formulated at the level of perturbation theory, a consistent computational framework allows one to make predictions on the outcomes of scattering experiments, once the numerical values of the coupling constants appearing in the Standard Model action (at an appropriate energy scale) have been inferred from previous experiments.
Needless to say, these predictions have been extraordinarily successful at the energy scales that are currently within reach of experiment, lending strong support to the Standard Model.
As with all interacting quantum field theories, the values of the coupling constants are associated to an energy scale and their dependence on this scale is governed by a renormalization group flow. 
Thanks to the perturbative renormalizability of the Standard Model, renormalization can be implemented within the fixed finite-dimensional space of values for the couplings constants of the model.
As a consequence the predictive power of perturbation theory persists in principle to arbitrarily high energies.

The situation is different when, in an attempt to predict quantum effects in the gravitational force,  perturbation theory is applied to the field corresponding to the Lorentzian metric on spacetime that features in Einstein's theory of general relativity.
This perturbation theory is not renormalizable, meaning that the renormalization group flow cannot be confined to a finite-dimensional space comprising of Newton's constant, the cosmological constant and a finite number of couplings constants associated to higher-order curvature invariants.
This is problematic from a fundamental physics point of view, because in order to fully characterize the theory one would have to experimentally determine the values of an infinite number of coupling constants.
A way out of this unfortunate situation would be if by some principle the values of the infinitely many couplings constants were a priori constrained to lie on a finite-dimensional subspace.
In that case only a finite number of experiments would suffice to determine the theory, and predictivity would be restored.

This is precisely what the Asymptotic Safety scenario for gravity would accomplish if realized.
In this scenario, proposed by Weinberg in the seventies \cite{Weinberg_Ultraviolet_1979}, the non-perturbative quantum field theory possesses an ultraviolet fixed point, meaning that the renormalization group flow of the (dimensionless) couplings of the theory approach fixed values in the limit of high energies.
If the so-called critical surfaces, i.e. the union of renormalization trajectories ending at this ultraviolet fixed point, is finite dimensional, then constraining the model to be on this surface provides precisely the restoration of predictivity needed. 
The main support for the scenario is provided by functional renormalization group methods \cite{Reuter_Nonperturbative_1998}, which rely on studying truncations of the renormalization group flow onto finite-dimensional subspaces. 
Consistently these have provided evidence for the existence of non-perturbative ultraviolet fixed point satisfying the desired properties (see \cite{Reuter_Quantum_2012,Reuter_Quantum_2019} for an overview). 
While from a mathematical viewpoint it is hard to evaluate the significance of these results for a full-fledged quantum field theory of the spacetime metric, it is certainly a scenario to be reckoned with. 

What would the Asymptotic safety scenario mean for the geometry of spacetime? Let us focus on the situation of a Euclidean Quantum Field Theory of a four-dimensional Riemannian metric, which has been the subject of most functional renormalization group investigations in gravity. 
One might hope that such a theory exists in a constructive sense as a family of probability measures, one for each choice of coupling constants, on some space of suitably generalized Riemannian metrics on $\mathbb{R}^4$.
In such a realization, the renormalization group flow would amount to a pull-back of the measure under an overall scaling of the metrics, and a fixed point of this flow would thus correspond to a probability measure that is invariant under scaling.
Leaving aside the question of what the appropriate notion of generalized metric structure is in the constructive framework, such a scale-invariant probability distribution can be interpreted as a \emph{universality class of random geometry}. 
Vice versa, an explicit construction of such a universality class would be a candidate for a fixed point in a yet to be determined family of quantum field theories, such that on the critical surface each of these features the scale-invariant geometry in the high-energy or small-distance limit.

\subsection*{Acknowledgments}

This work is part of the START-UP 2018 programme with project number 740.018.017, which is financed by the Dutch Research Council (NWO). TB also acknowledges support from the VIDI programme with project number VI.Vidi.193.048, which is financed by the Dutch Research Council (NWO). L.L. is now supported by the European Research Council (ERC) under the European Union’s Horizon 2020 research and innovation program (grant agreement No818066) and by Deutsche Forschungsgemeinschaft (DFG, German Research Foundation) under Germany's Excellence Strategy  EXC-2181/1 - 390900948 (the Heidelberg STRUCTURES Cluster of Excellence). 

\section{Model introduction and main results}

\subsection{Triple-trees}\label{sec:triptrees}
\begin{figure}[h]
	\centering
	\includegraphics[width=\linewidth]{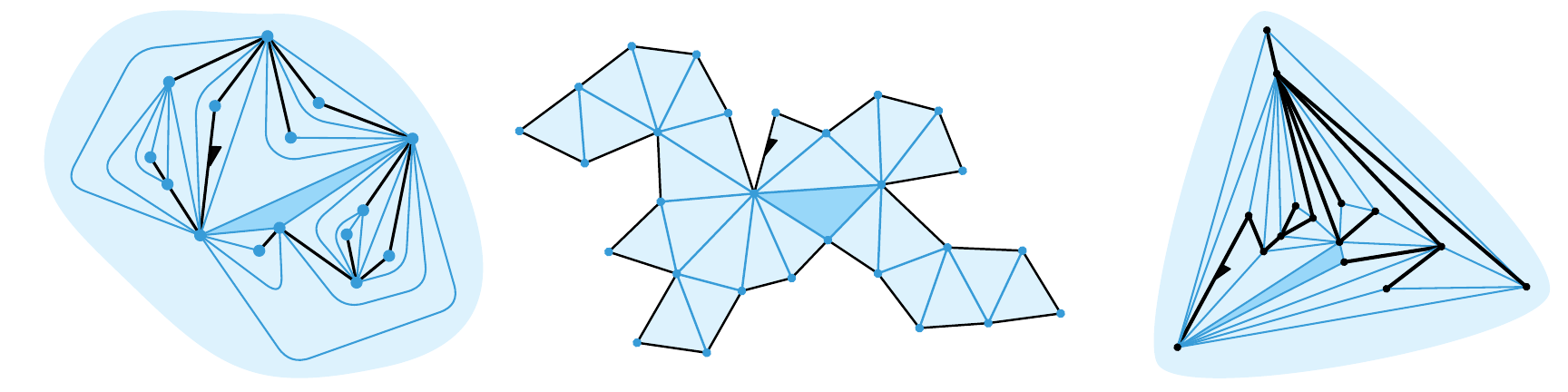}
	\caption{An example of an outerplanar triangulation in $\mathcal{O}_n$ with $n=14$ (middle) that admits a gluing into a hierarchical triangulation (left) and an Apollonian triangulation (right). The darker shaded triangles are in correspondence and are an example of the triangles $\theta_{_\mathrm{H}}$, $\theta$, $\theta_{_\mathrm{A}}$ introduced in the text. \label{fig:hierarchicalapollonian}}
\end{figure}

We start with several definitions.
A \emph{planar map} is a connected multi-graph, i.e.\ a connected graph with loops and multiple edges allowed, drawn on the $2$-sphere without intersecting edges and viewed modulo orientation-preserving homeomorphisms.
A \emph{rooted} planar map is a planar map together with a distinguished oriented edge. 
The connected regions in the complement of a planar map are called \emph{faces} and the \emph{degree} of a face is the number of sides of edges that bound the face.
A face of degree $k$ is \emph{simple} if it is incident to $k$ distinct vertices.
A planar map in which all faces have degree three is called a \emph{planar triangulation} (examples appear on the left and right of Fig.~\ref{fig:hierarchicalapollonian}).

An \emph{outerplanar triangulation of the $n$-gon} is an $n$-vertex rooted planar map whose face on the right of the root, the \emph{outer} face, is simple and has degree $n$ while all other faces are of degree three, implying that all  vertices of an outerplanar triangulation are incident to the outer face (see the example in the middle of Fig.~\ref{fig:hierarchicalapollonian}). An outerplanar triangulation with $n$ vertices has $n-2$ triangles. 
Let us denote the set of outerplanar triangulations with $2n$ vertices by $\mathcal{O}_n$.
A non-crossing pairing of size $n$ is a pairing of $\{1,2,\ldots,2n\}$ such that there exist no $a < b < a' < b'$ such that $a$ is paired with $a'$ and $b$ with $b'$. 
The set of non-crossing pairings of size $n$ is denoted by $\mathcal{P}_n$.
Given an outerplanar triangulation $t\in \mathcal{O}_n$ and a non-crossing pairing $\pi \in \mathcal{P}_n$ one may naturally construct a planar triangulation with $n+1$ vertices by labeling the $2n$ edges on the outer face (called boundary edges) of $t$ from $1$ to $2n$ (in clockwise order starting from the root edge) and pairwise gluing these edges according to $\pi$.
After gluing, these $n$ edges form a spanning tree in the resulting triangulation.
Denoting by $\mathcal{S}_n$ the set of rooted planar triangulations with $n+1$ vertices and a distinguished spanning tree (it is said to be \emph{tree-decorated}) that includes the root edge, this determines a bijection 
\[
 \mathsf{Glue} : \mathcal{O}_n \times \mathcal{P}_n \to \mathcal{S}_n.
 \]
Elements of $\mathcal{S}_n$ will be called \emph{(rooted) tree-decorated triangulations}.

We will consider two natural classes of triangulations that are obtained by recursive subdivision starting from the unique loopless triangulation with two triangles.
Repeatedly selecting an edge and inserting a pair of triangles glued along two of their sides one obtains the family that we call \emph{hierarchical triangulations} (also called \emph{melonic}  triangulations \cite{Gurau2017}), while repeated star-division of a triangle into triples of triangles yields the family of \emph{Apollonian triangulations}.
Let $\mathcal{H}_n \subset \mathcal{S}_n$ be the set of tree-decorated hierarchical triangulations with $n+1$ vertices with the extra condition that the spanning tree cannot contain edges that are contained in a cycle of length $2$.
Let $\mathcal{A}_n \subset \mathcal{S}_n$ be the set of tree-decorated Apollonian triangulations with $n+1$ vertices without further restrictions.

\begin{figure}[h]
	\centering
	\includegraphics[width=.5\linewidth]{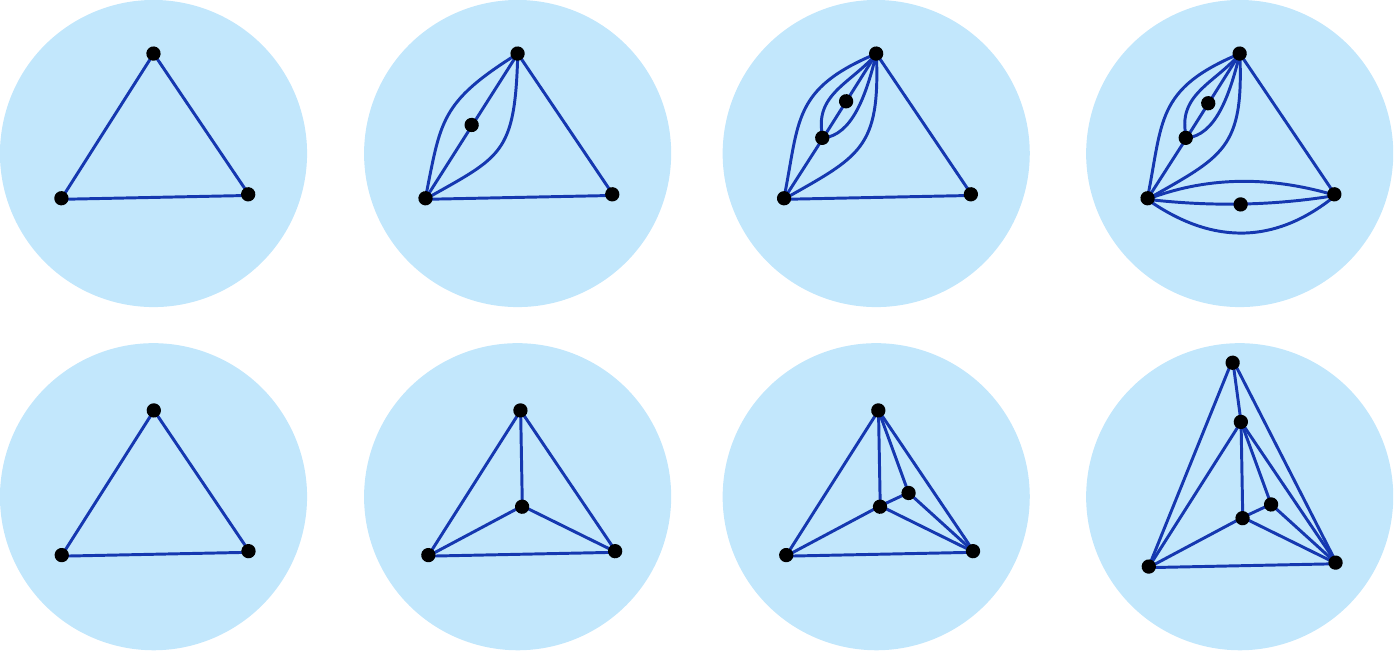}
	\caption{Two families of triangulations obtained by recursive subdivisions: hierarchical triangulations (top) and Apollonian triangulations (bottom).}
\end{figure}

We are now ready to introduce our main combinatorial family:
\begin{equation}
\label{eq:def-triple-trees}
	\mathcal{M}_n = \bigl\{ (t,\pi_{_\mathrm{H}},\pi_{_\mathrm{A}}) \in \mathcal{O}_n \times \mathcal{P}_n \times \mathcal{P}_n : \mathsf{Glue}(t,\pi_{_\mathrm{H}}) \in \mathcal{H}_n\text{ and }\mathsf{Glue}(t,\pi_{_\mathrm{A}}) \in \mathcal{A}_n\bigr\},
\end{equation}
i.e.\ the set of outerplanar triangulations with a pair of non-crossing pairings such that gluing according to the first gives a hierarchical triangulation, while gluing along the second gives an Apollonian triangulation.
We refer to an element of $\mathcal{M}_n$ as a \emph{triple-tree}, because both an outerplanar triangulation and a non-crossing pairing have the structure of a rooted plane tree. An example is shown in Fig.~\ref{fig:hierarchicalapollonian}.

\subsection{Associated 3-dimensional triangulations.}
\begin{figure}[h!]
	\centering
	\includegraphics[width=\linewidth]{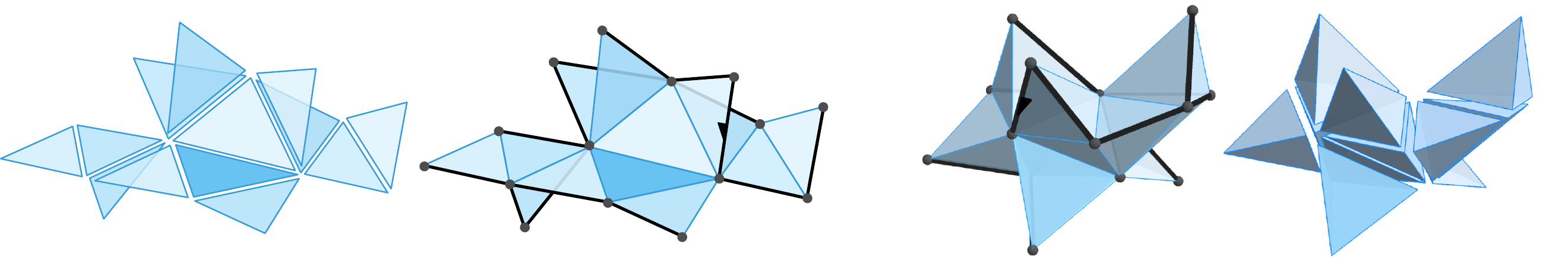}
	\caption{The hierarchical triangulation of Fig.~\ref{fig:hierarchicalapollonian} represented as a gluing of triangles (left) and the Apollonian triangulation as a gluing of tetrahedra (right). 
	\label{fig:apollonian3D}}
\end{figure}

In both cases there is a natural three-dimensional perspective on these maps. Informally for the moment, it is useful to think of the 2-sphere on which a planar map lives as embedded in a 3-sphere, such that each face of the map naturally acquires a top side and a bottom side induced by the orientation of the 2-sphere. Then an Apollonian triangulation can be viewed as the boundary of a three-dimensional simplicial ball consisting of stacked tetrahedra, called a \emph{tree of tetrahedra}, as illustrated on the right-hand side of Fig.~\ref{fig:apollonian3D} (the Apollonian triangulation can be wrapped around the tree of tetrahedra with the bottom side of triangles facing the tetrahedra). 

A \emph{tree of triangles} is a 2-dimensional simplicial complex made of triangles glued along their edges in a tree-like fashion and embedded in the 3-sphere  (so that there is a cyclic ordering of the triangles around every edge), as illustrated on the left-hand side of Fig.~\ref{fig:apollonian3D}. 
A hierarchical triangulation $h$ can be viewed as the contour of a unique tree of triangles in the sense that $h$ can be wrapped around it, but this time with the top sides of the triangles on the inside. It is because of the latter convention that, the Fig.~\ref{fig:apollonian3D} appears mirrored compared to the hierarchical triangulation in Fig.~\ref{fig:hierarchicalapollonian}.
The tree of triangles corresponding to $h$ is denoted by $\mathsf{Id}[h]$, because informally it results from \emph{identifying} pairs of triangles of $h$ into the triangles of the tree.
Note in particular that each triangle of $\mathsf{Id}[h]$ has two pre-images in the hierarchical triangulation.

In a triple-tree $(t, \pi_{_\mathrm{H}}, \pi_{_\mathrm{A}})$, each triangle $\theta$ of $t$ is at the same time a triangle $\theta_{_\mathrm{A}}$ on the boundary $ \mathsf{Glue}(t,\pi_{_\mathrm{A}})$ of a tree of tetrahedra, and a triangle  $\theta_{_\mathrm{H}}$ around the tree of triangles $\mathsf{Id}[h]$, for $h=\mathsf{Glue}(t,\pi_{_\mathrm{H}})$, as illustrated in Fig.~\ref{fig:apollonian3D} for the example of Fig.~\ref{fig:hierarchicalapollonian}.

\begin{figure}[h!]
	\centering
	\includegraphics[width=.7\linewidth]{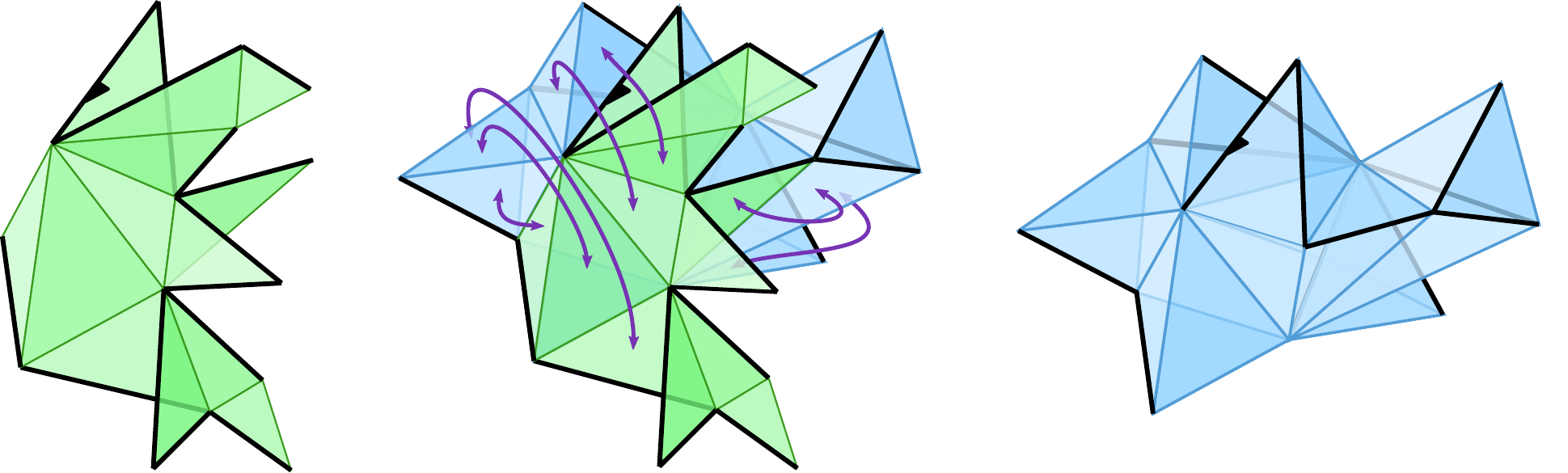}
	\caption{A partial gluing of the tree of triangles and tree of tetrahedra of Fig.~\ref{fig:apollonian3D}. 
	\label{fig:treeglue}}
\end{figure}

If for each triangle $\theta$ of $t$, we glue together its two copies $\theta_{_\mathrm{A}}$ and $\theta_{_\mathrm{H}}$ while matching the duplicated edges and vertices (Fig.~\ref{fig:treeglue}), we obtain a closed 3-dimensional triangulation $T$. The edges of the distinguished spanning trees of $ \mathsf{Glue}(t,\pi_{_\mathrm{A}})$ and of $\mathsf{Id}[ \mathsf{Glue}(t,\pi_{_\mathrm{H}})]$ are identified in groups in $T$, 
each group resulting in a distinguished edge of $T$. We will show that this subset of distinguished edges $E$ forms a spanning tree of $T$.
Three trees therefore also appear in $T$:
\begin{itemize}
\item A spanning tree of edges $E$, 
\item A spanning tree of tetrahedra $T_0$ (whose boundary $\partial T_0$ is $ \mathsf{Glue}(t,\pi_{_\mathrm{A}})$), 
\item If in $T$ we remove $T_0$ and then split open all the edges of $E$, then the result is the tree of triangles $\mathsf{Id}[ \mathsf{Glue}(t,\pi_{_\mathrm{H}})]$.
\end{itemize}
The last assertion is not obvious and will be shown in Sec.~\ref{sec:origin-in-3d}.

\begin{figure}[h!]
	\centering
	\includegraphics[width=.7\linewidth]{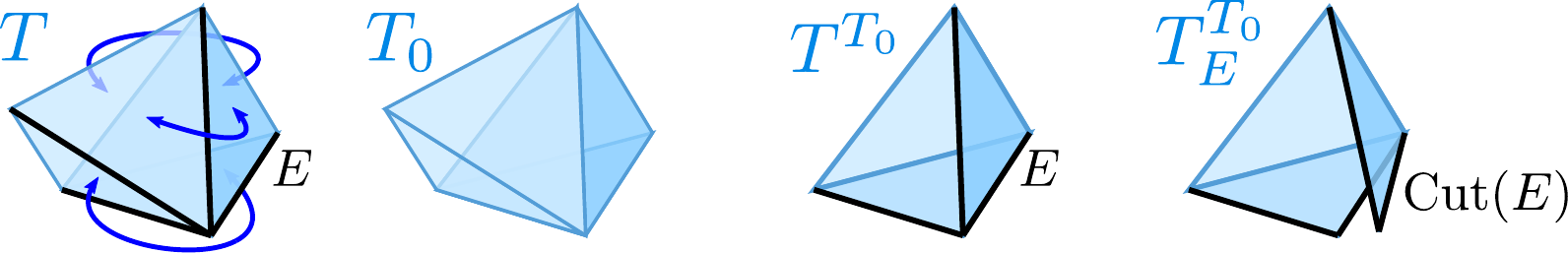}
	\caption{A triangulation $T$ of the 3-sphere with two tetrahedra and a spanning tree of edges $E$ and a spanning tree of tetrahedra $T_0$. The 2-complexes $T^{T_0}$ and $T^{T_0}_E$ each consist of three triangles. This example corresponds to the top-right triple tree of size $n=4$ in Fig.~\ref{fig:tripletreesn4n5}.  
	\label{fig:triangulationcutting}}
\end{figure}

More precisely, a spanning tree of tetrahedra of a 3-dimensional triangulation $T$ consists of all the tetrahedra, glued along the triangles dual to the edges of a spanning tree of the dual graph of $T$ (its internal triangles).  Replacing the tetrahedra and internal triangles of a spanning tree of tetrahedra $T_0$  of $T$ by a single 3-cell with the same boundary, one obtains a 2-complex embedded in the 3-sphere, which we denote by  $T^{T_0}$. 
If $E$ is a spanning tree of edges in $T$, then it determines one in $T^{T_0}$ as well, since $T^{T_0}$ has the same set of vertices and edges as $T$. We denote by $T^{T_0}_E$ the embedded 2-complex obtained from $T^{T_0}$ by duplicating every edge $e$ of $E$ so that every triangle of $T^{T_0}$ that contained $e$ contains a copy of $e$ in  $T^{T_0}_E$ and every copy of $e$ now belongs to a single triangle;  and only keep two triangles attached through a vertex if these triangles also share an edge containing this vertex. Let $\mathsf{Cut}(E)$ be the set of edges of $T^{T_0}_E$ that are duplicates of edges of $E$. 
See Fig.~\ref{fig:triangulationcutting} for an example.

\begin{figure}[h!]
	\centering
	\includegraphics[scale=0.65]{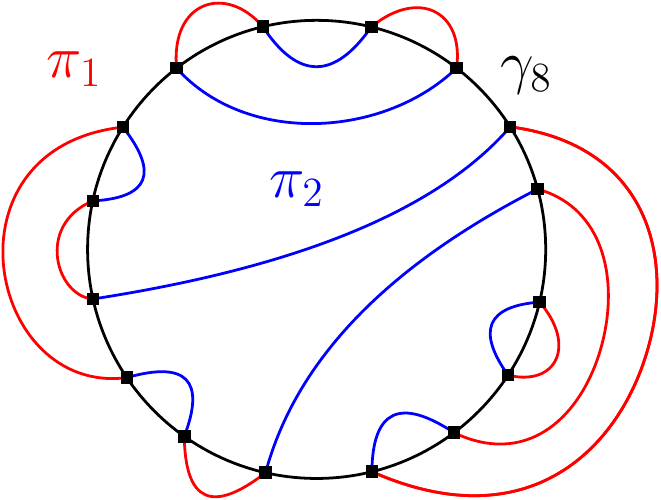}
	\caption{A meander system with 16 vertices and three loops.}
	\label{fig:ex-meander}
\end{figure}

Given two non-crossing pairings $\pi_1, \pi_2$ of $\{1,\ldots, 2n\}$ for $n\ge 1$, we also define the \emph{meander system} $[\pi_1, \pi_2]$ as the planar map obtained as follows (Fig.~\ref{fig:ex-meander}): consider a circle $\gamma_n$ drawn on the 2-sphere without crossings, with $2n$ vertices labeled from $1$ to $2n$  so that two consecutive vertices have consecutive labels. This delimitates two zones. In the interior of one of the  zones, draw an edge between every two vertices if the labels are paired in $\pi_1$ so that the edges don't cross, and do the same in the interior of the other zone for $\pi_2$, and then remove the labels. We call \emph{loops of the meander system} $[\pi_1, \pi_2]$ the connected components of the map obtained from $[\pi_1, \pi_2]$ by keeping only the edges corresponding to $\pi_1$ an $\pi_2$ (and not the edges of the cycle $\gamma_n$). 

\subsection{Main bijective result}

The following theorem summarizes some of the main results of the paper. A 3-dimensional triangulation is said to be rooted if it has a marked oriented edge, and a marked triangle not in $T_0$ (that is, which also belongs to $T^{T_0}$), among the triangles containing this edge.
\begin{theorem}
\label{th:thm-intro}
There is a bijection between triple-trees with $2n$ triangles,  and 3-dimensional triangulations $T$ with $n-1$ tetrahedra, with a marked spanning tree of tetrahedra $T_0$ and a marked spanning tree of edges $E$, rooted on an edge of $E$ and such that the 2-complex $T^{T_0}_E$ is a tree of triangles and $\mathsf{Cut}(E)$ is a spanning tree of $T^{T_0}_E$. 

Furthermore, letting $(t,\pi_{_\mathrm{H}}, \pi_{_\mathrm{A}})$ be a triple-tree  and $T$ decorated by $T_0$ and $E$ be its image by this bijection:
\begin{enumerate}[label=$\blacktriangleright$]
\item $\partial T_0= \mathsf{Glue}(t,\pi_{_\mathrm{A}})$\quad and\quad $T^{T_0}_E= \mathsf{Id}\bigl[\mathsf{Glue}(t,\pi_{_\mathrm{H}})\bigr]$;
\item the number of vertices of $T$ minus one is the number of loops of the meander system $[\pi_{_\mathrm{H}}, \pi_{_\mathrm{A}}]$;
\item $T$ has the topology of the 3-sphere. 
\end{enumerate}
\end{theorem}
A more detailed version of this theorem will appear as Thm.~\ref{thm:thm-1} below, with the exception of the last statement. 
Regarding the topology, we will show  that these triangulations form a subset of the so-called \emph{locally constructible triangulations} (Sec.~\ref{sub:triple-trees-are-LC}), as well as a subset of the triangulations that admit a \emph{discrete Morse gradient with two critical simplices} (Sec.~\ref{sub:Morse}), both known to triangulate the 3-sphere. In Sec.~\ref{sub:TALC} and Sec.~\ref{sub:triple-trees-are-TALC}, we will give a precise characterization of the subset of locally constructible triangulations in bijection with triple-trees. 
 
Let us comment on how this combinatorial family fits our guiding principles.
Since outerplanar triangulations and non-crossing pairings are exponentially bounded, the same is true for triple-trees, thus fulfilling \ref{item:exponentialsubset}. 
We will collect more precise bounds in the next paragraph.
The corresponding three-dimensional triangulations come with a certificate of 3-sphere topology in the sense of \ref{item:certifiedtopology}, in the form of a local construction or, equivalently, as a discrete Morse gradient with two critical simplices. 
As explained, the triangulations are fully encoded in a triple of plane trees satisfying rather explicit constraints, thus adhering to \ref{item:treeencoding}.
Finally, it is too early to tell whether large (uniform) random triple-trees differ substantially from (uniform) triangulated 3-spheres (criterion \ref{item:departDT}).
However, based on preliminary numerical results, we suspect that the restriction from general triangulated 3-spheres to locally constructible ones is relatively mild compared to the further restriction from locally constructible triangulations to triple trees.
Extensive Monte Carlo simulation results will be presented in a forthcoming paper.

\subsection{Enumeration bounds}
\begin{figure}[t]
	\centering
	\includegraphics[width=\linewidth]{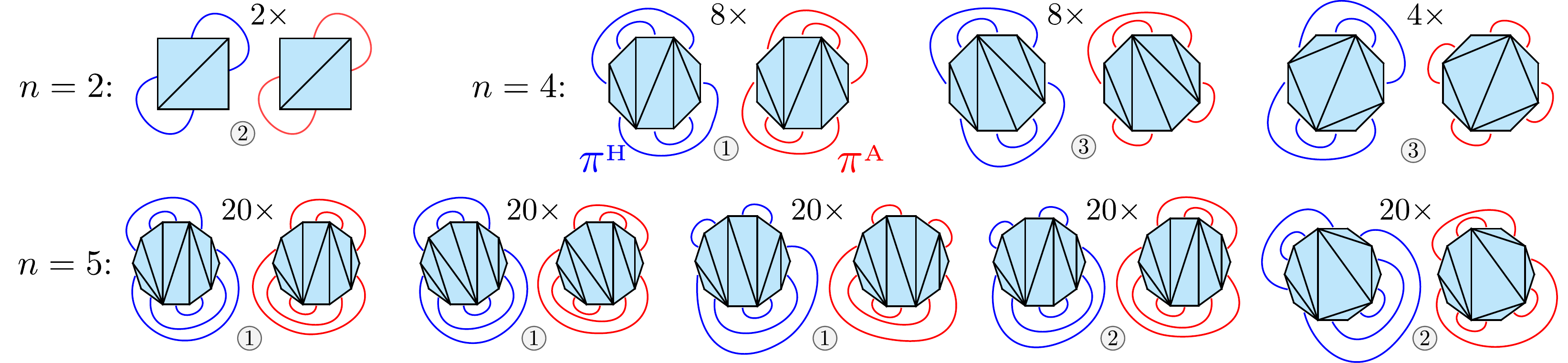}
	\caption{All triple trees of size $n\leq 5$ modulo rotations and reflections. The multiplicative factor indicates the number of inequivalent rotations and reflections, while the circled number counts the loops in the meander system $[\pi_{_\mathrm{H}}, \pi_{_\mathrm{A}}]$.}
	\label{fig:tripletreesn4n5}
\end{figure}

Since the triangulated 3-spheres under consideration are all locally constructible, and the number of locally constructible triangulations is exponentially bounded in $n$ \cite{Durhuus-Jonsson}, this holds for the triple-trees as well.
Although a precise enumeration of triple-trees is still out or reach, several explicit bounds can be obtained from the formulation in terms of trees.
The combinatorial families $\mathcal{O}_n$ (outerplanar triangulations) and $\mathcal{P}_n$ (non-crossing pairings) are straighforwardly enumerated in terms of the Catalan numbers $\operatorname{Cat}(n) = \frac{1}{n+1}\binom{2n}{n}$ as
\begin{align*}
	|\mathcal{O}_n| = \operatorname{Cat}(2n-2)  \sim \frac{16^n}{32\sqrt{2\pi} \,n^{3/2}},\qquad
	|\mathcal{P}_n| &= \operatorname{Cat}(n) \sim \frac{4^n}{\sqrt{\pi} \,n^{3/2}}.
\end{align*}
In Section \ref{sec:combinatorial-bounds}, Propositions \ref{prop:hierarchicalenum} and \ref{prop:apollonianenum}, we will show that the tree-decorated hierarchical triangulations $\mathcal{H}_n$ and tree-decorated Apollonian triangulations $\mathcal{A}_n$ have algebraic generating functions, from which one can explicitly determine the asymptotic enumeration to be 
\begin{align*}
	|\mathcal{H}_n| &= 2\cdot 3^{n-2}\operatorname{Cat}(n-2) \sim \frac{1}{72\sqrt{\pi}} \frac{12^n}{n^{3/2}},\\
	|\mathcal{A}_n| &\sim \frac{1}{(1752.10\ldots)} \frac{(28.43\ldots)^n}{n^{3/2}}.
\end{align*}
For triple trees, we introduce the bivariate generating function 
\begin{equation*}
	M(z,x) = \sum_{n=1}^\infty z^n M_n(x), \qquad M_n(x) = \sum_{(t,\pi_{_\mathrm{H}},\pi_{_\mathrm{A}})\in\mathcal{M}_n} x^{N(\pi_{_\mathrm{H}},\pi_{_\mathrm{A}})},
\end{equation*}
where $N(\pi_{_\mathrm{H}},\pi_{_\mathrm{A}})$ is the number of loops of the meander system $[\pi_{_\mathrm{H}},\pi_{_\mathrm{A}}]$.
By Theorem \ref{th:thm-intro}, it is also the generating function of our rooted and decorated triangulations of the $3$-sphere with a weight $z^{n+2}x^{N-1}$ where $n$ and $N$ are the number of tetrahedra and vertices respectively.
By brute force enumeration of the plane trees involved, we have found
\begin{align}
	M(z,x) &= 2  x^2 z^2 + (8 x + 12 x^3)z^4 + (60 x + 40 x^2) z^5 + (336 x + 996 x^2 + 420 x^3 + 618 x^4) z^6 \nonumber\\
	&\quad + (5460 x + 10416 x^2 + 6496 x^3 + 1652 x^4) z^7 \nonumber\\
	&\quad + (63344 x + 135776 x^2 + 150544 x^3 + 75360 x^4 + 46360 x^5) z^8 + \cdots. \label{eq:Mexpansion}
\end{align}
The first three terms are illustrated in Fig.~\ref{fig:tripletreesn4n5}.

Rather loose lower and upper bounds can be obtained for the radius of convergence $z_*(x)$ of $z\mapsto M(z,x)$, and thus on the exponential growth $\limsup_{n\to\infty}\frac1n\log M_n(x) = \frac{1}{z_*(x)}$ for fixed $x>0$,
\begin{align*}
	\frac{9}{2} \sqrt{x} \leq \frac{1}{z_*(x)} \leq 48 \max(1,x).
\end{align*}
The lower bound follows from an explicit family of triple trees obtained by recursive subdivisions (Section \ref{sec:specialclass}), while the upper bound simply uses that $|\mathcal{M}_n| \leq |\mathcal{H}_n|\,|\mathcal{P}_n|$.

\begin{figure}[t]
	\centering
	\includegraphics[width=\linewidth]{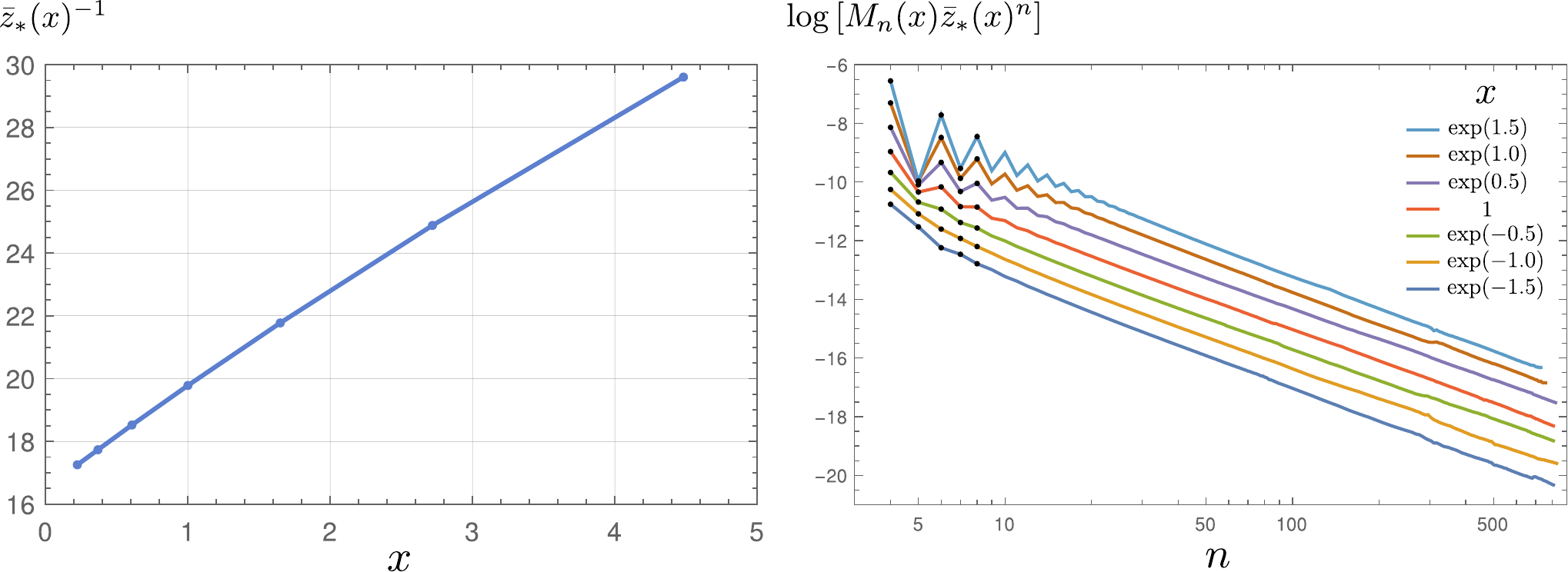}
	\caption{The left figure shows the Monte Carlo simulation estimate $\bar{z}_*(x)$ of the radius of convergence $z_*(x)$ of $M(z,x)$. In the right logarithmic plot the solid lines are the estimates for $M_n(x)$ for $7$ different values of $x$ normalized by $\bar{z}_*(x)^n$. The black dots correspond to the exact values given in \eqref{eq:Mexpansion} for small $n$, which are seen to be consistent with the simulation data. Statistical error bars are too small to display.  }
	\label{fig:simulation}
\end{figure}

Finally we provide a preview of simulation data that will be extensively discussed in an upcoming work. 
The data is based on simulations of Markov chains on triple trees $\mathcal{M}_n$ that are designed to have stationary distributions of the form $f(n) M_n(x)$ for appropriate functions $f(n)$.
From histograms of the random size $n$ one can deduce estimates of $M_n(x)$ and the radius of convergence $z_*(x)$. 
Using data up to $n\approx 1000$, these estimates are shown in Fig.~\ref{fig:simulation} for $x = 1, e^{\pm1/2}, e^{\pm1}, e^{\pm 3/2}$.
For the large values of $x$ the data is consistent with a scaling $M_n(x) z_*(x)^n \sim C\,n^{-3/2}$, as one would expect in the CRT / branched polymer universality class.
Although not clearly visible in these plots, upon decreasing $x$ the model appears to encounter a phase transition at some value $x_* \in (0,1)$, where the model starts to display characteristics that are qualitatively different from  Dynamical Triangulations. 
An extensive numerical analysis of these characteristics, including the order of the phase transition and scaling properties of the low-$x$ phase, will appear in a forthcoming paper.

\section{Origin in three-dimensional triangulations}
\label{sec:origin-in-3d}

\subsection{Complexes, triangulations}
\label{subsec:def-complexes}

A \emph{CW-complex} is a Hausdorff space with a cellular structure that is closure finite and has the weak topology relative to its cells.  The $k$-dimensional elements are called $k$-\emph{cells}. If $d=2$, the $2$-cells are called \emph{faces}. In a CW-complex, every $k$-cell $\delta_k$ is the image of a $k$-dimensional ball $b_k$ under a continuous characteristic map $\Phi[\delta_k]$ which induces an homeomorphism on the interior of $b_k$ and maps the boundary of $b_k$ to the union of lower dimensional cells contained in $\delta_k$. The CW-complexes are all assumed to be pure (every cell belongs to a $d$-cell). Two $k$-cells are allowed to share more than one $(k-1)$-cell, and a $(d-1)$-cell might belong to more than two $d$-cells. 

For every $p$-dimensional cell $\delta'_p$ of  $\delta_k$, $p<k$, the preimage by $\Phi[\delta_k]$  of the interior of  $\delta'_p$ is a collection of connected subsets of $b_k$. We say that a CW complex is \emph{non-degenerate} if  it is assumed that $\Phi[\delta_k]$ restricted to each such connected subset of $b_k$ induces a homeomorphism\footnote{This condition is to exclude any non-intuitive complexes, such as a characteristic map that sends the whole boundary of a $3$-ball to a point in the interior of a triangle and which is therefore not a $0$-cell (this is avoided by the surjectivity), or a $2$-complex with two triangles glued along two edges but so that the interior of one of the edges is glued several times to the interior of the other edge (this is avoided by the injectivity). Note that the complexes are not assumed to be regular: a tetrahedron might have two triangles identified for instance.}, and if it is also assumed that the preimage of a vertex is a collection of vertices.

In the following, by $d$-\emph{complex}, we mean a finite connected non-degenerate $d$-dimensional CW-complex,  considered up to homeomorphisms that preserve the cell structure (so that when attaching two complexes along a pair of cells, one only has to specify which sub-cell is identified with which). %

The $d$-cells of a $d$-complex are divided in two groups, called \emph{outer} and \emph{internal}. 
By $d$-dimensional \emph{triangulation}, we mean a $d$-complex for which the internal  $d$-cells and all the $k$-cells for $0\le k \le d-1$  are assumed to be simplices,
 and such that every $(d-1)$-simplex belongs to two (not necessarily distinct) $d$-cells. All triangulations/complexes are assumed to be oriented. 
The set of \emph{boundary simplices} is the set of $k$-simplices for $0\le k \le d-1$, that belong to an outer $d$-cell.  By topology of a $d$-complex $C$, we mean the topology of the complex obtained from $C$ by removing the interior of the outer $d$-cells (the space thus obtained is a CW-complex, and thereby a topological space).
By convention, a complex $C$ is said to be closed if it has no outer cells. 
A 2-dimensional triangulation is said to be planar if it has the topology of the $2$-sphere.

A \emph{tree of tetrahedra}, also called a  \emph{stacked triangulation} \cite{Kalai1987}, is a triangulation of the 3-ball whose partial dual graph, which has a vertex for each tetrahedron (but not the outer cell) and an edge between two vertices if the corresponding tetrahedra share a triangle), is a tree. Said otherwise, a tree of a single tetrahedron is the tetrahedron itself, and a tree of $n$ tetrahedra is obtained inductively from a tree of $n-1$ tetrahedra by attaching a new tetrahedron along one of the triangles of its boundary. An example is shown on the right of Fig.~\ref{fig:apollonian3D}. 

\begin{lemma}
\label{lemma:apollonian-tree-tetra}
There is a one-to-one correspondence between a tree of tetrahedra and its boundary: an Apollonian triangulation with $n+3$ vertices ($2n+2$ triangles)  is the boundary of a unique tree of $n$ tetrahedra.
\end{lemma}
\proof This is easily seen inductively: an Apollonian triangulation with 4 vertices is the boundary of a tetrahedron, and  an Apollonian triangulation with more than 4 vertices contains a vertex of valency three, which is necessarily associated to a tetrahedron with three boundary triangles. \qed

\subsection{Generalities on hierarchical triangulations}
\label{sub:hierarch}

We already encountered the notion of a hierarchical triangulation in the introduction. Here we provide a precise definition.
A \emph{hierarchical triangulation} is a planar triangulation in which for every triangle, there exists a companion triangle with which it shares three distinct vertices. By planarity, every triangle can only have a single such companion. For every hierarchical triangulation $h$, this  naturally defines a  partition $\Pi(h)$ of the triangles of $h$ in pairs (called a \emph{pairing}).

\begin{lemma}
\label{lem:reduce-triangles-hierarch}
Consider a planar triangulation $h$ that contains two triangles $A,B$  that share three distinct vertices, and denote by $z_1, z_2, z_3$ the three zones they delimitate in $h$. If they are not empty, we let $h_1, h_2, h_3$ be the triangulations respectively obtained from each zone $z_1, z_2, z_3$ by identifying the two edges they share with the triangles $A$ and $B$ in $h$. 
Then $h$ is hierarchical if and only if $h_1$, $h_2$ and $h_3$ are hierarchical, and $\Pi(h)=\Pi(h_1)\cup \Pi(h_2)\cup \Pi(h_3) \cup \{A,B\}$. 
\end{lemma}

\proof If $h_1$, $h_2$ and $h_3$ are hierarchical, then it is clear. If now  $h$ is hierarchical, every triangle has a companion in $\Pi(h)$ with which it shares three distinct vertices. If a triangle is in $z_1$, then by planarity its companion in $\Pi(h)$ must necessarily be in $z_1$ as well, and the same goes independently for each zone $z_i$, so that $h_1$, $h_2$ and $h_3$ are hierarchical and $\Pi(h)=\Pi(h_1)\cup \Pi(h_2)\cup \Pi(h_3) \cup \{A,B\}$. \qed 

\

For completeness we show that the definition agrees with the informal definition provided in the introduction.

\begin{figure}[h!]
	\centering
	\includegraphics[scale=0.7]{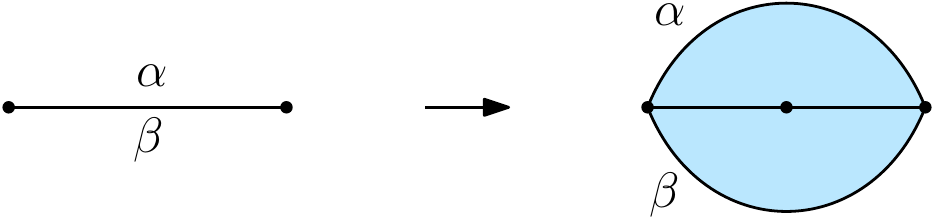}
	\caption{A hierarchical triangulation can always from  two triangles that share three edges by recursively replacing   edges by  pairs of triangles that share two edges.}
	\label{fig:Hierarch-def}
\end{figure}

\begin{lemma}
\label{lem:hierarch-recursive-construction}
A triangulation is hierarchical if and only if it can be obtained from the unique hierarchical triangulation with two triangles by recursive insertion of two triangles that share two edges as in Fig.~\ref{fig:Hierarch-def}. 
\end{lemma}

\proof Consider a planar triangulation $h$ obtained by recursive insertion of pairs as is Fig.~\ref{fig:Hierarch-def}, starting from two triangles that share three edges. Whenever a pair of triangles is inserted, they do share three vertices, and the remaining insertions of pairs of triangles in the construction do not change this property. $h$ is therefore hierarchical.

We now prove by induction on the number of triangles that a hierarchical triangulation admits such a recursive construction.  For two triangles there is nothing to prove. Let $h$ be a hierarchical triangulation with more than two triangles. Every pair of triangles $\{A,B\}$ of $\Pi(h)$ delimitates three possibly empty zones of $h$, $z_1, z_2, z_3$, and if $V(z_i)$ is the number of triangles in the zone $z_i$, we let: 
$$
w(A,B) = \min\left( V(z_1) + V(z_2), V(z_2) + V(z_3), V(z_1) + V(z_3)\right).
$$
We choose $\{A,B\}$ of $\Pi(h)$ for which $w(A,B)$ is minimal among all pairs of $\Pi(h)$, and let $z_1, z_2$ be two zones such that $w(A,B) =  V(z_1) + V(z_2)$. If $z_1$ or $z_2$ is not empty, say $z_1$, we let $A'$ be the triangle of $z_1$ adjacent to $A$, and we let $B'$ be it's companion in $\Pi(h)$. By planarity and since $A'$ and $B'$ share three distinct vertices, $B'$  must also be in $z_1$. We let $z_1', z_2'$ be the two zones delimited by $A'$ and $B'$ and which do not contain $A$ and $B$. Then, since $z_1'\cup z_2'\cup A'\cup B' \subset z_1$:
$$
w(A',B') < V(z_1') + V(z_2') + 2 \le V(z_1) \le w(A,B),
$$
which contradicts the minimality of $w(A,B)$. Both  $z_1$ and $z_2$ must therefore be empty. We can therefore perform the inverse of the operation of Fig.~\ref{fig:Hierarch-def} on $A$ and $B$ and use the induction hypothesis. \qed

\

Finally, the following simple lemma will be used later on.

\begin{lemma}\label{lem:2cycles-vs-shared}In a hierarchical triangulation,  an edge is not in any cycle of length two if and only if it is shared by two paired triangles.\end{lemma}
\proof Two paired triangles share three distinct vertices, and for each choice of two of these vertices, these two triangles have either one shared edge linking them, or two distinct edges linking them, in which case these edges form a cycle of length $2$. If an edge $e$ is not in any cycle of length $2$, consider the triangle $A$ on one side of $e$, and the triangle $B$ paired to $A$. Since $e$ cannot form a cycle of length $2$ with the edge of $B$ that links the same two vertices, then $e$ must be shared by $A$ and $B$. 
Reciprocally, if $e$ is shared by two paired triangles $A$ and $B$, then these two triangles sharing three distinct vertices $v_1, v_2, v_3$ such that $e=(v_1, v_2)$ and we call $z_1$ and $z_2$ the (possibly empty) zones respectively delimitated by the edges of $A$ and $B$ that link $v_2$ and $v_3$ on one hand, and $v_1$ and $v_3$ on the other hand. All the edges different from $e$ are either in $z_1$ or $z_2$, and by planarity the edges of $z_1$ that have $v_2$ as one endpoint cannot have $v_1$ as the other endpoint, and similarly the edges of $z_2$ cannot link $v_1$ and $v_2$. Thus $e$ is not in any cycle of length $2$. \qed

\subsection{Operations on 2-complexes}
\label{sub:operations-2-complexes}

In this section, we introduce the notions and tools necessary to enunciate and prove the first main theorem, Thm.~\ref{thm:thm-1}.

\paragraph{The gluing map ($\mathsf{Glue}$).} This map has already been defined in a particular case in the introduction, but we will need to use it in a slightly more general case. Let $t$ be an outerplanar triangulation and $\pi$ a pairing of its boundary edges (possibly with crossings). Since $t$ is oriented, we may orient its boundary edges in clockwise direction and distinguish for each boundary edge the vertex at its head and tail. 
We define $\mathsf{Glue}(t, \pi)$ as the closed triangulation obtained from $t$ by removing the outer face, and attaching for every element of $\pi$ the two paired edges with opposite orientation. The result $\mathsf{Glue}(t, \pi)$ is not planar if $\pi$ has crossings.

Letting $E_0(t, \pi)$ be the set of distinguished edges of $\mathsf{Glue}(t,\pi)$ (the former boundary edges of $t$), we record the following easy fact for future purposes.
\begin{lemma}
\label{lem:distinguished-iff-boundary}
Let $\pi$ be a pairing of the boundary edges of $t$. An edge is on the boundary of $t$ if and only if it belongs to an edge of  $E_0(t, \pi)$.
\end{lemma}

\paragraph{Embedded 2-complexes.} In the following, by \emph{embedded 2-complex} $C$, we mean a 2-complex $c$ in which all cells are simplices, together with a non-degenerate continuous attachment map $\phi$, from a (not necessarily connected) oriented 2-dimensional manifold $m$ to $c$, such that the interior of each triangle of $c$ consists of  exactly two open discs on $m$. We say that $c$ is the 2-complex \emph{underlying} $C$. 
If $m$ is a collection of 2-spheres, then $C$ can be seen as a 3-complex by taking the cones of the connected components of $m$ and extending the attachement map $\phi$ to the interiors of the resulting 3-balls to define a set of characteristic maps (recall that we see the complexes up to homeomorphisms preserving the cell structure).

Because of the non-degeneracy, \emph{an embedded 2-complex cannot have two triangles  that  share a vertex but that don't share any edge containing that vertex}: 
the preimage of that vertex under $\phi$ must be a 1-dimensional subspace of  $m$, while the non-degeneracy imposes the preimage to be a set of vertices. 

\begin{figure}[h!]
	\centering
	\includegraphics[scale=0.65]{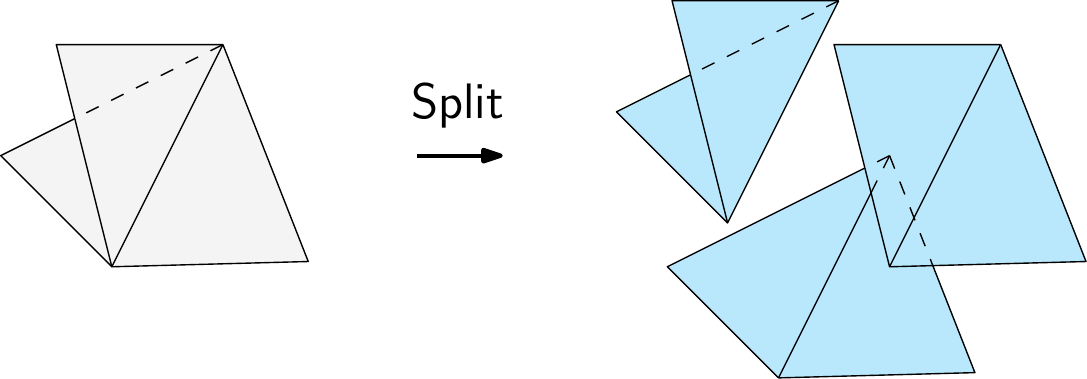}
		\caption{Given an embedded  2-complex $C$, $\mathsf{Split}(C)$ is a collection of 2-dimensional triangulations. Left: a local portion of $C$ is represented, which consists of an edge with three incident triangles and thus three corners. Right: this edge is duplicated into three edges in $\mathsf{Split}(C)$. } 
	\label{fig:Splitting-2-complex}
\end{figure}

Also because of the non-degeneracy, \emph{the connected components of $m$ inherit a cell-structure and can be seen as a collection of 2-dimensional triangulations} $\mathsf{Split}(C)$: 
If $e$  is an edge of $C$ contained in $p$ triangles, the preimage of the interior of $e$ (the edge minus its two extremities) by $\phi$ is a collection of $p$ open segments $c_1, \ldots, c_p$ on $m$. 
By continuity of $\phi$ an open curve $c_i$ has one of the preimages in $m$ of one of the triangles of $C$ containing $e$ on one side and one of the preimages of another triangle containing $e$ on the other side. Similarly, the preimage by $\phi$ of the vertices is a collection of vertices, and by continuity, the endpoints of the preimages on $m$ of the open edges of $C$ must coincide with the preimages of the vertices. 
An embedded 2-complex $C$ thus has a well-defined notion of \emph{corners}, indexed by the edges of $\mathsf{Split}(C)$: two triangles sharing an edge $e$ define a corner around $e$ if among the preimages of their interiors there is a pair that meet at one of the $c_i$ (Fig.~\ref{fig:Splitting-2-complex}).

An example has already been mentioned in the introduction: given a 3-dimensional triangulation $T$ and  a spanning tree of tetrahedra $T_0$, an example of embedded 2-complex denoted by $T^{T_0}$ is obtained by  removing the internal triangles of $T_0$, thus replacing $T_0$ by a 3-cell with the same boundary. Indeed, in doing so, we lose the internal structure of $T_0$, but only keep the information on the  attachment map from its boundary to the 2-complex underlying $T^{T_0}$ (because the complexes are considerd up to homeomorphisms preserving the cell-structure). From the discussion above:
\begin{equation}
\label{eq:T0-inverse-id}
\mathsf{Split}(T^{T_0})=\partial T_0. 
\end{equation}

\begin{figure}[h!]
	\centering
	\includegraphics[scale=0.75]{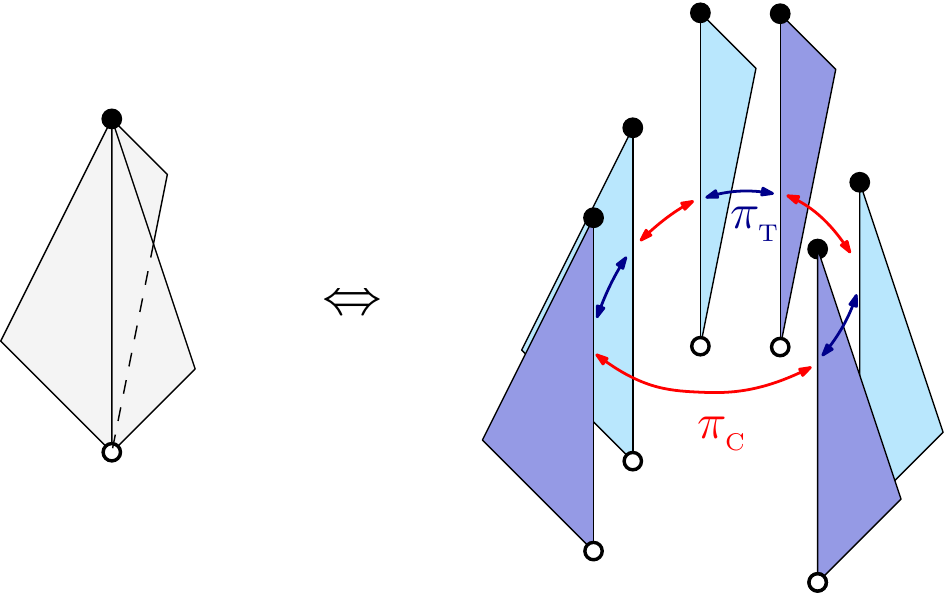}
	\caption{An embedded  2-complex $C$ can be encoded by a set of oriented triangles and two pairings of their edges. Left: a local portion of $C$ consisting of an edge and three incident triangles. Right: This edge corresponds to an orbit of $\pi_{_\mathrm{T}}\circ\pi_{_\mathrm{C}}$ (and $\pi_{_\mathrm{C}}\circ\pi_{_\mathrm{T}}$). }
	\label{fig:encoding-pairings}
\end{figure}

\paragraph{Encoding with pairings.} Given two oriented triangles (with a top side and a bottom side) we call a pairing of their edges non-crossing if the pairs are encountered in the same order when going around one of the triangle clockwise and the other counterclockwise.
\begin{lemma}
\label{lem:encoding-pairings-2comp}
An embedded 2-complex with $n$ triangles is uniquely encoded by:
\begin{enumerate}[label=-]
\item a set $\mathcal{T}$ of $2n$ oriented triangles with distinguishable edges $\mathcal{E}(\mathcal{T})$,
\item two pairings $\pi_{_\mathrm{T}}, \pi_{_\mathrm{C}}$ of these edges, 
\end{enumerate}
where two edges may only be paired by $\pi_{_\mathrm{T}}$ if they are in different triangles, and whenever two triangles have two edges paired in $\pi_{_\mathrm{T}}$, then their other four edges must also be paired in a non-crossing manner.
\end{lemma}
In the following we will sometimes view a pairing alternatively as a fixed-point free involution, meaning that we shall write
$$
\pi(e)=e' \quad \Leftrightarrow \quad \pi(e')=e \quad \Leftrightarrow \quad \{e,e'\}\in \pi. 
$$ 
\proof Consider a collection of $2n$ triangles and two pairings of their edges as in the lemma. Since the triangles are oriented, $(\mathcal{T}, \pi_{_\mathrm{C}})$ uniquely defines a collection of 2-dimensional triangulations $m$ by attaching the triangles along the paired edges $\pi_{_\mathrm{C}}$, while respecting the orientations.  
We now define the 2-complex $c$: it has one triangle for every pair of triangles of $\mathcal{T}$ whose edges are paired by $\pi_{_\mathrm{T}}$, an edge for each orbit of $\pi_{_\mathrm{T}}\circ\pi_{_\mathrm{C}}$ (Fig.~\ref{fig:encoding-pairings}), and if two triangles of $c$ share an edge, then they also share the vertices at its extremities. The pairing $\pi_{_\mathrm{T}}$ then defines an attachment  map $\phi$ from $m$ to $c$.

Reciprocally, given an embedded 2-complex $C$, we define $\mathcal{T}$ as the set of triangles of $m=\mathsf{Split}(C)$, and $\pi_{_\mathrm{C}}$ as the pairing of their edges induced by  $\mathsf{Split}(C)$. Every triangle $\theta$ of $c$ has exactly two preimages $\theta_1$ and $\theta_2$ on $m$. Let $e$ be an edge of $\theta$ in $C$ and  $U(e)$ be an open subset of the interior of $\theta$  bounded on one side by $e$. From the non-degeneracy, the preimage of $U(e)$ by $\phi$ consists of exactly two open subsets of $m$, one in $\theta_1$ and  one in $\theta_2$, each partially bounded by one of the preimages of the interior of $e$ on $m$. These two segments correspond to two edges of $\theta_1$ and $\theta_2$, which we pair. Repeating this for every pair consisting of a triangle and an edge of $C$ leads to a  pairing  $\pi_{_\mathrm{T}}$ of the edges of the triangles in $\mathcal{T}$. 
 \qed
 
 \

In the proof above, we have also shown that: 
\begin{lemma} 
\label{lem:cons-encoding-pairings-2comp}
Let $C$ be an embedded 2-complex and $\mathcal T$, $\pi_{_\mathrm{T}}, \pi_{_\mathrm{C}}$ encoding it. Then:
\begin{itemize}
\item $\mathsf{Split}(C)$ is the 2-triangulation obtained attaching the triangles of $\mathcal T$ using $ \pi_{_\mathrm{C}}$ while matching the orientations, 
\item The edges of $C$ contained in  $p$ triangles are in bijection with the orbits of $\pi_{_\mathrm{T}}\circ\pi_{_\mathrm{C}}$ (and of $\pi_{_\mathrm{C}}\circ\pi_{_\mathrm{T}}$) with $p$ elements.
\end{itemize}
\end{lemma}

The encoding of Lemma~\ref{lem:encoding-pairings-2comp} allows for a simple proof of the following:

\begin{lemma}
\label{lem:Apo-Split-gives-T}
Let $C$ be an embedded 2-complex with $n$ triangles such that $\mathsf{Split}(C)$ is Apollonian. Then there exists a unique 3-dimensional triangulation $T$ with $n-1$ tetrahedra and with a distinguished spanning tree of tetrahedra $T_0$ such that $T^{T_0}=C$. 
\end{lemma}
\proof Consider $C$ such that $\mathsf{Split}(C)$ is Apollonian with $2n$ triangles and  a 3-triangulation $T$ with a spanning tree of tetrahedra $T_0$ satisfying $T^{T_0}=C$, and consider $\pi_{_\mathrm{T}}$ from the unique encoding of $C$ of Lemma~\ref{lem:encoding-pairings-2comp}.
Then from \eqref{eq:T0-inverse-id}, $\partial T_0 = \mathsf{Split}(C)$, but from Lemma~\ref{lemma:apollonian-tree-tetra}, there exists a unique  such tree of tetrahedra $T_0$ (with $n$ triangles).  $T$ must therefore be the only 3-triangulation obtained by attaching the triangles on the boundary of  $T_0$ using the attachment map $\pi_{_\mathrm{T}}$. \qed

\ 

Given two \emph{corners} in two (not necessarily distinct) embedded 2-complexes around two oriented edges, we may \emph{attach the corresponding edges along these corners} while matching the orientations of the edges, obtaining a new embedded 2-complex. Two corners are given by two pairs $(e_1, e_2)$ and $(e_1', e_2')$ of $\pi_{_\mathrm{C}}$ in different orbits of $\pi_{_\mathrm{T}}\circ\pi_{_\mathrm{C}}$, where $e_1$ preceeds $e_2$ in clockwise direction around the edge seen from its tail (and similarly for $e_1'$ and $e_2$). 
Attaching the corresponding edges together along these corners amounts to replacing these two pairs in $\pi_{_\mathrm{C}}$  by $(e_1, e_2')$ and $(e_1', e_2)$, leading to a new $\pi'_{_\mathrm{C}}$.

Reciprocally, given an edge $e$ of an embedded 2-complex and two corners around that edge, we might detach the triangles in two groups while keeping the cyclic ordering of triangles around the two duplicates of $e$. With the same notations, the operation is exactly the same as above, with the only difference that $(e_1, e_2)$ and $(e_1', e_2')$ belong to the same orbit of $\pi_{_\mathrm{T}}\circ\pi_{_\mathrm{C}}$

 In both cases, seeing these two pairings as permutations, this amounts to composing $\pi_{_\mathrm{C}}$ with the transpositions $(e_1, e_1')$ and $(e_2, e_2')$.

\begin{figure}[h!]
	\centering
	\raisebox{1.2ex}{\includegraphics[scale=0.7]{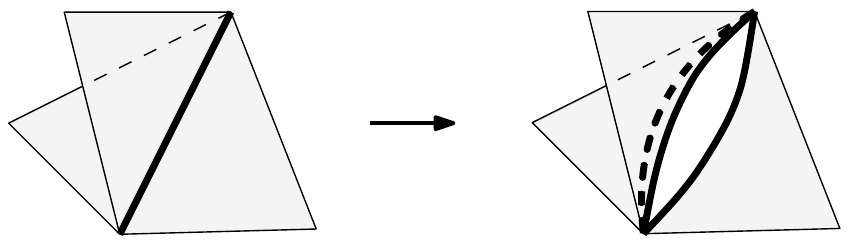}}\hspace{2cm} \includegraphics[scale=0.95]{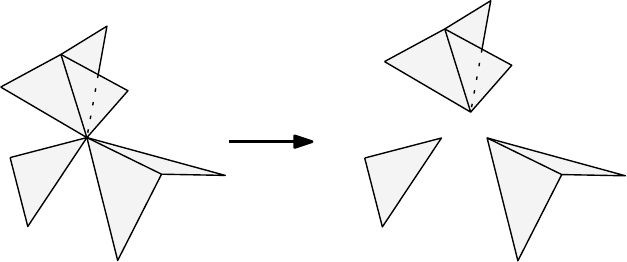}
	\caption{Steps for cutting a 2-complex along a distinguished subset of edges.   }
	\label{fig:Splitting-complement}
\end{figure}

\paragraph{Cutting a 2-complex along a subset of edges.}
Given  an embedded
2-complex $C$ and  a subset of edges $E$ of $C$, we define $C_E$  as the collection of embedded 2-complexes obtained from $C$ by detaching the triangles along the edges of $E$, and only keeping 
two triangles connected through a vertex if these two triangles share an edge containing this vertex. 
 That is, we first split open every edge $e$ of $E$  by duplicating $e$ into one edge per triangle containing it (but keeping the identifications at the vertices contained in $e$), as on the left of Fig.~\ref{fig:Splitting-complement}, and then if there are maximal components that only touch through a given vertex, we split these  components by introducing one duplicate of that vertex per component (right of Fig.~\ref{fig:Splitting-complement}). 
 We say that  $C_E$  is obtained from $C$ by \emph{cutting it along} $E$. We will denote by $\mathsf{Cut}(E)$ the set of edges of $C_E$  obtained by duplicating the edges of $E$. 
 
We call \emph{free edge} of a complex an edge which is contained in a single triangle, and in no higher dimensional simplex. A free edge in an embedded 2-complex $C$ encoded by $\mathcal T$, $\pi_{_\mathrm{T}}, \pi_{_\mathrm{C}}$ corresponds to an orbit of $\pi_{_\mathrm{T}}\circ\pi_{_\mathrm{C}}$ with one element $e$: the two edges of triangles of $\mathcal T$,  $e$ and $\pi_{_\mathrm{C}}(e)$, are paired both in $\pi_{_\mathrm{T}}$ and in  $\pi_{_\mathrm{C}}$. Replacing an edge of $C$ contained in $p$ triangles by $p$ free edges amounts to replacing $\pi_{_\mathrm{C}}$ by $\pi_{_\mathrm{T}}$ on the subset of $\mathcal{E}(\mathcal T)$ that are  in the orbits of $\pi_{_\mathrm{T}}\circ\pi_{_\mathrm{C}}$ and of $\pi_{_\mathrm{C}}\circ\pi_{_\mathrm{T}}$.  Reciprocally, attaching $p$ free edges to form an edge with $p$ incident triangles amounts to replacing, for the corresponding elements of $\mathcal{E}(\mathcal T)$, the pairing $\pi_{_\mathrm{C}}$ - which on this subset coincides with $\pi_{_\mathrm{T}}$ - by another pairing.

\paragraph{Trees of triangles.}A \emph{tree of triangles} (left of Fig.~\ref{fig:apollonian3D}) is an embedded 2-complex  whose dual graph (with a disc vertex for each triangle, a square vertex for each edge, and an edge between a disc vertex and a square vertex if the corresponding triangle contains the corresponding edge) is a tree.

 A tree of one triangle is a triangle with three distinct edges and vertices. A tree of $n$ triangles can be recursively obtained from a tree of $n-1$ triangles by choosing a corner of an edge $e$ of that tree, and attaching a new triangle by identifying one of its edges with the edge $e$ along that corner.

A tree of triangles is said to be \emph{rooted} if it has a distinguished oriented corner, that is, a distinguished triangle $A$ around an oriented edge $e$ (if $B$ is the triangle 
 preceding $A$ around the edge $e$ in clockwise order when we look at $e$ from its origin, then the oriented corner can be chosen to be $B \rightarrow A$).

\begin{lemma}
\label{lem:hierarch-iff-triang-tree}
There is a bijection $\mathsf{Id}$ between rooted hierarchical triangulations $h$ with $2n$ triangles and rooted trees of $n$ triangles, that maps edges of $h$  to the corners of $\mathsf{Id}[h]$. 

In particular, the vertices of $h$ are mapped to the vertices of $\mathsf{Id}[h]$, and an edge of $\mathsf{Id}[h]$ incident to $k\ge 1$ triangles corresponds to $k$ edges of $h$ that link the same two vertices. 
\end{lemma}

\proof This is obvious for $n=2$ and we prove it by induction on $n$. Let $h$ be a rooted (oriented) hierarchical triangulation with $2n$ triangles, $A$ the triangle on the right of the root edge, and $B$ its companion triangle in $\Pi(h)$. We  attach the two triangles $A$ and $B$ on their top sides  to form a single triangle $A'$, so that the edges of $A$ and $B$ that share the same  two vertices in $h$ are identified and so that the three vertices shared by $A$ and $B$ are still three distinct vertices of $A'$. 
The resulting space is denoted by $C$. 

We now label the root edge $e_1$, and orient the two other edges of $A'$ to form an oriented cycle, and label them $e_2$ and $e_3$ in increasing order along that cycle. From Lemma~\ref{lem:reduce-triangles-hierarch},  $C$ has the following form: attached to each edge $e_i$ is a hierarchical triangulation $h_i$, rooted on the edge $e_i$. 
We may apply the induction hypothesis for $h_1$, $h_2$ and $h_3$: they are in bijection with three rooted trees of triangles $C_1$, $C_2$ and $C_3$. From $C$ we keep only the triangle $A'$, and for each $1\le i\le 3$, we attach the edge $e_i$ of this single triangle $A'$ with the edge along the root corner of $C_i$, while matching the orientations of the two edges. The resulting space $C'$ is a tree of triangle, which is rooted at the oriented edge  $e_1$ with the distinguished triangle $A'$. This mapping is clearly invertible and thus defines a bijection.  

The statement regarding the vertices is clear by induction. The statement regarding the edges of $\mathsf{Id}[h]$ is true by induction on all edges except for the edges  $e_1, e_2, e_3$, which have one additional incident triangle. If in $h_i$, the preimage by $\mathsf{Id}$ of the edge $e_i$ consists of $k_i$ edges linking the same two vertices, then in $h$ these same two vertices are linked by $k_i+1$ edges (to obtain $h_i$ from $h$, one removes the triangles $A$ and $B$ and identifies the edges of $A$ and $B$ that share the same two vertices, see Lemma~\ref{lem:reduce-triangles-hierarch}). 
 \qed

\medskip

\begin{lemma}
\label{lem:tree-in-h-or-IDh}
Let $h$ be a hierarchical triangulation, $E_0$ a subset of edges of $h$ that do not belong to any cycle of length $2$, and $E$ the subset of free edges of $\mathsf{Id}[h]$ which is the image of $E_0$ by  $\mathsf{Id}$. Then $E_0$ is a spanning tree of edges of $h$ 
 if and only if $E$ is a spanning tree of 
 edges  of $\mathsf{Id}[h]$.  
\end{lemma}

\proof There is the same number of elements in $E_0$  and $E$, and two edges of $E_0$ share a vertex if and only if the same is true for their image by $\mathsf{Id}$.  \qed

\paragraph{The identification map  ($\mathsf{Id}$).}
The map $\mathsf{Id}$ of Lemma~\ref{lem:hierarch-iff-triang-tree} can therefore be seen as the inverse of $\mathsf{Split}$\footnote{The difference between the attachment map $\phi$ and $\mathsf{Id}$ is that $\phi$ is a map from the set of (not necessarily connected) 2-triangulations to \emph{the 2-complex} $c$ underlying $C$, while $\mathsf{Id}$ is a map from the same domain but to \emph{the space of embedded 2-complexes}. For an embedded 2-complex $C$ for which $\mathsf{Split}(C)$ is a collection of planar triangulations, $C$  can be seen as a 3-complex, and then rather than an attachment map,  $\mathsf{Id}$ can be seen as the map induced by the union of the characteristic maps of the 3-cells on the boundaries of their preimage.
} restricted to the set of trees of triangles:
\begin{align*}
	\mathsf{Id} = \mathsf{Split}^{-1}.
\end{align*}
Since $\mathsf{Split}$ is actually defined on the larger set of embedded 2-complexes, we may extend the definition domain of $\mathsf{Id}$ to the set of collections of non-necessarily planar 2-dimensional triangulations, if we provide a prescription on which triangles must be attached one to another, and for each such pair, which edges should be glued to which (given a collection of triangulations, from Lemmas \ref{lem:encoding-pairings-2comp} and \ref{lem:cons-encoding-pairings-2comp}, $\pi_{_\mathrm{C}}$ is known and we only need to provide a $\pi_{_\mathrm{T}}$ to obtain an embedded 2-complex).

If $t$ is an outerplanar triangulation and $\pi_{_\mathrm{H}}$ is a non-crossing pairing of its boundary edges such that $\mathsf{Glue}(t,\pi_{_\mathrm{H}})$ is hierarchical, then for any other pairing $\pi$ of the boundary edges of $t$ (not necessarily non-crossing) we can use the data provided by $\pi_{_\mathrm{H}}$ to build an embedded 2-complex out of $\mathsf{Glue}(t, \pi)$ (that is, define the corresponding $\pi_{_\mathrm{T}}$ of Lemma \ref{lem:encoding-pairings-2comp}): we can attach the triangles of $\mathsf{Glue}(t, \pi)$ if they are paired in $\mathsf{Glue}(t,\pi_{_\mathrm{H}})$, and identify the edges of  two paired triangles of $\mathsf{Glue}(t, \pi)$ if these edges share the same two vertices in the hierarchical triangulation $\mathsf{Glue}(t,\pi_{_\mathrm{H}})$. This defines an embedded 2-complex  
$$
\mathsf{Id}_{\pi_{_\mathrm{H}}}\bigl[\mathsf{Glue}(t, \pi)\bigr],
$$ 
so that the map $\mathsf{Id}$ of Lemma \ref{lem:hierarch-iff-triang-tree} is a particular case of $\mathsf{Id}_{\pi_{_\mathrm{H}}}$:
$$
\mathsf{Id}_{\pi_{_\mathrm{H}}}\bigl[\mathsf{Glue}(t, {\pi_{_\mathrm{H}}})\bigr] = \mathsf{Id}\bigl[\mathsf{Glue}(t, {\pi_{_\mathrm{H}}})\bigr].
$$
 From the construction,
\begin{equation}
\label{eq:split-inverse-id}
\mathsf{Split}\left(\mathsf{Id}_{\pi_{_\mathrm{H}}}\bigl[\mathsf{Glue}(t, \pi)\bigr]\right) = \mathsf{Glue}(t,\pi).
\end{equation}

We recall that $E_0(t, \pi)$ denotes the set of distinguished edges of $\mathsf{Glue}(t, \pi)$ resulting from the identification of the boundary edges of $t$. 
\begin{lemma}
\label{lem:distinguished-edges-on-Id-Glue}
With these notations and considering $\mathcal T$, $\pi_{_\mathrm{T}}, \pi_{_\mathrm{C}}$  encoding $\mathsf{Id}_{\pi_{_\mathrm{H}}}\bigl[\mathsf{Glue}(t, \pi)\bigr]$, 
let $e\in \mathcal{E}(\mathcal{T})$ be an edge of a triangle of $\mathcal{T}$ and $e' = \pi_{_\mathrm{T}}(e) \in \mathcal{E}(\mathcal{T})$ . Then the following are equivalent: 
\begin{enumerate}[label=(\roman*)]
\item $e$ belongs to an edge of $E_0(t, \pi)$,
\item $e$ and $e'$ belong to an edge of $E_0(t, \pi)$,
\item $e$ and $e'$ are on the boundary of $t$ and they are paired in $\pi_{_\mathrm{H}}$. 
\end{enumerate}
\end{lemma}
\proof Let us assume that $e$ belongs to an edge of $E_0(t, \pi)$. From Lemma~\ref{lem:distinguished-iff-boundary}, it belongs to an edge $e_{_\mathrm{H}}$ of $E_0(t, \pi_{_\mathrm{H}})$ on $\mathsf{Glue}(t, \pi_{_\mathrm{H}})$. From Lemma~\ref{lem:2cycles-vs-shared}, since $e_{_\mathrm{H}}$ is not in any cycle of length $2$, it must be shared by two paired triangles. By definition, if $\theta$ is the triangle of $\mathcal{T}$ containing $e$, then $e'$ is the edge of $\mathcal{E}(\mathcal{T})$ that belongs to the triangle $\theta'$ of $\mathcal{T}$ paired with $\theta$ in $\Pi(\mathsf{Glue}(t, \pi_{_\mathrm{H}}))$, and with which $e$ shares two vertices on $ \mathsf{Glue}(t, \pi_{_\mathrm{H}})$. Therefore, $e'$ must be the other edge of $\mathcal{E}(\mathcal{T})$ to which $e$ is attached to form $e_{_\mathrm{H}}$, so that in particular $e'$ belongs to an edge of  $E_0(t, \pi_{_\mathrm{H}})$, and from Lemma~\ref{lem:distinguished-iff-boundary}, $e'$ belongs to an edge of  $E_0(t, \pi)$, so that \emph{(i)}$\Rightarrow$\emph{(ii)}, and therefore \emph{(i)}$\Leftrightarrow$\emph{(ii)}. 
But we have actually shown that \emph{(i)} is equivalent to $e$ and $e'$ belonging to the same edge of $E_0(t, \pi_{_\mathrm{H}})$, or said otherwise, \emph{(iii)}.\qed

\ 

Let us comment on Lemma~\ref{lem:distinguished-edges-on-Id-Glue}: from the point \emph{(ii)} of the Lemma, we know that if $e$ is an edge of the embedded 2-complex $C=\mathsf{Id}_{\pi_{_\mathrm{H}}}[\mathsf{Glue}(t, \pi)\bigr)]$, then either \emph{all} of its corners on  $\mathsf{Split}(C) = \mathsf{Glue}(t,\pi)$ are distinguished edges in $E_0(t, \pi)$, or \emph{none} of them are in $E_0(t, \pi)$.

\paragraph{The spanning tree $E$.} We may therefore define $E(t, \pi_{_\mathrm{H}},\pi)$ as the set of edges of $\mathsf{Id}_{\pi_{_\mathrm{H}}}[\mathsf{Glue}(t, \pi)]$ whose corners are in $E_0(t, \pi)$ on $\mathsf{Glue}(t,\pi)$.
\begin{lemma}
\label{lem:pi-reduce-to-piH}
Let $t\in\mathcal{O}_{n+1}$ and $\pi_{_\mathrm{H}}\in\mathcal{P}_{n+1}$ such that $\mathsf{Glue}(t,\pi_{_\mathrm{H}})\in\mathcal{H}_{n+1}$, and let 
$\pi$ be a pairing of the boundary edges of $t$.  
Take $\mathcal T$,  
$\pi_{_\mathrm{T}}, \pi_{_\mathrm{C}}$ encoding $C=\mathsf{Id}_{\pi_{_\mathrm{H}}}[\mathsf{Glue}(t, \pi)\bigr)]$, and consider the subset $\mathcal{E}_0$  of $\mathcal{E}(\mathcal{T})$ of edges of triangles in $\mathcal{T}$ which belong to an edge in $E(t, \pi_{_\mathrm{H}},\pi)$. Then, restricted to this subset $\mathcal{E}_0$:
 $${\pi_{_\mathrm{T}}} _{\lvert_{\mathcal{E}_0}} = \pi_{_\mathrm{H}},\qquad\mathrm{ and } \qquad \pi_{_\mathrm{C}}{_{\lvert_{\mathcal{E}_0}}} =\pi.$$
\end{lemma}
\proof From Lemma~\ref{lem:cons-encoding-pairings-2comp}, attaching the elements of $\mathcal T$ according to the pairing $\pi_{_\mathrm{C}}$ while respecting the orientation leads to $\mathsf{Split}(C)$.  For the edges of $E_0(t, \pi)$ on $\mathsf{Split}(C) = \mathsf{Glue}(t,\pi)$, this pairing coincides with $\pi$. The statement regarding  $\pi_{_\mathrm{T}}$ is just point \emph{(iii)} of Lemma~\ref{lem:distinguished-edges-on-Id-Glue}. \qed

\

\begin{figure}[h!]
	\centering
	\includegraphics[scale=0.7]{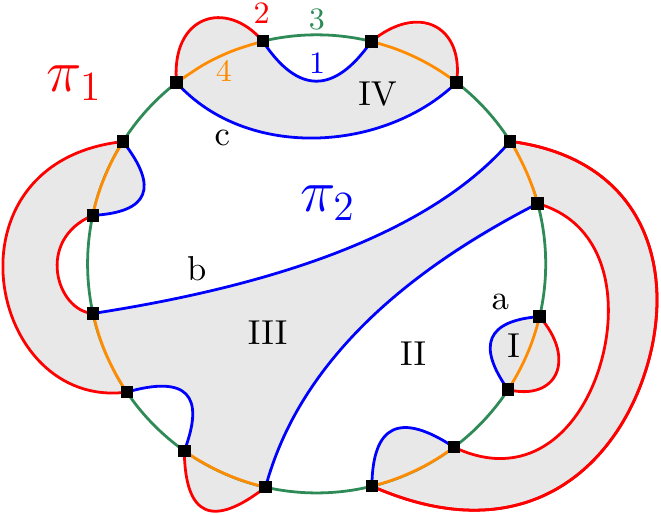}\hspace{1cm} \includegraphics[scale=0.8]{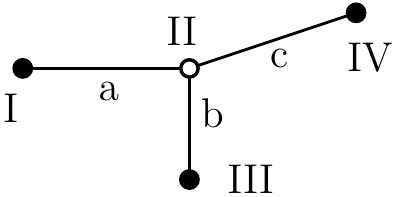}
	\caption{Zones and colors on a meander system. The four zones are labeled with Roman numbers and the three loops with letters. The corresponding adjacency tree $\mathcal{G}(\pi_1, \pi_2)$ is shown on the right.}
	\label{fig:ex-meander-zones}
\end{figure}

Recall the definition of a meander system from the introduction. 
Given two non-crossing pairings $\pi_1, \pi_2$ of $\{1,\ldots, 2n\}$ for $n\ge 1$, the \emph{meander system} $[\pi_1, \pi_2]$ is the planar map obtained from a cycle $\gamma_n$ with $2n$ vertices by drawing arcs between pairs of vertices according to $\pi_1$ on the inside and $\pi_2$ on the outside. The \emph{loops of the meander system} $[\pi_1, \pi_2]$ are the connected components of the map after dropping all edges of $\gamma_n$.
The loops of a meander system partition $[\pi_1, \pi_2]$ into \emph{zones}, which are regions of the map delimited by the loops and that only have edges of $\gamma_n$ in their interior. 
To a meander system $[\pi_1, \pi_2]$, we associate the adjacency tree $\mathcal{G}(\pi_1, \pi_2)$ of the zones. It has a vertex for each zone of $[\pi_1, \pi_2]$ and an edge connecting the two adjacent zones for each loop of $[\pi_1, \pi_2]$.
Since the loops are simple and disjoint, the resulting graph is a tree (see Fig.~\ref{fig:ex-meander-zones}).

\begin{proposition}
\label{prop:E-spanning-tree-meanders}
With the notations of the previous lemma, if in addition $\pi$ is non-crossing, then there is a bijection between $E(t, \pi_{_\mathrm{H}},\pi)$ and the tree $\mathcal{G}( \pi_{_\mathrm{H}},\pi)$ associated to the meander system $[\pi, \pi_{_\mathrm{H}}]$, which respectively maps the edges and the vertices of $E(t, \pi_{_\mathrm{H}},\pi)$  to the loops and zones of the meander system $[\pi, \pi_{_\mathrm{H}}]$. 

In particular, $E(t, \pi_{_\mathrm{H}},\pi)$ is a spanning tree of $\mathsf{Id}_{\pi_{_\mathrm{H}}}[\mathsf{Glue}(t, \pi)]$. 
\end{proposition}

\proof Let $\mathcal{T}$, $\pi_{_\mathrm{T}}, \pi_{_\mathrm{C}}$ encode $\mathsf{Id}_{\pi_{_\mathrm{H}}}[\mathsf{Glue}(t, \pi)]$ (Lemma~\ref{lem:encoding-pairings-2comp}). From Lemma~\ref{lem:cons-encoding-pairings-2comp}, the edges in $E(t, \pi_{_\mathrm{H}},\pi)$  correspond to the orbits of $\pi_{_\mathrm{T}}\circ\pi_{_\mathrm{C}}$ (and of $\pi_{_\mathrm{C}}\circ\pi_{_\mathrm{T}}$) restricted to the subset of $\mathcal{E}(\mathcal{T})$ of edges of triangles in $\mathcal{T}$ which belong to an edge in $E(t, \pi_{_\mathrm{H}},\pi)$. From  Lemma~\ref{lem:pi-reduce-to-piH}, these correspond to the orbits of $\pi_{_\mathrm{H}}\circ\pi$ (and of $\pi\circ\pi_{_\mathrm{H}}$), which in turn correspond to the loops of the meander system $[\pi, \pi_{_\mathrm{H}}]$.

Let us color the vertices of $t$ in black and white, so that the origin of the root edge is black and every edge on the boundary of $t$ links a black and a white vertex. If $\pi$ is a non-crossing pairing of the boundary edges of $t$, then the gluing of the boundary edges according to $\pi$ respects the vertex colors. An edge of $E(t, \pi_{_\mathrm{H}},\pi)$  consists of boundary edges of $t$ that are glued together using both $\pi$ and $\pi_{_\mathrm{H}}$, and in both cases this is done with the same convention as for $\mathsf{Glue}(t, \cdot)$, so that the white vertices are glued together and so are the black vertices, as shown in Fig.~\ref{fig:encoding-pairings}. 

Two edges $\epsilon$, $\epsilon'$  of $E(t, \pi_{_\mathrm{H}},  \pi)$ share a \emph{white (resp.~black)} vertex if and only if two boundary edges $e, e'$ of $t$ respectively  belonging to $\epsilon$ and a $\epsilon'$ are linked by a sequence $e_1, \ldots, e_k$, $k\geq 2$ of boundary edges of $t$ such that $e_1=e$, $e_k=e'$, and two consecutive edges $e_i, e_{i+1}$ in the sequence satisfy one of the following:

-  $e_i, e_{i+1}$ are paired in $\pi$,

- $e_i, e_{i+1}$ are paired in $\pi_{_\mathrm{H}}$,

- $e_i, e_{i+1}$ share a \emph{white  (resp.~black)}  vertex on the boundary of $t$. 
 
We may introduce two additional pairings $\pi_\circ$ and $\pi_{\bullet}$ that pair two boundary edges of $t$ if they share a white resp.~black vertex respectively, and then seeing all the pairings as maps which to an edge associate the other edge in the pair, the \emph{white (resp.~black) vertices of $E(t, \pi_{_\mathrm{H}},  \pi)$ correspond to the transitivity classes of boundary edges of $t$ under the action of $\pi$, $\pi_{_\mathrm{H}}$ and $\pi_\circ$ (resp.~$\pi_{\bullet}$)}. 

 \begin{figure}[h!]
	\centering
	\includegraphics[scale=0.6]{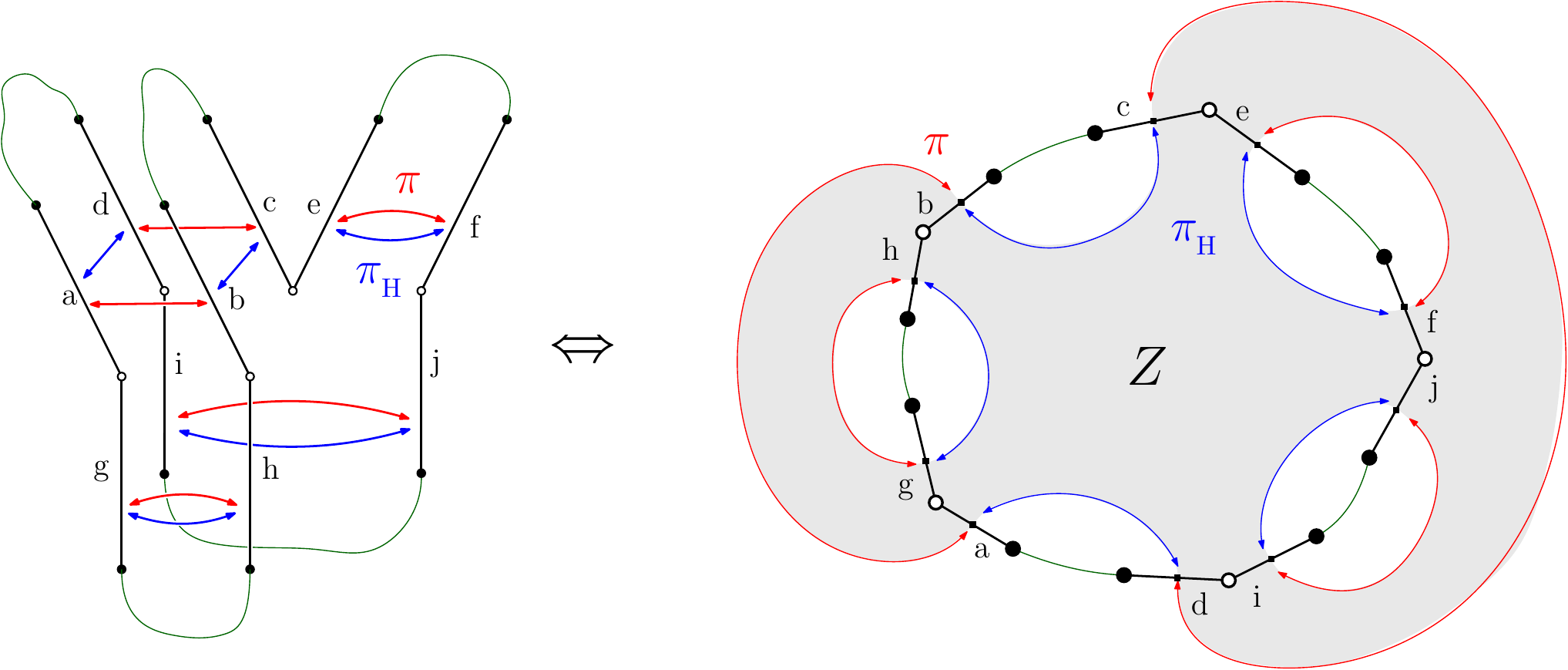}
	\caption{Left: The boundary of $t$ is schematized, and we have represented all the boundary edges of $t$ that belong to the edges incident to a given white vertex of $E(t, \pi_{_\mathrm{H}},\pi)$ of valency 4. This corresponds to a transitivity class under the action of $\pi$, $\pi_{_\mathrm{H}}$ and $\pi_\circ$. Right: This class is in bijection with a zone $Z$ of the meander system $[\pi, \pi_{_\mathrm{H}}]$ bounded by 4 loops.}
	\label{fig:edges-of-E}
\end{figure}

Note that the vertices of the meander system naturally correspond to the boundary edges of $t$ and that the vertices of $t$ included in a zone of $[\pi, \pi_{_\mathrm{H}}]$ are either all black or all white.
Given an all-white zone $Z$ and one of its incident vertices in $[\pi, \pi_{_\mathrm{H}}]$ corresponding to a boundary edge $e$ of $t$. 
Then the orbit of $e$ under $\pi$ and $\pi_{_\mathrm{H}}$ is the loop of the meander system containing $e$. 
The image under $\pi_\circ$ of any edge on this loop will again be an edge adjacent to $Z$, potentially on another loop delimiting $Z$.
Hence, the orbit of $e$ under $\pi$, $\pi_{_\mathrm{H}}$ and $\pi_\circ$ is equal to the set of vertices of $[\pi, \pi_{_\mathrm{H}}]$ that are incident to $Z$.
Hence, the transitivity classes $\pi$, $\pi_{_\mathrm{H}}$ and $\pi_\circ$ are in bijection with the white zones of the meander system $[\pi, \pi_{_\mathrm{H}}]$. 
Analogously the transitivity classes $\pi$, $\pi_{_\mathrm{H}}$ and $\pi_\bullet$ are in bijection with the black zones of the meander system $[\pi, \pi_{_\mathrm{H}}]$. 
The vertices, respectively edges, of $E(t, \pi_{_\mathrm{H}},  \pi)$ and the adjacency graph $\mathcal{G}( \pi_{_\mathrm{H}},\pi)$ of the zones are thus in one-to-one correspondence.

It only remains to prove that $E(t, \pi_{_\mathrm{H}}, \pi)$  spans  $\mathsf{Id}_{\pi_{_\mathrm{H}}}[\mathsf{Glue}(t, \pi)]$: $E_0(t, \pi)$ is a spanning tree of $\mathsf{Glue}(t, \pi)$. In order to recover the embedded 2-complex $\mathsf{Id}_{\pi_{_\mathrm{H}}}\bigl[\mathsf{Glue}(t, \pi)\bigr]$, one has to add the pairing $\pi_{_\mathrm{T}}$, to attach the triangles of $\mathsf{Glue}(t, \pi)$ two-by-two. The image $E(t, \pi_{_\mathrm{H}},\pi)$  of $E_0(t, \pi)$ by this operation is still a spanning subset of edges. \qed

\subsection{Triple-trees encode three-dimensional triangulations
}
%
We recall that the embedded 2-complex $T^{T_0}_E$ is obtained from $T$ by removing the internal triangles of $T_0$, and then cutting along the edges of $E$,  and that $T$ is said to be rooted if it has a marked oriented edge, and a marked triangle not in $T_0$ (a triangle of $T^{T_0}$), among the triangles containing this edge.

\begin{theorem}
\label{thm:thm-1}
Let $\mathbb{T}_n$ be the set of 3-dimensional triangulations $T$ with $n$ tetrahedra, with a marked spanning tree of tetrahedra $T_0$ and a marked spanning tree of edges $E$, such that the 2-complex $T^{T_0}_E$ is a tree of triangles and $\mathsf{Cut}(E)$ is a spanning tree of $T^{T_0}_E$. $T$ is rooted on an edge of $E$.

For $n\ge2$, there is a bijection between $\mathcal{M}_{n+1}$ and $\mathbb{T}_{n-1}$, which  to a triple-tree $(t,\pi_{_\mathrm{H}}, \pi_{_\mathrm{A}})$ assigns a unique 3-dimensional triangulation $T$ with a distinguished spanning tree of tetrahedra $T_0$ satisfying
\begin{equation}
\label{eq:bijection}
T^{T_0} = \mathsf{Id}_{\pi_{_\mathrm{H}}}\bigl[\mathsf{Glue}(t,\pi_{_\mathrm{A}})\bigr].
\end{equation}
The distinguished spanning tree of edges of $T$ is $E=E(t, \pi_{_\mathrm{H}},\pi_{_\mathrm{A}})$, in bijection with  the tree $\mathcal{G}( \pi_{_\mathrm{H}},\pi_{_\mathrm{A}})$. In particular,   the number of vertices of the triangulation $T$ minus one is the number of loops of the meander system $[\pi_{_\mathrm{H}}, \pi_{_\mathrm{A}}]$. 
\end{theorem}
We let $\mathbb{T}=\cup_n \mathbb{T}_n$. In the following, an element of $\mathbb{T}$ will therefore be called a \emph{triple-tree triangulation}. 
\proof Consider a triple-tree $(t,\pi_{_\mathrm{H}}, \pi_{_\mathrm{A}})\in \mathcal{M}_{n+1}$. Since $\mathsf{Glue}(t,\pi_{_\mathrm{A}})$ is Apollonian with $2n$ triangles, we know from Lemma~\ref{lem:Apo-Split-gives-T} that there is a unique 3-triangulation $T$ with $n-1$ tetrahedra and with a marked spanning tree of tetrahedra $T_0$ such that \eqref{eq:bijection}. 
If we set $E=E(t, \pi_{_\mathrm{H}},\pi_{_\mathrm{A}})$, we know from Prop.~\ref{prop:E-spanning-tree-meanders} that $E$ is a spanning tree of edges of $T^{T_0}$ and therefore of $T$. 
Since $\mathsf{Glue}(t,\pi_{_\mathrm{A}})=\mathsf{Split}(T^{T_0})$ has an oriented root edge in $E_0(t, \pi_{_\mathrm{A}})$, $T^{T_0}$ has an oriented root edge in $E$, and a marked corner around it. 
From Prop.~\ref{prop:E-spanning-tree-meanders} again, we know that $E(t, \pi_{_\mathrm{H}},\pi_{_\mathrm{A}})$ is in bijection with the tree $\mathcal{G}( \pi_{_\mathrm{H}},\pi_{_\mathrm{A}})$. The number of vertices of $T$ correspond to the number of zones of the meander system $[\pi_{_\mathrm{H}}, \pi_{_\mathrm{A}}]$, which is its number of loops plus one. 
Finally, as explained in the paragraph where $C_E$ is introduced, going from $T^{T_0}$ to $T^{T_0}_E$  amounts to replacing  for the set $\mathcal{E}_0$ of elements of $\mathcal{E}(\mathcal{T})$ that belong to an edge of $E$, the pairing ${\pi_{_\mathrm{C}}}_{\lvert_{\mathcal{E}_0}}$ by the pairing ${\pi_{_\mathrm{T}}}_{\lvert_{\mathcal{E}_0}}$, that is, from Lemma~\ref{lem:pi-reduce-to-piH}, by replacing $\pi_{_\mathrm{A}}$ by $\pi_{_\mathrm{H}}$. The resulting embedded 2-complex is 
\begin{equation}
\label{eq:TT0E}
T^{T_0}_{E}=\mathsf{Id}\bigl[\mathsf{Glue}(t,\pi_{_\mathrm{H}})\bigr], \quad \mathrm{and} \quad \mathsf{Cut}(E)=E(t, \pi_{_\mathrm{H}},\pi_{_\mathrm{H}}).
\end{equation}
$T^{T_0}_{E}$ is therefore a tree of $n$ triangles, and from Lemma~\ref{lem:tree-in-h-or-IDh} or Prop.~\ref{prop:E-spanning-tree-meanders}, $\mathsf{Cut}(E)$ is a spanning tree of $T^{T_0}_E$. 

\ 

Reciprocally, consider a 3-dimensional triangulation $T$ with a marked spanning tree of tetrahedra $T_0$ and a marked spanning tree of edges $E$, such that the 2-complex $T^{T_0}_E$ is a tree of triangles and $\mathsf{Cut}(E)$ is a spanning tree of $T^{T_0}_E$. Lemma~\ref{lem:hierarch-iff-triang-tree} implies the existence of a hierarchical triangulation $h$ with $2n$ triangles such that $T^{T_0}_E = \mathsf{Id}[h]$.  From Lemma~\ref{lem:tree-in-h-or-IDh}, since $\mathsf{Cut}(E)$ is a spanning tree of free edges of $\mathsf{Id}[h]$,  the preimage of $\mathsf{Cut}(E)$ by $\mathsf{Id}$ is a distinguished spanning tree of edges that do not belong to any cycle of length $2$. The marked oriented edge and triangle of $T^{T_0}$ (from the rooting of $T$) single out an oriented free edge of $T^{T_0}_E$, so $h$ is rooted on its distinguished spanning tree, and therefore $h\in\mathcal{H}_{n+1}$. Consequently, there exists $t\in\mathcal{O}_{n+1}$ and $\pi_{_\mathrm{H}}\in\mathcal{P}_{n+1}$ such that $h=\mathsf{Glue}(t,\pi_{_\mathrm{H}})$,  and the preimage of $\mathsf{Cut}(E)$ by $\mathsf{Id}$ is $E_0(t, \pi_{_\mathrm{H}})$, so that $\mathsf{Cut}(E)=E(t, \pi_{_\mathrm{H}},\pi_{_\mathrm{H}})$, that is, \eqref{eq:TT0E}.

Consider $\mathcal T$, $\pi_{_\mathrm{T}}, \pi_{_\mathrm{C}}$ encoding $T^{T_0}_E$ (Lemma~\ref{lem:encoding-pairings-2comp}). The embedded 2-complex $T^{T_0}$ is recovered from $T^{T_0}_E$ by  attaching the elements of $\mathsf{Cut}(E)$ (which are free edges) in groups with a prescribed cyclic ordering of the edges in each group.  
As discussed above in the paragraph where $C_E$ is introduced, this amounts to replacing, for the set $\mathcal{E}_0$ of elements of $\mathcal{E}(\mathcal{T})$ that belong to an edge of $\mathsf{Cut}(E)$, the pairing ${\pi_{_\mathrm{C}}}_{\lvert_{\mathcal{E}_0}}$ by another pairing $\pi_{_\mathrm{A}}$. But since  $\mathsf{Cut}(E)=E(t, \pi_{_\mathrm{H}},\pi_{_\mathrm{H}})$, from Lemma~\ref{lem:pi-reduce-to-piH}, this amount to replacing 
$$
{\pi_{_\mathrm{C}}}_{\lvert_{\mathcal{E}_0}} = \pi_{_\mathrm{H}}
$$
by another pairing $\pi_{_\mathrm{A}}$, so that from Lemma~\ref{lem:cons-encoding-pairings-2comp}, $\mathsf{Split}(T^{T_0})$ is the 2-triangulation obtained by attaching the triangles on the boundary of $t$ using $\pi_{_\mathrm{A}}$ instead of $\pi_{_\mathrm{H}}$ while $\pi_{_\mathrm{T}}$ is left unchanged, that is:
\begin{equation}
\label{eq:TT0}
T^{T_0}=\mathsf{Id}_{\pi_{_\mathrm{H}}}\bigl[\mathsf{Glue}(t,\pi_{_\mathrm{A}})\bigr],  \quad \mathrm{and} \quad E=E(t, \pi_{_\mathrm{H}},\pi_{_\mathrm{A}}).
\end{equation}

Applying $\mathsf{Split}$ to both sides of  the leftmost equation of \eqref{eq:TT0}, we see from \eqref{eq:split-inverse-id} and from \eqref{eq:T0-inverse-id} that  
\begin{equation}
\partial T_0 = \mathsf{Glue}(t,\pi_{_\mathrm{A}}), 
\end{equation}
which is planar, so that $\pi_{_\mathrm{A}}$ must be a non-crossing pairing satisfying 
$
\mathsf{Glue}(t,\pi_{_\mathrm{A}}) \in \mathcal{A}_{n+1}.
$
Therefore,  $(t,\pi_{_\mathrm{H}}, \pi_{_\mathrm{A}})$ is a triple-tree, which concludes the proof. 
\qed

\

At this stage, in order to complete the proof of Thm.~\ref{th:thm-intro} it remains to show that a triangulation satisfying these conditions has the topology of $S^3$.

\section{Triple-tree triangulations and local constructibility}

In this subsection, we define locally constructible triangulations, originally introduced in \cite{Durhuus-Jonsson}, and show that $\mathbb{T}$ is a subset of the locally constructible triangulations that we will explicitly characterize.

\subsection{Triple-tree triangulations are locally constructible}
\label{sub:triple-trees-are-LC}

\subsubsection{Locally constructible triangulations}

\begin{figure}[h!]
	\centering
	\includegraphics[scale=0.55]{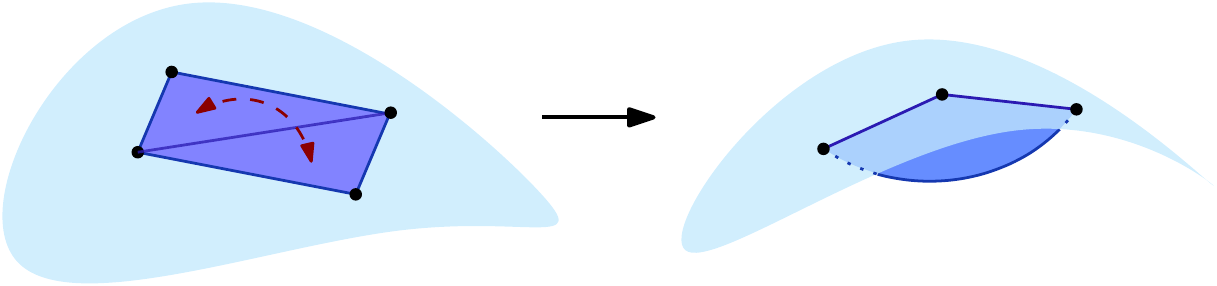}\hspace{1.5cm}\includegraphics[scale=0.55]{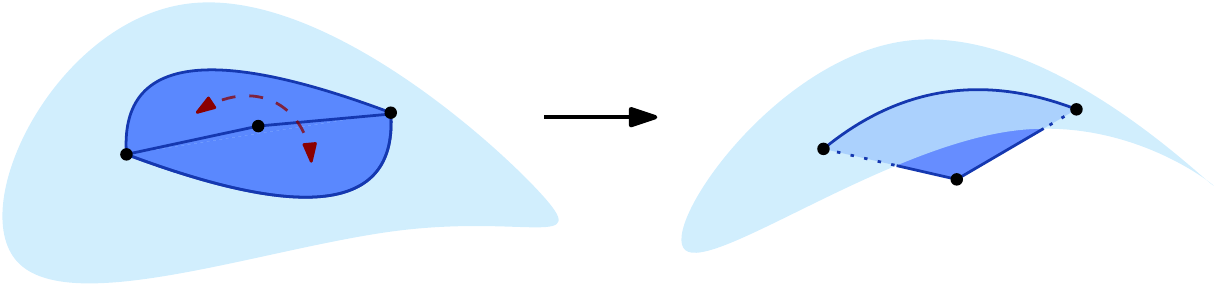}	\caption{The two possibilities for the gluing of two two boundary triangles in a step of a local construction. The obtained triangle is in the interior of the resulting triangulation.}
	\label{fig:Glue-triangles}
\end{figure}

\paragraph{Definition.}  Consider a tree of tetrahedra and two triangles on its boundary that share at least one edge (they are said to be \emph{adjacent}). One may glue these two triangles together, respectively identifying the pairs of edges of these triangles that share a vertex (Fig.~\ref{fig:Glue-triangles}). Doing so, the boundary is reduced by two triangles, but the number of tetrahedra remains constant. One may choose again  two adjacent triangles on the boundary and repeat the operation, and so on until the boundary is empty. 

Starting from a tree of tetrahedra and recursively choosing an edge on the boundary shared by two triangles and gluing these triangles until the boundary is empty is by definition a \emph{local construction} of the resulting closed three-dimensional triangulation.

We usually denote by $T_s$ the triangulation at step $s$, that is after having glued two-by-two the first $s$ pairs of triangles, 
so that $T_0$ is the initial tree of tetrahedra, and if $T_0$ has $2n$ boundary triangles, $T_n=T$, and for every $s$, $T_{s+1}$ is obtained from $T_s$ by gluing two triangles $A_{s+1}$ and $B_{s+1}$ on the boundary of $T_s$. By convention, if $A_{s+1}$ and $B_{s+1}$ share more than one edge, we also chose one of these edges\footnote{This is not the usual convention, but it allows for simpler statements of the results of the paper and does not change the notion of locally constructible triangulation. }.

A triangulation that admits a local construction is said to be \emph{locally constructible}.
Note that a locally constructible triangulation $T$ may a priori admit different local constructions: it may a priori be obtained starting from different trees of tetrahedra, corresponding to different spanning trees on the dual graph of $T$, and for each tree of tetrahedra, there might be different local constructions leading to the same final result (for instance for different choices of edges shared by two triangles, or carrying out the gluings of triangles in a different order, see Sec.~\ref{subsub:Morse-and-ineq}). In particular, a given local construction provides a distinguished spanning tree of tetrahedra. 
By abuse of notation, if  
 a local construction of $T$ is initiated with the tree of tetrahedra $T_0$, we also denote by $T_0$ this spanning tree of tetrahedra of $T$. We say that the local construction 
 is \emph{based on $T_0$}.
 
 \paragraph{Some history.}  
Locally constructible triangulations were introduced by Durhuus and Jonsson in \cite{Durhuus-Jonsson} for two important reasons. First of all the topology of such a triangulation is fully characterized: every locally constructible triangulation has the topology of the 3-sphere. Second, the family of triangulations is exponentially bounded in the number of tetrahedra. 
More precisely, there exists a $C>0$ such that the number of locally constructible 3-dimensional triangulations with $n$ tetraedra is bounded from above by $C^n$.
If all 3-sphere triangulations were locally constructible, as the authors in \cite{Durhuus-Jonsson} conjectured, this would imply a positive answer to the problem of exponential boundedness of triangulated 3-spheres.
However, this turns out not to be the case: Benedetti and Ziegler have shown that not all triangulations of the 3-sphere are locally constructible \cite{Benedetti-Ziegler}, leaving the problem thus wide open.

\subsubsection{Local constructibility and collapsibility}
\begin{figure}[h]
	\centering
	\includegraphics[width=.8\linewidth]{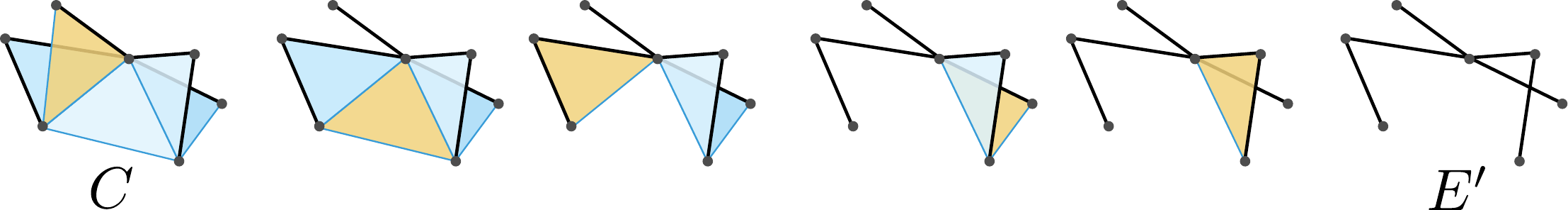}
	\caption{An example of a collapsing sequence $C \searrow E'$, where in this case $C$ is a tree of triangles and $E'$ a tree of free edges of $C$ (the setting of Lem.~\ref{TrOfTr-exists-collapse} below). The triangle $\Sigma$ removed at each step is indicated in orange.\label{fig:collapse}}
\end{figure}

An \emph{elementary collapse} is the removal from a complex $C$ of a $k$-simplex $\Sigma$ and a $(k-1)$-simplex $\sigma$, such that $\sigma$ is contained in $\Sigma$, and $\sigma$ is contained in no other cell of $C$ of dimension $k$ or larger ($\sigma$ is called a \emph{free simplex} of $C$). A \emph{collapsing sequence} is a sequence of elementary collapses. A complex $C$ is said to \emph{collapse onto} a complex $C'$ (denoted by  $C\searrow C'$), if $C'$ can be obtained from $C$ by a collapsing sequence. See Fig.~\ref{fig:collapse} for an example.

The following is a refinement of some developments  of \cite{Benedetti-Ziegler}.
\begin{theorem}
\label{thm:LC-and-collapse}
Let $T$ be a 3-dimensional triangulation, $T_0$ a spanning tree of tetrahedra. 
There is a bijection between:
\vspace{-0.2cm}
\begin{itemize}
\item A local construction of $T$ based on $T_0$.
\item  A collapsing sequence of $T^{T_0}$ onto $E$, where $E$ is a spanning tree of edges of $T$. 
\end{itemize}
\end{theorem}
In the following, given a local construction of $T$ based on $T_0$, we refer to the unique spanning tree $E$ of the corresponding collapsing sequence of $T^{T_0}$ onto $E$  as  \emph{its critical tree.}

\proof  Given a local construction, at every step we glue two adjacent triangles $A_{s+1}$ and $B_{s+1}$ on the boundary of $T_s$, and for every $s$, we have a chosen  edge  $\sigma_{s+1}$ shared by $A_{s+1}$ and $B_{s+1}$ on the boundary of $T_s$. We  let $\Sigma_{s+1}$ be the triangle of $T^{T_0}$ resulting from the gluing of $A_{s+1}$ and $B_{s+1}$. 

Let $C_s$ be the subcomplex of $T^{T_0}$ obtained by removing $\sigma_1, \Sigma_1, \ldots, \sigma_s, \Sigma_s$. 
At a step $s$ of the local construction, the triangles of $C_s$ correspond to the subset of triangles of $T^{T_0}$ that have not been glued together yet and are still on the boundary of $T_s$, while the triangles $\Sigma_1, \ldots, \Sigma_s$ removed correspond to the triangles of $T^{T_0}$ that are in the interior of $T_s$. The key observation is that an edge of $T_s$ is shared by exactly two boundary triangles  if and only if  all the other triangles incident to it are in the interior of $T_s$. Since  $\sigma_{s+1}$ is shared by exactly two boundary triangles on $T_s$ and since we have removed all the triangles of $T^{T_0}\cap \mathrm{Int}(T_s)$, $\sigma_{s+1}$  is therefore a free simplex of $C_s$. While gluing $A_{s+1}$ and $B_{s+1}$ on the boundary of $T_s$, we remove $\sigma_{s+1}$ and $\Sigma_{s+1}$ from $C_s$, and 
this is an elementary collapse from $C_s$ onto $C_{s+1}$. A local construction and a choice of an edge share by the two triangles to be glued at every step thus corresponds to a collapsing sequence from $C_0 = T^{T_0}$ to $E:=C_n$, where $n$ is the number of triangles of $T^{T_0}$. $E$ does not have any triangles, so it is a 1-dimensional subcomplex of $T^{T_0}$. Since $T$ is a triangulated sphere, $T^{T_0} $ must be contractible \cite[Prop.~2.4]{Benedetti-Ziegler}, and therefore $E$ must be contractible, which for a 1-dimensional complex is equivalent to $E$ being a tree. It must be spanning in $T$, since we never removed any vertex of  $T^{T_0}$ in the collapse, and the vertices of $T$ all belong to $T^{T_0}$. 

Reciprocally, let $C_1, \ldots, C_p$  be a sequence of subcomplexes of $T^{T_0}$ such that $C_p:=E$ is a spanning tree of edges of $T$, and $C_{s+1}$ is obtained from $C_s$ by the elementary collapse of  $\sigma_{s+1}$ and $\Sigma_{s+1}$. Since $E$ is spanning, all the vertices of  $T^{T_0}$ are still in $E$, so that for $1\le s \le p$, $\sigma_{s}$ is an edge (and not a vertex), so that $p=n$. Let $T_s$ be obtained from $T$ by ungluing all the triangles of $T_0$ in $C_s$, as well as all the vertices and edges if they do not belong to two triangles of $T_0$ glued in $T^{T_0}\setminus C_s$. At every step, $\sigma_s$ belongs to a single triangle $\Sigma_s$  in  $C_s$, so $\sigma_s$ belongs to exactly two boundary triangles $A_{s+1}$ and $B_{s+1}$ in $T_s$, and these two triangles are adjacent  on the boundary of $T_s$. This defines a local construction of $T$ based on $T_0$.
\qed

\subsubsection{Triple-tree triangulations are locally constructible}

We now focus on the situation where $T^{T_0}_E$ is a tree of triangles. 
The following lemma is illustrated in Fig.~\ref{fig:collapse}.

\begin{lemma}
\label{TrOfTr-exists-collapse}
Let $C$ be a tree of triangles with a marked spanning tree of free edges $E'$. Then $C\searrow E'$.
\end{lemma}
\proof 
We proceed by induction: if $C$ has a single triangle, $E'$ consists of two of its edges, and removing the third edge and the triangle is a collapsing sequence from $C$ onto $E'$. Otherwise, the set of free edges of a tree of $n\ge 2$ triangles has cycles, and since $E'$ is a tree,  $C$ has a free edge $\sigma$ not in $E'$, and we let $v$ be the vertex of that triangle not contained in $\sigma$. After an elementary collapse of $\sigma$ and the triangle it belongs to, we have two smaller trees  of triangles $C_1$ and $C_2$ respectively with spanning trees of free edges $E'_1$ and $E'_2$ attached at $v$. $C_i$ might have no triangle, in which case it is just an edge which is in $E'_i$. From the induction, there exist collapsing sequences from $C_1$ and $C_2$ to $E'_1$ and $E'_2$, which together with the elementary collapse of $\sigma$ form a collapsing sequence from $C$ onto $E'$.
\qed

\begin{proposition}
\label{prop:Triple-Tree-implies-LC}
Let $ T\in \mathbb{T} $, with marked spanning trees of tetrahedra $T_0$ and edges  $E$. There exists a local construction based on $T_0$ whose critical tree is $E$  .  \end{proposition}
In particular, setting aside the information on $T_0, E$ and the root, this implies the weaker statement:
\begin{align}
\mathbb{T} \subset \{\textrm{locally constructible triangulations}\} \subset \{\textrm{triangulations of the 3-sphere}\}.
\end{align}

\proof Let $ T\in \mathbb{T} $. From Lemma~\ref{TrOfTr-exists-collapse},  $T^{T_0}_E\searrow \mathsf{Cut}(E)$. Consider a collapsing sequence from $T^{T_0}_E$ onto $\mathsf{Cut}(E)$.  Going from $T^{T_0}_E$ to $T^{T_0}$, only edges of $\mathsf{Cut}(E)$ are attached to each other, which are not removed in that collapsing sequence,
so that  the latter induces a collapsing sequence of $T^{T_0}$ onto $E$. From Thm.~\ref{thm:LC-and-collapse}, it corresponds to a local construction based on $T_0$ whose critical tree is $E$.  \qed

\subsection{Tree-avoiding local constructions}
\label{sub:TALC}

\begin{figure}[h!]
	\centering	
	(a)\hspace{-0.1cm}\includegraphics[scale=0.55]{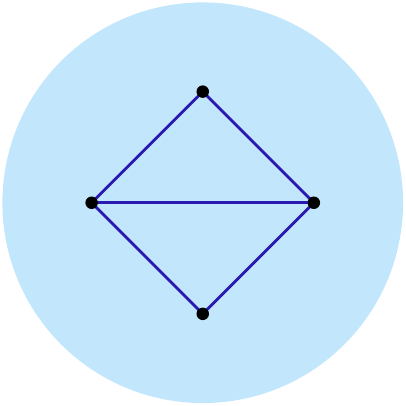}\hspace{0.2cm} 
	(b)\hspace{-0.1cm}\includegraphics[scale=0.55]{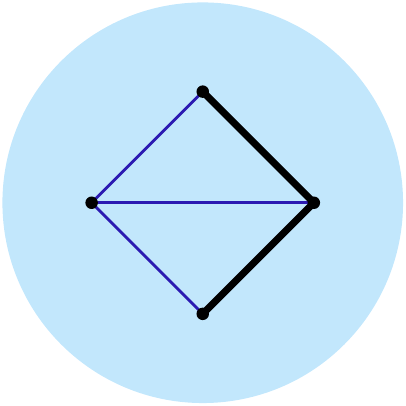}\hspace{0.2cm} 	
	(c)\hspace{-0.1cm}\includegraphics[scale=0.55]{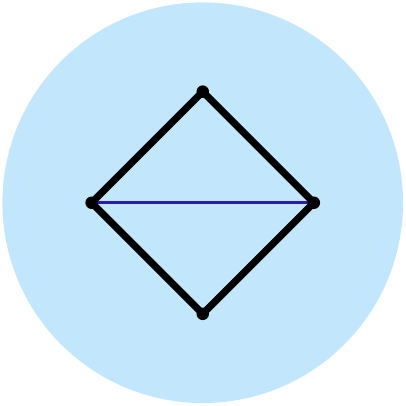}\hspace{0.2cm} 		
	(d)\hspace{-0.1cm}\includegraphics[scale=0.55]{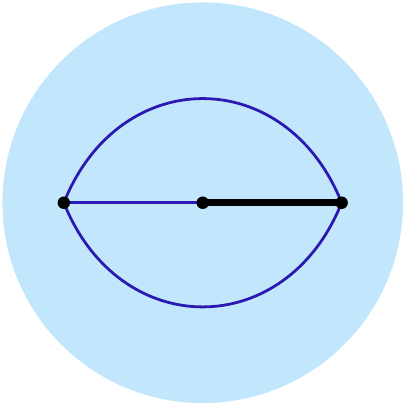}\hspace{0.2cm} 
	(e)\hspace{-0.1cm}\includegraphics[scale=0.55]{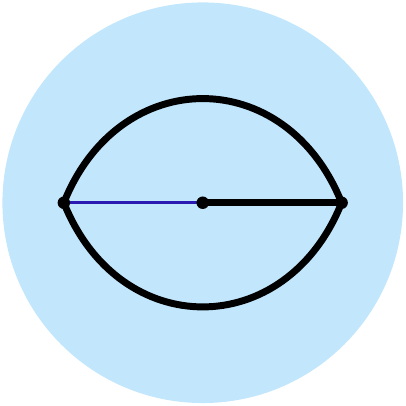}\\
	(f)\hspace{-0.1cm}\includegraphics[scale=0.55]{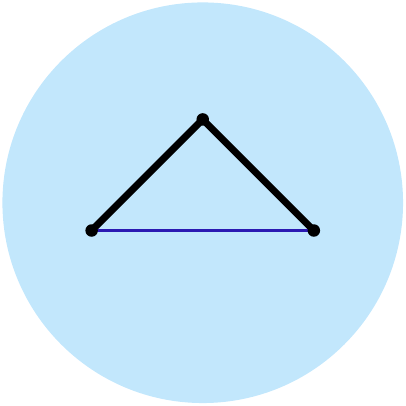}\hspace{0.2cm} 
	(g)\hspace{-0.1cm}\includegraphics[scale=0.55]{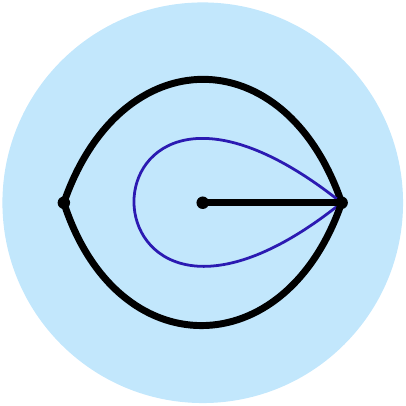}\hspace{0.2cm} 
	(h)\hspace{-0.1cm}\includegraphics[scale=0.55]{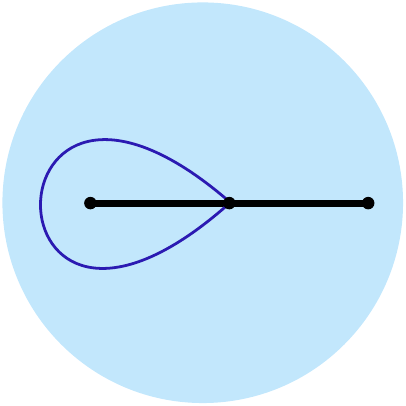}
\caption{Eight possibilities for two admissible boundary triangles. The bold edges belong to $E_s$. The represented edges are  distinct, but the vertices might not be.  Cases (f) and (h) may only occur when the boundary component is reduced to two triangles.}
	\label{fig:Possibilities-triangles-3trees}
\end{figure}

\begin{figure}[h!]
	\centering
	(i)\hspace{-0.1cm}\includegraphics[scale=0.5]{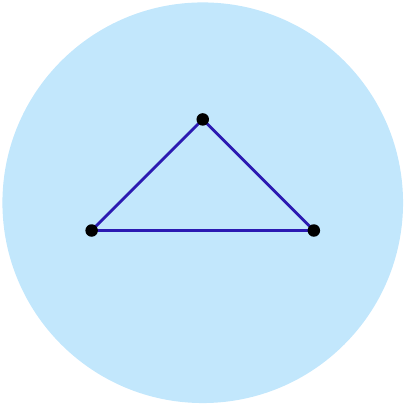}\hspace{0.15cm} 
	(j)\hspace{-0.1cm}\includegraphics[scale=0.5]{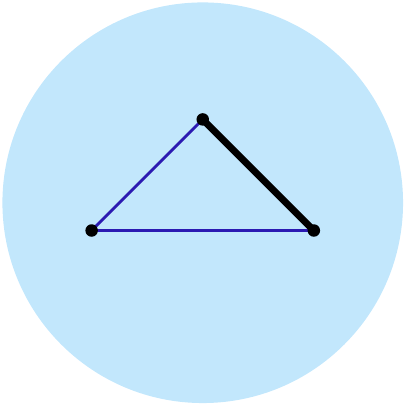}\hspace{0.15cm} 
	(k)\hspace{-0.1cm}\includegraphics[scale=0.5]{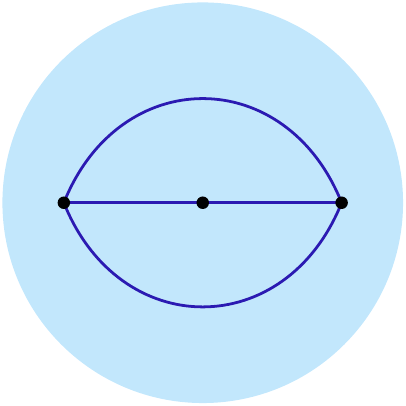}\hspace{0.15cm} 	
	(l)\hspace{-0.1cm}\includegraphics[scale=0.5]{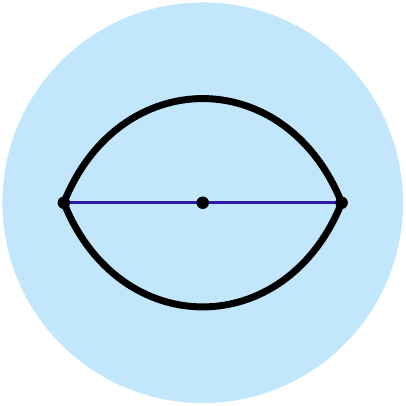}\hspace{0.15cm} 
	(l)\hspace{-0.1cm}\includegraphics[scale=0.5]{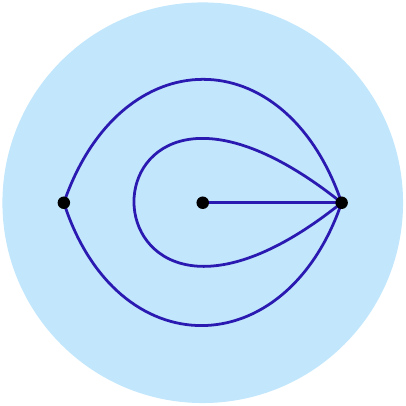}\hspace{0.15cm} 
	(l)\hspace{-0.1cm}\includegraphics[scale=0.5]{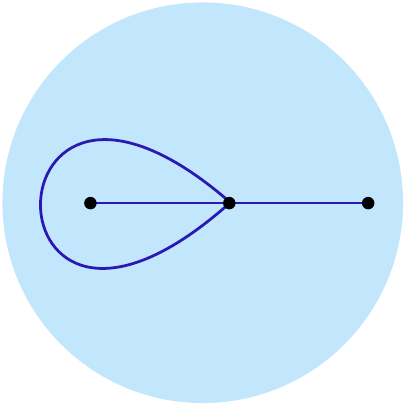}
\caption{Although {\it a priori} admissible, these six possibilities never occur.
}
	\label{fig:NonPossibilities-triangles-3trees}
\end{figure}

Given a local construction of a closed 3-dimensional triangulation $T$ in $S$ steps, we denote by $T_s$ the triangulation at step $s\in\{1,\ldots, S\}$ (so that $T_0$ is the initial tree of tetrahedra and $T_{S}=T$). 
A subset of edges $E'$ of a triangulation is said to be \emph{spanning} if each vertex in the triangulation has at least one incident edge in $E'$. Note that $E'$ need not be connected.
Let $E_0$ be a distinguished spanning subset of edges of $T_0$. 
We will introduce a restriction on the possible triangles that can be glued in the local constructions. This new rule for gluing adjacent triangles uses a spanning 
subset $E_s$ of edges of $T_s$ for $s\ge 2$, which is defined recursively, starting from $E_0$: 
\begin{itemize}
\item \textit{Rule for gluing adjacent triangles}: 
at every step $s$, the gluing of a pair of adjacent boundary triangles is only allowed if one of the edges they share is \emph{not} in $E_s$, and if the pair(s) of edges that are to be  identified upon gluing are of matching type (both in $E_s$ or both not). In this case the pair of triangles is said to be \emph{admissible} (at step $s$). 
If we differentiate the cases based on which edges the two triangles share, this {\it a priori} leaves fourteen possibilities, shown in Fig.~\ref{fig:Possibilities-triangles-3trees} and Fig.~\ref{fig:NonPossibilities-triangles-3trees}. 

\item \textit{Definition of $E_s$}: The set $E_{s+1}$ of distinguished edges in $T_{s+1}$  is then defined as the image of $E_s$ after the step $s+1$, that is: 
\begin{enumerate}[label=-] 
\item Any edge of $E_s$ that is not identified with another edge in the gluing also belongs to $E_{s+1}$; 
\item If two edges in $E_s$ are identified in the gluing, then the resulting edge belongs to $E_{s+1}$.
\end{enumerate}
\end{itemize}

A local construction for which at each step the pair of glued triangles is admissible is said to \emph{avoid the spanning 
set} $E_0$. A local construction that avoids a spanning \emph{tree} $E_0$ is called a \emph{tree-avoiding local construction}. 
Note that there exist locally constructible triangulations that avoiding a spanning set $E_0$ even when $E_0$ is not a tree (one may for instance have a number of triangles as in case (c) of Fig.~\ref{fig:Possibilities-triangles-3trees} on the boundary of $T_0$). 

\ 

For $s\ge 1$, we denote by $\mathrm{Int}(T_{s})$ the interior of $T_{s}$, and by $\partial T_{s}$ the boundary of $T_s$.

\begin{lemma}
\label{lemma:distinguished-set-spanning-connected}
If $E_s$ is connected and spanning in $T_s$, then $E_{s+1}$ is spanning in $T_{s+1}$, and $E_{s+1}\cap\partial T_{s+1}$ is spanning on  $\partial T_{s+1}$. In particular, the cases of Fig.~\ref{fig:NonPossibilities-triangles-3trees} never occur in a 
set-avoiding local construction.
\end{lemma}

\proof The first property is clear from the construction. The second property is shown recursively: assuming the properties of the lemma to be true at step $s$,
then the cases of  Fig.~\ref{fig:NonPossibilities-triangles-3trees} cannot occur; 
for the case (a) all the vertices stay on the boundary, they all have incident edges in $E_s\cap\partial T_s$ (which is connected) and after the gluing these edges belong to  $E_{s+1}\cap\partial T_{s+1}$ (and the number of connected components cannot increase); 
for the case (b) the same occurs but the leftmost vertex and at least one of the remaining vertices must be linked by a path of edges that belong to $E_s\cap\partial T_{s}$, 
which is still true after the gluing,
and similarly for the case (d) with the leftmost and the rightmost vertices. For the cases (c), (e), (f), (g) and (h) there is nothing to prove. \qed

\ 

From now on, two adjacent triangles are only said to be admissible if they are as in one of the cases of Fig.~\ref{fig:Possibilities-triangles-3trees}. As a consequence, there is always a single non-distinguished edge shared by two admissible triangles, and we always choose this edge when gluing two admissible triangles.

Note  that only for the cases  (d), (e),  (f) and (h) of Fig.~\ref{fig:Possibilities-triangles-3trees}, does an edge in $E_s$ on the boundary of $T_s$ become an edge in $E_{s+1}$ in the interior of $T_{s+1}$. For the  other cases  of Fig.~\ref{fig:Possibilities-triangles-3trees},   $E_{s}\cap \mathrm{Int}(T_{s})=E_{s+1}\cap \mathrm{Int}(T_{s+1})$. 

\

If $E$ is a set of edges of $T$ and $T$ is obtained from a tree of tetrahedra $T_0$ by attaching two-by-two the triangles on its boundary, then we call \emph{preimage of $E$ on the boundary} $\partial T_0$ of $T_0$ the set of edges of $\partial T_0$ which are attached together  in $T$ to form the edges of $E$. 
\begin{proposition}
\label{prop:TALC_iff_LC_and_preimage}
Consider a local construction  of $T$ based on $T_0$ and whose critical tree is $E$. Then this local construction avoids  the preimage of $E$ on $\partial T_0$.
This local construction is tree-avoiding if and only if the preimage of $E$ on $\partial T_0$ is a tree.
\end{proposition}
\proof It is a simple consequence of the definitions of preimage and of critical tree. Regarding the first statement:  the critical tree $E$ of a local construction of $T$ based on $T_0$ is the set of edges of $T$ whose preimages $E_0$ on $\partial T_0$ consist of edges which are not selected at any step of the local construction. The set of edges in the complement of $E_0$ is partitioned in groups, such that the edges of a groups are identified together to form an edge which at some step $s$ will be a non-distinguished edge shared by  two triangles to be glued, which will result in an edge in the complement of $E$ in $T$.   The rule for the admissibility is always verified by two triangles to be glued by definition of the preimage: edges of the boundary at step $s$ which are identified are either both in the preimage  of $E$ on $T_s$ or they are either not. 

Regarding the second statement: given a local construction of $T$ that avoids a tree $E_0$, the image of $E_0$ in $T$ is the critical tree for that local construction, so its preimage is $E_0$, a tree. Reciprocally, if the preimage $E_0$ of $E$ on $\partial T_0$ is a tree, then from what we have just proven, this local construction avoids $E_0$, a tree.
\qed

\

\begin{proposition}
\label{prop:Triple-tree-are-TALC}
Let $ T\in \mathbb{T} $, with marked spanning trees of tetrahedra $T_0$ and edges  $E$. There exists a tree-avoiding local construction of $T$ based on $T_0$ whose critical tree is $E$. 
\end{proposition}
In particular, setting aside the information on $T_0, E$ and the root, this implies the weaker statement:
\begin{align*}
\mathbb{T} \subset \{\textrm{tree-avoiding locally constructible triangulations}\} \subsetneq \{\textrm{locally constructible triangulations}\}.
\end{align*}

\proof Let $T\in\mathbb{T}$ for distinguished $T_0, E$, and $(t, \pi_H, \pi_A)$ be the corresponding triple-tree (Thm.~\ref{thm:thm-1}). Then we saw in Prop.~\ref{prop:Triple-Tree-implies-LC} that there exists a  local construction based on $T_0$ whose critical tree is $E$. On the other hand, we know from Thm.~\ref{thm:thm-1} that  $E=E(t, \pi_{_\mathrm{H}},\pi_{_\mathrm{A}})$, so that the preimage of $E$ on $\partial T_0$ is $E_0(t, \pi_{_\mathrm{A}})$, a tree. We conclude using Prop.~\ref{prop:TALC_iff_LC_and_preimage}. \qed

\subsection{Tree-avoiding local-constructions are triple-tree triangulations}
\label{sub:triple-trees-are-TALC}

\subsubsection{Reduction sequences of outerplanar triangulations}

Consider a 3-dimensional triangulation $T$ with a distinguished spanning tree of tetrahedra $T_0$, and a distinguished rooted spanning tree of edges $E_0$ on  $\partial T_0$, such that there exists a tree-avoiding local construction of $T$ based on $T_0$ that avoids $E_0$, and consider such a local construction. We denote by $T_s$ the triangulation at step $s$, and by $(A_s, B_s)$, $1\le s \le n$, the admissible pair of triangles glued on the boundary of $T_{s-1}$ to obtain $T_s$. We denote by   $\Pi$ the pairing $\{(A_s, B_s)\}_{1\le s \le n}$ of the boundary triangles  of $T_0$ and $\vec o$ the ordering of $\Pi$ given by the label $s$.

From Lemma~\ref{lemma:apollonian-tree-tetra}, the pair $(T_0, E_0)$  can be expressed as  $\mathsf{Glue}(t_0,\pi_{_\mathrm{A}}) \in \mathcal{A}_{n+1}$, where  $t_0\in\mathcal{O}_{n+1}$, $\pi_{_\mathrm{A}}\in\mathcal{P}_{n+1}$.  
Since $\mathsf{Glue}(t_0,\pi_{_\mathrm{A}})$ is identified with the boundary of $T_0$, $\Pi$ is a pairing of the triangles of $t_0$.

\begin{figure}[h!]
	\centering
	\includegraphics[scale=0.77]{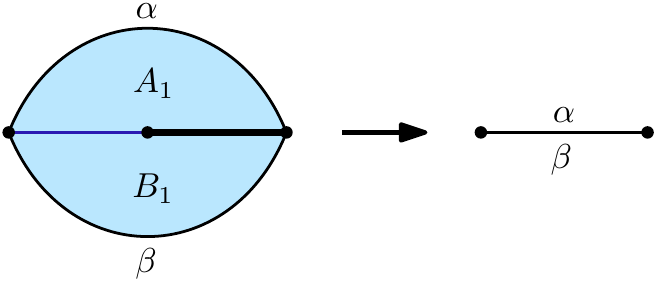}\hspace{1cm}\raisebox{6ex};\hspace{1cm}\includegraphics[scale=0.77]{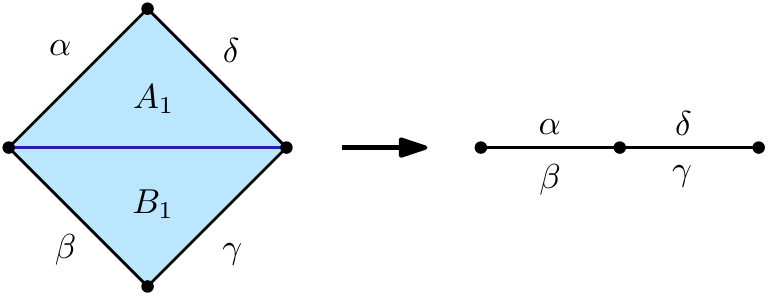}
		\caption{The gluing of $A_1$ and $B_1$ from the viewpoint of $\partial T_0$.
	For the second case, the edges adjacent  
	to   $\delta$ and $\gamma$ 
	are 
	either distinguished of not.}
	\label{fig:Operation-a}
\end{figure}

At the first step of the local construction, the triangles $A_1$ and $B_1$  are glued together on $\partial T_0$. From the viewpoint of $\partial T_0$, this corresponds to one of the two possibilities illustrated in Fig.~\ref{fig:Operation-a}. For the second case, the edges adjacent 
to   $\delta$ and $\gamma$ are 
either distinguished  (i.e.~in $E_0$) or not.

\begin{figure}[h!]
	\centering
	\includegraphics[scale=0.75]{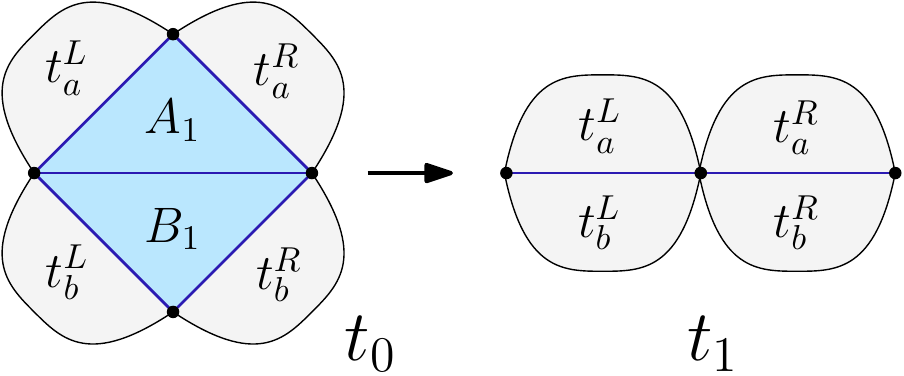}	\hspace{0.6cm}\raisebox{7ex};\hspace{0.6cm} \raisebox{6ex}{\includegraphics[scale=0.7]{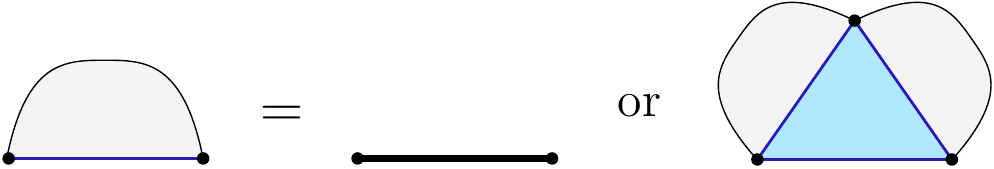}}
	\caption{The gluing of $A_1$ and $B_1$ viewed on $t_0$, and inductive structure of the ``blobs''.}
	\label{fig:Operation-b}
\end{figure}

This can be seen as performing the operation on the left of Fig.~\ref{fig:Operation-b} directly on $t_0$, where the ``blobs'' $t_a^L$,  $t_a^R$, $t_b^L$,  $t_b^R$  represent possibly empty parts of $t_0$ as illustrated on the right of Fig.~\ref{fig:Operation-b}.

It results in a new planar map $t_1$ whose 
outer face might
not be simple anymore (left of Fig.~\ref{fig:Operation-b}).
In $t_1$, the blobs $t_a^L$ and  $t_b^L$ (resp.~$t_a^R$ and  $t_b^R$) become adjacent, forming a planar map $t_1^L$ (resp.~$t_1^R$), so that $t_1^L$ and $t_1^R$ have one common vertex in $t_1$.
$t_1^L$ or $t_1^R$ might be reduced to one edge linking two vertices.  If  $t_1^L$ or $t_1^R$ is not reduced to an edge, then it is a (possibly non-rooted) outerplanar triangulation.

We then proceed to a second step: gluing  $A_2$ and $B_2$ on the boundary of $T_1$ can be seen locally\footnote{The outer face is not necessarily simple, and after the first step, there can be other portions of the map attached to the four vertices of the two triangles $A_s$ and $B_s$.} as performing the operation of Fig.~\ref{fig:Operation-a}  for $A_2$ and $B_2$ on $t_1$, obtaining a planar map $t_2$, and so on, yielding a sequence $t_0, t_1, \ldots, t_n$ of planar maps, where for $1\le s \le n$, $t_s\in \mathcal{N}_{n-s}$, where $\mathcal{N}_{k}$ is the set of rooted planar maps with a total of $2k+1$ faces, such that all the faces apart from the outer face are triangles,  and so that all the vertices lie on the outer face. 

\paragraph{Reduction sequence.} Two adjacent triangles of $t\in \mathcal{N}_k$ are said to be \emph{admissible on $t$} if they have either no boundary edges, or one boundary edge each so that the two boundary edges share a vertex, or two boundary edges each. 
Let $t_0$ be an outerplanar triangulation with $2n\ge 2$ triangles, $\Pi$ a pairing of its triangles and $\vec o$ an ordering of $\Pi$ from $1$ to $n$. 
We say that $(t_0, \Pi, \vec o)$ is a \emph{reduction sequence} of $t_0$ if at each step $1\le s\le n$, the pair of triangles of $\Pi$ numbered $s$ in $\vec o$, which we denote by $\{A_s, B_s\}$, satisfies: 
\begin{itemize}
\item $A_s$ and $B_s$ are admissible on $t_{s-1}$, 
\item $t_{s}$ is obtained from $t_{s-1}$ by performing locally the operation of Fig.~\ref{fig:Operation-b} on the pair of triangles $A_s$ and $B_s$ on $t_{s-1}$. As a consequence, $t_s\in\mathcal{N}_{n-s}$.
\end{itemize}

\begin{lemma}
\label{lem:choosing-pairs} 
A pair of triangles is admissible on the boundary of $T_s$ if and only if it is admissible on $t_s$.
\end{lemma}

\proof We consider two admissible triangles on $\partial T_s$, and cut along the distinguished subset of edges by splitting the distinguished edges in two. In that process: 
two triangles as in case (a) of Fig.~\ref{fig:Possibilities-triangles-3trees} don't have any boundary edges on $t_s$;
two triangles as in as in one of the cases (b) or (d) of Fig.~\ref{fig:Possibilities-triangles-3trees}  have one boundary edge each that share a vertex on $t_s$; 
and two triangles as in as in one of the cases (c), (e), (f), (g) or (h) of Fig.~\ref{fig:Possibilities-triangles-3trees}  have two boundary edges each on $t_s$.  \qed

\

Let us formalize the discussions of this section in a lemma:
\begin{lemma}
\label{lem:reduction-seq-vs-TALC} 
With these notations,  if starting from $T_0$  and gluing the triangles of $\Pi$ in the order $\vec o$ defines a tree-avoiding local construction, then $(t_0, \Pi, \vec o)$ is a reduction sequence.
\end{lemma}

The converse is also true, but we do not prove it. Finally, we will need the following (with the convention that a reduction sequence on a map consisting of a single edge is an empty reduction sequence). 
\begin{lemma}
\label{lem:induce-red-seq} 
With the notations above, a reduction sequence on $t_0$ induces two independent reduction sequences on $t_1^L$ and $t_1^R$ respectively.
\end{lemma}
 \proof Indeed, there are no triangle of $t_1^L$ paired with a triangle of $t_1^R$ in $\Pi$, because the recursive application of the operation of Fig.~\ref{fig:Operation-b} will never make a triangle in $t_1^L$ and a triangle in $t_1^R$ adjacent again. The ordered pairing $(\Pi, \vec o)$ therefore induces an ordered pairing $(\Pi_L, \vec o_L)$ (resp.~$(\Pi_R, \vec o_R)$) on $t_1^L$ (resp.~$t_1^R$), which defines a reduction sequence on $t_1^L$ (resp.~$t_1^R$), with $\Pi=\Pi_L\cup \Pi_R\cup \{A_1, B_1\}$. \qed

\subsubsection{Tree-avoiding local constructions are triple-tree triangulations}

The family $\mathbb{T}$ has been defined in Thm.~\ref{thm:thm-1}.

\begin{theorem} 
We fix a 3-dimensional triangulation $T$, with $T_0$ and $E$ some spanning trees of tetrahedra and edges of $T$. Then the following assertions are equivalent:
\begin{enumerate}[label=(\roman*)]
\item $T\in \mathbb{T}$ for the distinguished spanning trees  $T_0$ and $E$. 
\item There exists a tree-avoiding local construction of $T$ based on $T_0$ whose critical tree is $E$.
\end{enumerate}
\end{theorem}
\noindent
In particular, setting aside the information on $T_0, E$ and the root, this implies the weaker statement:
\begin{align*}
\mathbb{T} = \{\textrm{tree-avoiding locally constructible triangulations}\} \subsetneq \{\textrm{locally constructible triangulations}\}.
\end{align*}
As a corollary of this theorem, denoting by $E_0$ the preimage of $E$ on the boundary of $T_0$:
\begin{equation}
 E_0\,\textrm{ is a tree and }\, T^{T_0}\searrow E\qquad \Leftrightarrow \qquad T^{T_0}_E\,\textrm{ is a tree of triangles, and }\,\mathsf{Cut}(E)\,\textrm{ is a tree.}
\end{equation}

\proof The implication $(i)\Rightarrow (ii)$ has been proven in Prop.~\ref{prop:Triple-tree-are-TALC}, and we must now prove $(ii)\Rightarrow (i)$. 
Consider a tree-avoiding local construction of $T$ based on $T_0$ whose critical tree is $E$.   From Prop.~\ref{prop:TALC_iff_LC_and_preimage}, the avoided tree $E_0$ is the preimage of $E$ on $\partial T_0$. Since $T$ is rooted at a marked oriented edge of $E$ and a marked triangle of $T^{T_0}$ containing this edge, it singles out one of the preimages of this edge on $\partial T_0$, so that $E_0$ is rooted.
Assuming that $T$ has $n-1$ tetrahedra, $n\ge 2$, we may therefore identify $\partial T_0$ equipped with $E_0$ as $\mathsf{Glue}(t_0, \pi_{_\mathrm{A}})\in \mathcal{A}_{n+1}$,  $t_0\in\mathcal{O}_{n+1}$ and $\pi_{_\mathrm{A}}\in\mathcal{P}_{n+1}$ (so that $E_0=E_0(t_0,\pi_{_\mathrm{A}})$). 
Let $T_s$ be the triangulation at step $s$ of the local construction, and $(A_s, B_s)$, $1\le s \le n$ be the pair of admissible triangles glued on the boundary of $T_{s-1}$ to obtain $T_s$.

We want to show the existence of a non-crossing pairing $\pi_{_\mathrm{H}}$ of the boundary edges of $t_0$ such that  $\mathsf{Glue}(t_0,\pi_{_\mathrm{H}})\in \mathcal{H}_{n+1}$, and $T^{T_0}= \mathsf{Id}_{\pi_{_\mathrm{H}}}[\mathsf{Glue}(t_0, \pi_{_\mathrm{A}})]$. Let us start with the first point. To simplify the discussion, we respectively denote by $\Hnr_n$
and by  $\Onr_n$ the sets of \emph{unrooted} elements of $\mathcal{H}_{n}$
and $\mathcal{O}_{n}$.

\begin{lemma}
\label{lem:red-sequence-to-hierarch}
Consider $t_0\in\Onr_{n+1}$ with $2n\ge2$ triangles and a pairing $\Pi$ of its triangles such that there exists a reduction sequence $(t_0, \Pi, \vec o)$. Then we can construct a non-crossing pairing $\pi = \pi (t_0, \Pi, \vec o)$ of its boundary edges such that the triangulation $\mathsf{Glue}(t_0,\pi)$ is in $\Hnr_{n+1}$.
\end{lemma}

\begin{figure}[h!]
	\centering
	\raisebox{5ex}{\includegraphics[scale=1]{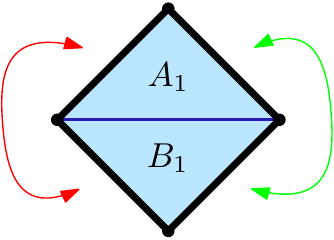}}\hspace{1.2cm}\includegraphics[scale=1]{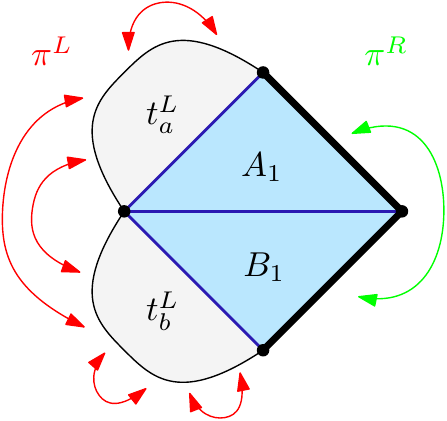}	\hspace{1cm}\includegraphics[scale=01]{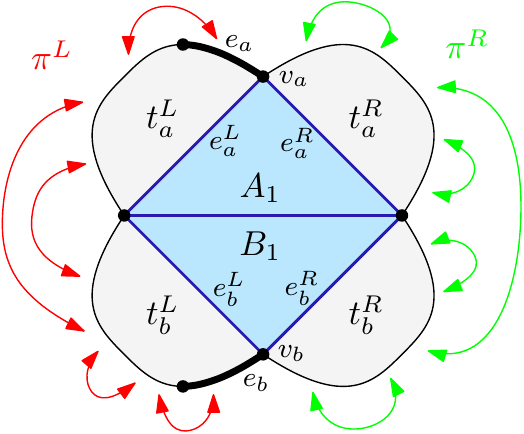}	
	\caption{The triangles $A_1$ and $B_1$ on $t_0$, and the pairing $\pi$.}
	\label{fig:Def-pi}
\end{figure}

\proof We construct this pairing recursively on the number of triangles in $t_0$. If $t_0$ has only two triangles, its boundary edges are paired in $ \pi (t_0, \Pi, \vec o)$ as on the left of Fig.~\ref{fig:Def-pi}, and $\mathsf{Glue}(t_0,\pi)$ is the only element of $\Hnr_{2}$ (case (f) of Fig.~\ref{fig:Possibilities-triangles-3trees}).

Otherwise we assume the property to be true for $t_0$ with $2k\ge2$ triangles, $k\le n-1$, and consider  $t_0$ with $2n\ge4$ triangles.
We perform the operation of Fig.~\ref{fig:Operation-b} on the first pair of triangles $(A_1, B_1)$ in the reduction sequence, obtaining $t_1$. 
With the notations of Fig.~\ref{fig:Operation-b}, in $t_1$, the blobs $t_a^L$ and  $t_b^L$ (resp.~$t_a^R$ and  $t_b^R$) become adjacent (Fig.~\ref{fig:Operation-b}), forming a planar map $t_1^L$ (resp.~$t_1^R$), so that $t_1^L$ and $t_1^R$ have one common vertex in $t_1$. 
We recall that $t_1^L$ or $t_1^R$ might be reduced to one edge linking two vertices (but not both since $t_0$ has more than two triangles).  If  $t_1^L$ is not reduced to an edge, then it is an outerplanar triangulation, and from Lemma~\ref{lem:induce-red-seq},  it admits a reduction sequence $(t_1^L, \Pi_L, \vec o_L)$, and the same goes for $t_1^R$.

If  $t_1^L$ (resp.~$t_1^R$) is not reduced to an edge, we may therefore apply the recursion hypothesis: we can construct a non-crossing pairing $\pi^L =  \pi (t_1^L, \Pi_L, \vec o_L)$ (resp.~$\pi^R=  \pi (t_1^R, \Pi_R, \vec o_R)$)  of its boundary edges such that $\mathsf{Glue}(t_1^L,\pi^L)\in\Hnr_{k+1}$  (resp. $\mathsf{Glue}(t_1^R,\pi^R)\in\Hnr_{k+1}$)  for some $k\le n-1$.
\

We define a non-crossing pairing $\pi=\pi(t_0, \Pi, \vec o)$  of the boundary edges of  
$t_0$ as follows (see Fig.~\ref{fig:Def-pi}): 
\begin{enumerate}[label=-]
\item if $t_1^L$ (resp.~$t_1^R$) is not reduced to an edge in $t_1$, then from the recursion hypothesis, we can pair the boundary edges of $t_a^L$ and  $t_b^L$ using $\pi^L$ (resp.~pair the boundary edges of $t_a^R$ and  $t_b^R$ using $\pi^R$). 
\item if $A_1$ and $B_1$ have one boundary  edge each in $t_0$ (so that $t_1^L$ or $t_1^R$ are reduced to a single edge), then these two boundary  edges are paired in $\pi$, and we respectively call $\pi^L$ or $\pi^R$ this pairing,
 \end{enumerate}

We respectively let $v_a$ and $v_b$ be the vertices of $t_0$ contained in $A_1$ and $B_1$ and  which do not belong to the edge shared by $A_1$ and $B_1$ on $t_0$ (see the notations on the right of Fig.~\ref{fig:Def-pi}).

 \begin{lemma}
 \label{lem:A0B0-three-vert}
$v^a$ and $v^b$ are identified in $\mathsf{Glue}(t_0,\pi)$.
\end{lemma}
\proof If $A_1$ and $B_1$ have one boundary  edge each, this is clear since these two edges are paired in $\pi$. We now assume this not to be the case, and call $v$ the only vertex that belongs to both  $t_1^L$ and $t_1^R$ on $t_1$ (so that $v$ results from the identification of $v_a$ and $v_b$ when going from $t_0$ to $t_1$), and $e_a$, $e_b$ the boundary edges  of  $t_a^L$ and  $t_b^L$ respectively incident to $v_a$ and $v_b$ (see the notations on the right of Fig.~\ref{fig:Def-pi}).

We call \emph{chain} of length $2k$ on $(t_1^L, \pi^L)$ a sequence $(e_1,  e_2, \ldots, e_{2k})$, where  the $e_i$ are boundary  edges of $t_1^L$ and  for which for each $1\le i \le k$, $\{e_{2i-1}, e_{2i}\}$ are paired in $\pi^L$ and for $1\le i \le k-1$, $\{e_{2i}, e_{2i+1}\}$,  share a vertex on $t_1^L$.
 There necessarily exists a chain on $(t_1^L, \pi^L)$ for which $e_1$ is $e_a$ and $e_{2k}$ is $e_b$ (this chain forms a vertex of $\mathsf{Glue}(t_1^L,\pi^L)$). This chain also exists on $(t_0, \pi)$, so that $v_a$ and $v_b$ must be the same vertex in $\mathsf{Glue}(t_0,\pi)$. \qed

 \
 
 \begin{figure}[h!]
	\centering
\includegraphics[scale=0.9]{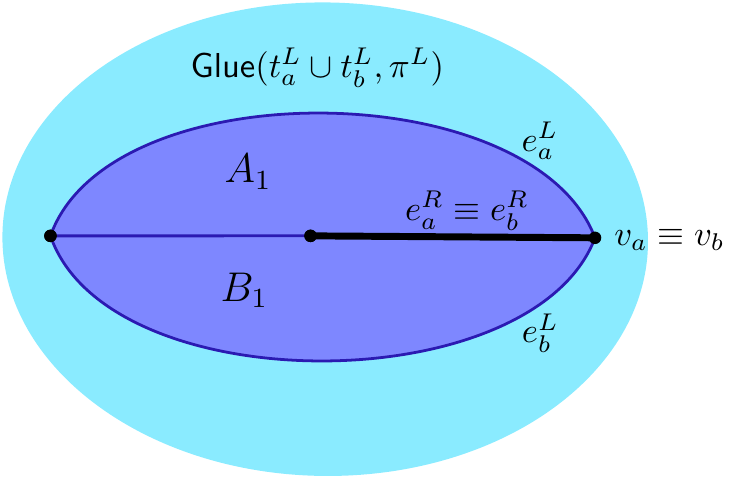}\hspace{1cm}\includegraphics[scale=0.9]{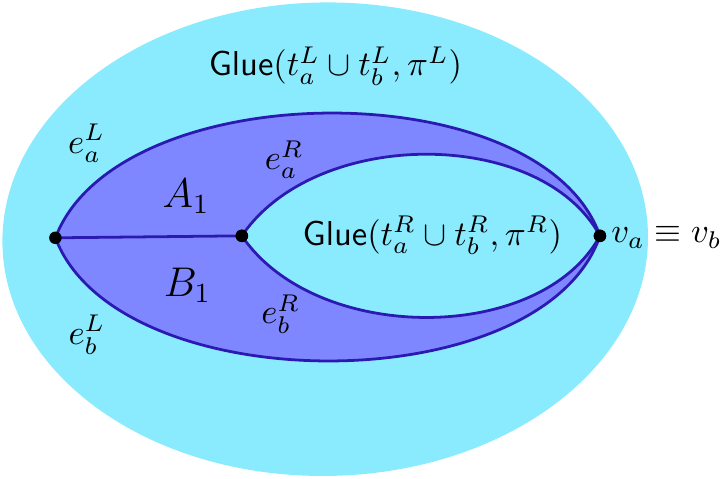}	
	\caption{$\mathsf{Glue}(t_0,\pi)$ for the middle case of Fig.~\ref{fig:Def-pi} (left) and for the case on the right of Fig.~\ref{fig:Def-pi} (right).}
	\label{fig:glue-t0-pi}
\end{figure}

 On $t_0$, we respectively call $e_a^L$,  $e_b^L$, $e_a^R$ and  $e_b^R$, the edges of $A_1$ and $B_1$ contained in $t_a^L$,   $t_b^L$,  $t_a^R$ and $t_b^R$ (Fig.~\ref{fig:Def-pi}). 
 From Lemma~\ref{lem:A0B0-three-vert}, we see that on $\mathsf{Glue}(t_0,\pi)$ (Fig.~\ref{fig:glue-t0-pi}):
\begin{itemize}
\item If $t_1^L$ (resp.~$t_1^R$) is reduced to a single edge, then $e_a^L$ and  $e_b^L$ (resp.~$e_a^R$ and  $e_b^R$) are the same distinguished edge in $\mathsf{Glue}(t_0,\pi)$. 
\item If $t_1^L$ (resp.~$t_1^R$) is not reduced to an edge, then $e_a^L$ and  $e_b^L$ (resp.~$e^a_R$ and  $e_b^R$) are distinct and form a cycle of length two that splits  $\mathsf{Glue}(t_0,\pi)$ in two parts, one of which is  $\mathsf{Glue}(t_a^L\cup t_b^L,\pi^L)$ (resp.~$\mathsf{Glue}(t_a^R\cup t_b^R,\pi^R)$).
 \end{itemize}
 The second point can be seen as a particular case of the first point, since in that case $\mathsf{Glue}(t_a^L\cup t_b^L,\pi^L)$ (resp. $\mathsf{Glue}(t_a^R\cup t_b^R,\pi^R)$) is reduced to a distinguished edge.
 
We now show that  $\mathsf{Glue}(t_0,\pi)$ is hierarchical, that is, every triangle has a companion in $\mathsf{Glue}(t_0,\pi)$ with which it shares three distinct vertices: we have already shown this to be true for $A_1$ and $B_1$ (Lemma~\ref{lem:A0B0-three-vert}), and since if they are not empty,  $\mathsf{Glue}(t_1^L,\pi^L)\in\Hnr_{k_L+1}$  and $\mathsf{Glue}(t_1^L,\pi^L)\in\Hnr_{k_R+1}$ for some $k_L, k_R\ge 1$, it is also true  that every triangle of  $\mathsf{Glue}(t_a^L\cup t_b^L,\pi^L)$ and $\mathsf{Glue}(t_a^R\cup t_b^R,\pi^R)$ has a companion with which it shares three distinct vertices. 

It remains to see that $\mathsf{Glue}(t_0,\pi)$ has no distinguished edge in a cycle of length $2$: this is true for the edges represented in Fig.~\ref{fig:glue-t0-pi}, and it is also true for the other edges of  $\mathsf{Glue}(t_a^L\cup t_b^L,\pi^L)$ and $\mathsf{Glue}(t_a^R\cup t_b^R,\pi^R)$, because it is true in $\mathsf{Glue}(t_1^L,\pi^L)\in\Hnr_{k_L+1}$  and $\mathsf{Glue}(t_1^L,\pi^L)\in\Hnr_{k_R+1}$. This concludes the proof of Lemma~\ref{lem:red-sequence-to-hierarch}. \qed

\ 

With the notations at the beginning of the proof of the theorem, from Lemma~\ref{lem:reduction-seq-vs-TALC},  the tree-avoiding local construction induces a reduction sequence $(t_0, \Pi, \vec o)$, where  $\Pi$ is the pairing $\{(A_s, B_s)\}_{1\le s \le n}$  and $\vec o$ the ordering of $\Pi$ given by the label $s$. Applying Lemma~\ref{lem:red-sequence-to-hierarch}, we construct a non-crossing pairing $\pi _{_\mathrm{H}}= \pi (t_0, \Pi, \vec o)$ of its boundary edges such that  $\mathsf{Glue}(t_0,\pi_{_\mathrm{H}})\in\mathcal{H}_{n+1}$. It remains to prove that $T^{T_0}= \mathsf{Id}_{\pi_{_\mathrm{H}}}[\mathsf{Glue}(t_0, \pi_{_\mathrm{A}})]$, that is, we must show that if $\mathcal{T},\pi_{_\mathrm{T}},\pi_{_\mathrm{C}} $ encode the embedded 2-complex $T^{T_0}$ (Lemma~\ref{lem:encoding-pairings-2comp}), then any two triangles of  $\partial T_0$ that are glued together in $T^{T_0}$ are paired in $\Pi(\mathsf{Glue}(t_0,\pi_{_\mathrm{H}}))$, and two edges of these triangles in $\mathcal{E}(\mathcal T)$ are identified  (i.e.~are paired in $\pi_{_\mathrm{T}}$) if they link the same two vertices in $\mathsf{Glue}(t_0,\pi_{_\mathrm{H}})$.

\begin{figure}[h!]
	\centering
	\includegraphics[scale=1]{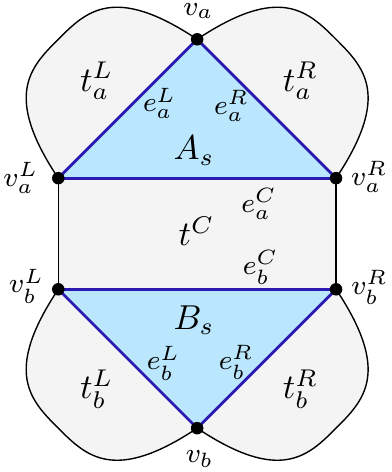}
		\caption{The triangles $A_s$ and $B_s$ on $t_0$.}
	\label{fig:notations-AsBs}
\end{figure}

Consider the pair $\{A_s, B_s\}\in \Pi$ on $t_0$. There exists a possibly empty portion $t^C$ of $t_0$ that separates these two triangles, so that they are as in Fig.~\ref{fig:notations-AsBs}, where  the notations are for the edges of $\mathcal{E}(\mathcal{T})$ that belong to the triangles $A_s, B_s$  in $\mathcal{T}$.
If $t^C$ is empty, $e_a^C$ and $e_b^C$ belong to the same edge, and if $t_a^L$ is empty, $e_a^L$ is a boundary edge, and similarly for the other ``blobs''.

\begin{lemma}
\label{lem:piT-of-LC}
With the notations above, $(e_a^C, e_b^C)\in\pi_{_\mathrm{T}}$, $(e_a^L, e_b^L)\in\pi_{_\mathrm{T}}$, and $(e_a^R, e_b^R)\in\pi_{_\mathrm{T}}$.
\end{lemma}

\proof There is a step $s$ of the local construction, at which $A_s$ and $B_s$ are adjacent on $\partial T_s$. At this step, $A_s$ and $B_s$ must also be adjacent on $t_s$. The first steps of the reduction sequence can only bring $e_a^C, e_b^C$ to be part of the same edge of $t_s$ (by progressively reducing $t^C$). When gluing $A_s$ and $B_s$ on the boundary of $T_s$,  it is therefore $e_a^L$ which is identified to $e_b^L$, and $e_a^R$ to $e_b^R$.  \qed

\begin{lemma}
\label{lem:the-LC-is-Id}
Let $(t_0, \Pi, \vec o)$ be a reduction sequence, and $\pi _{_\mathrm{H}}= \pi (t_0, \Pi, \vec o)$ constructed inductively as in the proof of Lemma~\ref{lem:red-sequence-to-hierarch}. Then $\Pi = \Pi\bigl(\mathsf{Glue}(t_0,\pi_{_\mathrm{H}})\bigr)$, and for $i\in \{C,R,L\}$, $e_a^i$ and  $e_b^i$ 
link the same two vertices in $\mathsf{Glue}(t_0,\pi_{_\mathrm{H}})$.  
\end{lemma}

\proof

We show the lemma inductively on the number of triangles of $t_0$. It is clear if $t_0$ has two triangles. Otherwise we consider $\{A_1, B_1\}$ the first pair in $(\Pi,\vec o)$. We have shown in the proof of Lemma~\ref{lem:red-sequence-to-hierarch} that $\{A_1, B_1\}$ share three distinct vertices in $\mathsf{Glue}(t_0,\pi_{_\mathrm{H}} )$, from which we deduce that these triangles are paired in $\Pi\bigl(\mathsf{Glue}(t_0,\pi _{_\mathrm{H}})\bigr)$, and that for this pair, $e^a_L$ and  $e^b_L$ (resp.~$e^a_R$ and  $e^b_R$) 
link the same two vertices in $\mathsf{Glue}(t_0,\pi_{_\mathrm{H}})$.  
Performing the operation of  Fig.~\ref{fig:Operation-b} on $A_1$ and $B_1$ on $t_0$, we get $t_1^L$ and $t_1^R$ (possibly empty), so that (Lemma~\ref{lem:induce-red-seq}) there exist two reduction sequences $(t_1^L, \Pi_L, \vec o_L)$ and $(t_1^R, \Pi_R, \vec o_R)$ (possibly trivial), where 
\begin{equation}
\label{eq:decomp-pi}
\Pi=\Pi_L\cup \Pi_R\cup \{A_1, B_1\}.
\end{equation} 
From the induction hypothesis, 
\begin{equation}
\label{eq:using-rec-inversion}
\Pi_L = \Pi\bigl(\mathsf{Glue}(t_1^L,\pi^L)\bigr), \qquad \Pi_R = \Pi\bigl(\mathsf{Glue}(t_1^R,\pi^R)\bigr),
\end{equation} 
 where  $\pi^L=\pi(t_1^L, \Pi_L, \vec o_L)$, and similarly  for $\pi^R$. Performing the operation of  Fig.~\ref{fig:Operation-b} on $A_1$ and $B_1$ on $t_0$ amounts to performing the operation described in Lemma~\ref{lem:reduce-triangles-hierarch} on these triangles on $\mathsf{Glue}(t_0,\pi_{_\mathrm{H}} )$, so that from Lemma~\ref{lem:reduce-triangles-hierarch}, 
$$
\Pi\bigl(\mathsf{Glue}(t_0,\pi )\bigr) = \Pi\bigl(\mathsf{Glue}(t_1^R,\pi^R)\bigr) \cup \Pi\bigl(\mathsf{Glue}(t_1^R,\pi^R)\bigr) \cup \{A_1, B_1\} = \Pi_L \cup \Pi_R  \cup \{A_1, B_1\}=\Pi,
$$
where we have used \eqref{eq:using-rec-inversion} for the second equality, and \eqref{eq:decomp-pi} for the third one. 

Therefore, any pair of triangles $\{A_s, B_s\}\in \Pi$ must share three distinct vertices in $\mathsf{Glue}(t_0,\pi_{_\mathrm{H}})$. By planarity, with the notations of Fig.~\ref{fig:notations-AsBs}, we must have in $\mathsf{Glue}(t_0,\pi_{_\mathrm{H}})$ $v_a \equiv v_b$, $v_a^L \equiv v_b^L$, and $v_a^R \equiv v_b^R$, which concludes the proof of the lemma. \qed

\

Lemma~\ref{lem:piT-of-LC} and Lemma~\ref{lem:the-LC-is-Id} show that for the  $\pi _{_\mathrm{H}}= \pi (t_0, \Pi, \vec o)$ of Lemma~\ref{lem:red-sequence-to-hierarch}, $T^{T_0}= \mathsf{Id}_{\pi_{_\mathrm{H}}}[\mathsf{Glue}(t_0, \pi_{_\mathrm{A}})]$, which concludes the proof of the theorem. \qed

\subsection{Morse gradients and weakly equivalent local constructions
\label{sub:Morse}
}

We have seen that triple trees naturally give rise to a local construction of the corresponding three-dimensional triangulation, thereby providing a certificate of 3-sphere topology.
In this section we discuss an alternative formulation of such a certificate in the context of Morse theory. 
It is a classical result of Reeb that if $M$ is a compact manifold and $f$ a smooth real function on $M$ with only two critical points that are both non-degenerate (non-singular Hessian), then $M$ is homeomorphic to a sphere.
In the setting of differentiable $3$-manifolds, such a function thus serves as a convenient certificate of $3$-sphere topology.
Here we focus on an analogous result in the setting of discrete Morse theory on cell complexes \cite{Forman_Morse_1998,Forman2002}, and show how it is related to local constructions.

\subsubsection{Discrete Morse spheres}

\paragraph{Discrete vector fields.} 
Let $T$ be a closed 3-dimensional triangulation, and for $k\in\{0,1,2,3\}$, let $S_k(T)$ be the set of $k$-dimensional simplices of $T$ and $S(T)=\cup_k S_k(T)$. A \emph{discrete vector field} on $T$ \cite[Def.~3.3]{Forman2002} is specified by choosing for every $k$-simplex of $T$, whether it is:
\begin{enumerate}[label=-]
\item oriented towards a $(k+1)$-simplex it belongs to,
\item or a $(k-1)$-simplex it contains is oriented towards it,
\item or none of the above, in which case this $k$-simplex is said to be \emph{critical}.
\end{enumerate}
A discrete vector field on $T$ therefore takes the form of a collection $\mathcal L$ of disjoint ordered pairs $(\eta,\eta')$, $\eta\in S_k$, $\eta'\in S_{k+1}$, $k\in\{0,1,2\}$, such that the simplex $\eta$ is contained in $\eta'$ and oriented towards it.
The critical simplices therefore are those simplices that do not occur in $\mathcal L$.

\begin{proposition} 
A closed 3-dimensional triangulation that admits a discrete vector field with the same number of odd-dimensional and even-dimensional critical simplices is a manifold. 
\end{proposition}
For example, this is the case if the vector field has the same number of critical vertices and of critical tetrahedra and no other critical simplices.

\proof We show that the Euler characteristics of that triangulation vanishes, which implies that it triangulates a manifold (see \cite[ch.~IX, \S 60, Thm.~I]{Seifert-Threlfall}). Indeed, for every $0\le k \le 3$, the number of $k$-simplices $n_k$ of a triangulation $T$ with a discrete vector field satisfies
$
n_k = n_k^c + n_k^\textrm{in} + n_k^\textrm{out},
$
where $n_k^c $ is the number of critical  $k$-simplices, and $n_k^\textrm{in} $ (resp.~$n_k^\textrm{out}$) the number of $k$-simplices with a $(k-1)$-simplex oriented towards it (resp. oriented towards a $(k+1)$-simplex). By definition of a discrete vector field, for every $0\le k \le 2$, $n_k^\textrm{out} = n_{k+1}^\textrm{in}$, and  $n_0^\textrm{in} =0$, and $n_3^\textrm{out}=0$. Moreover, by assumption, $n_0^c + n_2^c  =n_1^c  + n_3^c$. Therefore, $\chi(T) = n_3 - n_2 + n_1 - n_0 =0$. \qed

\paragraph{Walks.}  A walk on $(T, \mathcal L)$ is an ordered sequence $(\eta^1, \ldots, \eta^p)$ of simplices in $S(T)$ for some $p\ge 1$,  such that for $1\le i \le p-1$, if $\eta^i \in S_k(T)$ for some $k\in\{0,1,2,3\}$, then:
\begin{enumerate}[label=-]
\item either $\eta^{i+1}\in S_{k+1}(T)$ and $(\eta^i, \eta^{i+1})\in \mathcal L$,
\item or $\eta^{i+1}\in S_{k-1}(T)$ and $(\eta^{i+1}, \eta^{i})\notin \mathcal L$.
\end{enumerate}
That is, one can walk from a $k$-simplex towards a $(k+1)$-simplex if the former is oriented towards the latter, and one can walk from a $(k+1)$-simplex towards a $k$-simplex if the latter is not oriented towards the former.
 
A walk $(\eta^1, \ldots, \eta^p)$  is said to be a \emph{cycle} if $p>1$ and $\eta^1=\eta^p$. 
Note that a cycle necessarily alternates between $k$ and $k+1$ simplices, for some $k\in\{0,1,2\}$, because a step from a $k$ to a $k+1$ simplex must be followed by a step to a $k$ simplex again.
An \emph{acyclic discrete vector field} is a discrete vector field for which there are no cycles. 
These play an important in discrete Morse theory, because they are precisely the discrete vector fields that can appear as the gradients of discrete Morse functions \cite[Thm.~3.5]{Forman2002}.
For this reason we will refer to acyclic discrete vector fields as \emph{Morse gradients}.

An important result for our purpuses is the following by Forman.
\begin{theorem}[{\cite[Thm.~5.1]{Forman_Morse_1998}}] 
Let $T$ be a closed 3-dimensional triangulation with a Morse gradient that has two critical simplices. Then these critical simplices are a vertex and a tetrahedron, and $T$ triangulates the 3-sphere.
\end{theorem}

We call \emph{discrete Morse sphere} a 3-dimensional triangulation with a Morse gradient that has two critical simplices.

\subsubsection{A closer look on Morse gradients}
\label{subsub:Closer-on-Morse}
 
 Given a 3-dimensional triangulation $T$, there is no difficulty in building a discrete vector field that has no acyclic walks among $0$ and $1$ simplices or among $2$ and $3$ simplices, with a prescribed number of critical vertices or tetrahedra. Indeed:
 
 \begin{proposition} 
 \label{prop:acyclic-condition-on-top-bottom-layers}
A  discrete vector field on a 3-dimensional triangulation has no cycles among vertices and edges if and only if the subgraph of the 1-skeleton obtained by keeping only the edges that have a vertex oriented towards them is a forest with one critical vertex per connected component (which can be an isolated vertex). 

The vector field has no cycles among triangles and tetrahedra if and only if the subgraph of the dual graph obtained by keeping only the edges corresponding to triangles that are oriented towards a tetrahedron is a forest with one critical vertex per connected component. 
\end{proposition}

\proof By the definition of a discrete vector field, the edges that have a vertex oriented towards them naturally comprise a directed subgraph $G$ of the 1-skeleton of $T$, in which all vertices have out-degree $1$ except the critical vertices that have out-degree $0$.
We thus need to show that such a directed graph is acyclic if and only if it is a forest, with one vertex of out-degree $0$ per connected component.
One direction should be clear: if $G$ is a forest, each edge must be oriented towards the out-degree $0$ vertex in its connected component, and therefore $G$ is acyclic.
Conversely, denoting by $n_0$ the number of out-degree $0$ vertices of $G$, its number of vertices $V$ and edges $E$ satisfy $V = n_0 + E$. 
Hence the excess (first Betti number) of $G$ is $V-E+K = K - n_0 \geq 0$, where $K$ is the number of connected components. 
But for $G$ to be acyclic every connected component must have at least one out-degree $0$ vertex, so that $K \leq n_0$. 
Hence $K = n_0$. 
It follows that $G$ has $0$ excess and thus is a forest, and that indeed every connected component has exactly one out-degree $0$ vertex.

The same goes for the directed subgraph of the dual graph of $T$ in which an oriented edge corresponds to a triangle that is oriented towards a tetrahedron, but now all vertices have in-degree $1$ except for those corresponding to critical tetrahedra that have in-degree $0$.
The claimed statement reduces to the previous one upon reversal of the orientations.\qed

\

\begin{figure}[h!]
	\centering
\includegraphics[scale=0.6]{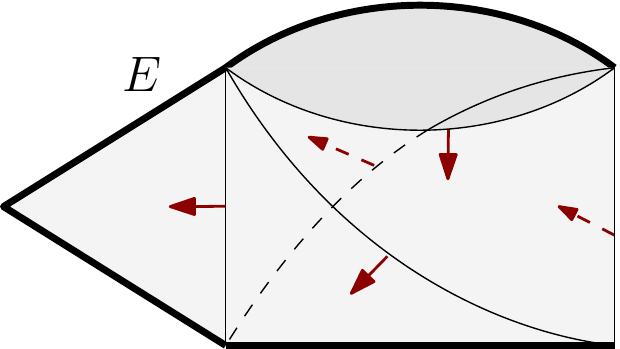}
	\caption{The 2-complex $T^{T_0}$ collapses onto a spanning tree $E$ (in bold) if and only if it admits a Morse gradient whose critical simplices are the vertices and edges of $E$.}
	\label{fig:collapsible-middle-layer}
\end{figure}

We can therefore see the existence of a Morse gradient with a prescribed number of critical simplices as a condition on certain triangles and edges in the complex. More precisely, for a discrete Morse sphere $T$, from Prop.~\ref{prop:acyclic-condition-on-top-bottom-layers}, one has a distinguished spanning tree $E$ of the 1-skeleton of $T$ which is \emph{pointed} (it has a distinguished vertex: the critical vertex) and a spanning tree of tetrahedra $T_0$ (pointed at a distinguished tetrahedron: the critical tetrahedron), and the condition that the vector field is acyclic and does not have any critical edge or triangle can be expressed as a condition on the 2-complex $T^{T_0}$ (Fig.~\ref{fig:collapsible-middle-layer}):

\begin{proposition}
Let $T$ be a 3-dimensional triangulation. The following assertions are equivalent:
\begin{enumerate}[label=(\roman*)]
\item $T$ admits a Morse gradient with two critical simplices.
\item For some spanning tree $T_0$ of tetrahedra, $T^{T_0}$ admits a Morse gradient whose critical simplices are the vertices and edges of some spanning tree $E$.
\end{enumerate}
\end{proposition}

\proof From Prop.~\ref{prop:acyclic-condition-on-top-bottom-layers}, $T$ admits a Morse gradient with a critical vertex and tetrahedron and no other critical simplices if and only if for some spanning tree $T_0$ of tetrahedra and some spanning tree $E$, it is possible to orient \emph{all} the edges of $T$ that are not in $E$  (whose set we denote by $N_1^r$) towards a triangle which is not in the interior of $T_0$  (whose set we denote by $N_2^r$) and to which they belong to, so that all the triangles in $N_2^r$ contain an edge in $N_1^r$ oriented towards them, and so that there are no cyclic walks among the elements of  $N_1^r$ and $N_2^r$. 
This is equivalent to $T$ admitting a Morse gradient whose critical simplices are all the tetrahedra, all the triangles in the interior of $T_0$, all the vertices, and all the edges in $E$, which is equivalent to $T^{T_0}$ admitting a Morse gradient whose critical simplices are all the vertices and all the edges in $E$. \qed

Two Morse gradients are said to related by \emph{re-pointing} if the only difference between the two is in the positions of the critical vertices (respectively critical tetrahedron) within the distinguished forest of the 1-skeleton (respectively dual graph). 

\subsubsection{Morse gradients and inequivalent local constructions}
\label{subsub:Morse-and-ineq}

Two local constructions are said to be \emph{weakly equivalent} if they are based on the same tree of tetrahedra $T_0$ and lead to the same 3-dimensional triangulation $T$ with the same critical tree $E$. Two weakly equivalent local constructions only differ by the order at which the admissible triangles are glued together.

\begin{proposition}
Let $T$ be a 3-dimensional triangulation. There is a bijection between:
\begin{itemize}
 \item Morse gradients on $T$ with two critical simplices up to re-pointing;
 \item local constructions of $T$ up to weak equivalence.
\end{itemize}
\end{proposition}

\proof 
From Thm.~\ref{thm:LC-and-collapse}, a local construction of $T$ based on $T_0$ with critical tree $E$ corresponds to a collapsing sequence from $T^{T_0}$ onto $E$. We use the notations of the proof of Thm.~\ref{thm:LC-and-collapse}: for every $s$, $\sigma_{s+1}$ is the chosen  edge  shared by $A_{s+1}$ and $B_{s+1}$ on the boundary of $T_s$ and  $\Sigma_{s+1}$ is the triangle of $T^{T_0}$ resulting from the gluing of $A_{s+1}$ and $B_{s+1}$, and  $C_s$ is the subcomplex of $T^{T_0}$ obtained by removing $\sigma_1, \Sigma_1, \ldots, \sigma_s, \Sigma_s$. 

We orient $\sigma_s$ towards $\Sigma_s$ in $T^{T_0}$, that is, we form the list of ordered pairs $\mathcal{L}_2=\{(\sigma_s, \Sigma_s)\}_{1\le s \le n}$. All the edges that are not in $E$ are oriented towards a triangle, and all the triangles have an edge oriented towards them, so that this defines a discrete vector field on  $T^{T_0}$ whose critical simplices are the vertices and edges of $E$. 
The same vector field is obtained from any other weakly equivalent local construction, since in each case the sequence collapses onto $E$. 
In order for an edge $\sigma_{s}$ to be a free simplex at some step $s$ of the collapse, it is necessary that all the triangles that it belongs to in  $T^{T_0}$ and towards which it is not oriented have been removed in $C_{s - 1}$. Said otherwise:
 $$
 \sigma_s \in \Sigma_{s'} \quad \textrm{and} \quad s'\neq s \quad \Rightarrow\quad  s' < s.
 $$ 
 This implies that the vector field we have defined is acyclic: if there was a cycle, starting from an edge $\sigma_{s}$ in that cycle and implementing the condition above along that cycle would lead to $s<s$, a contradiction.  To a weak equivalence class of local constructions of $T$ based on $T_0$, whose critical tree is $E$, we have therefore associated a Morse gradient  on $T^{T_0}$,  whose critical simplices are the vertices and edges of some spanning tree $E$. From the results and proofs in Sec.~\ref{subsub:Closer-on-Morse}, this corresponds to a Morse gradient on $T$, up to re-pointing of its tree of tetrahedra $T_0$ and its tree of edges $E$.
 
 \ 
 
 Reciprocally, consider a  Morse gradient on $T$ up to re-pointing of its trees of tetrahedra $T_0$ and edges $E$. We wish to construct a collapsing sequence from $T^{T_0}$ onto $E$, that is,   we need to show that there exists a free edge and that removing it and the triangle it belongs to, this property stays true, and so on.

\begin{lemma} 
\label{lem:exist-free-edge}
Let $C$ be a 2-complex with a Morse gradient with no critical triangle.
Then $C$ has a free edge that is oriented towards the triangle it belongs to.
\end{lemma}
\proof Consider a triangle $\Sigma$. Since it is not critical, there is an edge $\sigma$ oriented towards it. Either $\sigma$ is a free edge, or it is contained in a triangle $\Sigma^{(1)}$ different from $\Sigma$ (otherwise there would be a cyclic walk), in which case $\Sigma^{(1)}$ being non-critical, it has an edge $\sigma^{(1)}$ different from $\sigma$ and oriented towards it. Either $\sigma^{(1)}$ is a free edge, or it is contained in a triangle $\Sigma^{(2)}$, and so on. At no point can $\Sigma^{(k_1)}$ coincide with  $\Sigma^{(k_2)}$, as otherwise there would be a cyclic walk. $C$  being finite, this sequence must stop somewhere, that is, $C$ must have a free edge, which is oriented towards the triangle it belongs to. \qed

\ 

We may therefore apply this lemma to $T^{T_0}$, choose a free edge, name it $\sigma_1$ and name $\Sigma_1$ the triangle it belongs to, remove them in an elementary collapse  $T^{T_0}\searrow C_1$.  $C_1$ still has a Morse gradient with no critical triangle, so we may apply Lemma~\ref{lem:exist-free-edge}, and so on until we have removed the $n$ triangles of $T^{T_0}$. The remaining space $C_n$ consists of all the vertices of $T^{T_0}$ as well as the critical edges, that is, $E$. The different choices of ordering of the elementary collapses correspond to the different choices of orderings of the moves in the local construction, reproducing the different weakly equivalent local constructions. This concludes the proof of the theorem. \qed

\

\begin{proposition}
\label{prop:unicity}
Let $C$ be a 2-complex and $E$ a spanning tree. If $C$ admits a Morse gradient whose critical simplices are the vertices and edges of $E$, then it is unique.  
\end{proposition}
As a consequence, the number of  weakly equivalence classes of local constructions of $T$ is bounded from above by its number of spanning trees of tetrahedra and spanning trees of edges. 
\proof 
We prove the statement inductively on the number of triangles of $C$. If $C$ has a single triangle, since it is not critical it has one edge oriented towards it, and since the vector field is acyclic, this edge cannot be identified with any other edge of the triangle. Because of the definition of a discrete vector field, the other two (possibly non-distinct) edges of the triangle must be critical, that is, they must form $E$. Since there is only one edge not in $E$ and it must be oriented towards the only triangle, the Morse gradient is unique.

We now assume that $C$ has more than one triangle and that it admits a Morse gradient whose critical simplices are the vertices and edges of $E$. Then from Lemma~\ref{lem:exist-free-edge}, it must admit a free edge $\sigma$ oriented towards a triangle $\Sigma$.
In any other discrete vector field in which $\sigma$ is not critical, it will also be oriented towards $\Sigma$. 
Removing $\sigma$ and $\Sigma$ from $C$, we obtain a smaller 2-complex $C'$ with a Morse gradient whose critical simplices are the vertices and edges of $E$. By induction, this Morse gradient on $C'$ is unique. Hence the Morse gradient on $C$ is as well, since $\sigma$ must be oriented towards $\Sigma$ in $C$. \qed

\section{Enumerative bounds}
\label{sec:combinatorial-bounds}

\subsection{Enumeration of tree-decorated hierarchical triangulations}

\begin{figure}[t]
	\centering
	\includegraphics[width=.9\linewidth]{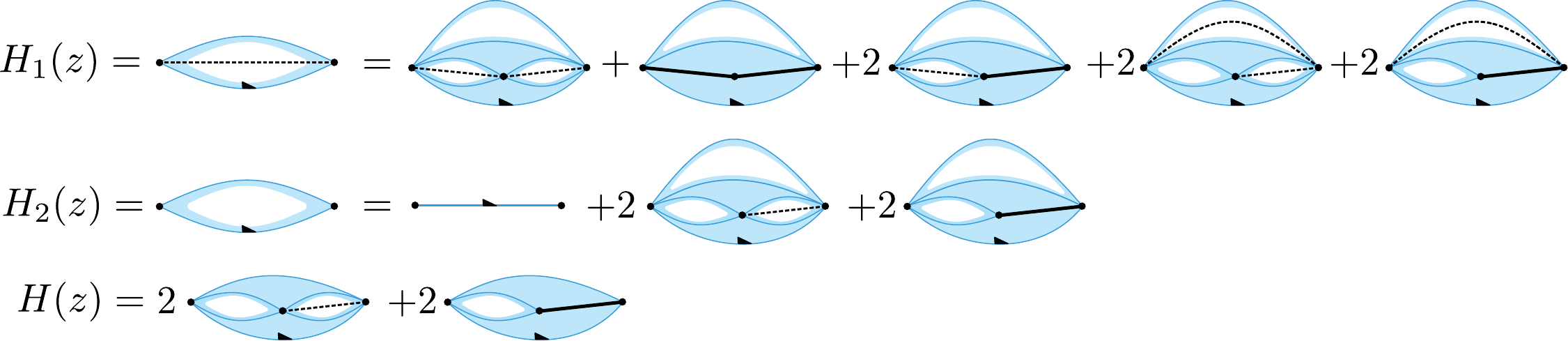}
	\caption{Pictorial representation of the system of equations satisfied by $H_1, H_2$. In each diagram the root face corresponds to the outside, the white regions represent arbitrary triangulations, and the dotted curves indicate the connectivity of the forest. The factors of $2$ account for the mirror images of the displayed diagrams.   \label{fig:hierarchicalrecurrence}}
\end{figure}

Recall that a hierarchical triangulation is a planar triangulation that can be obtained from the unique loopless triangulation with two triangles by repeatedly zipping open an edge and inserting a pair of triangles that are glued to each other along two of their sides.
The family $\mathcal{H}_n$ consists of rooted spanning-tree-decorated hierarchical triangulations $n+2$ vertices and the requirement that the root edge is in the spanning tree and that the spanning tree does not contain edges that are in length-two cycle.
In order to enumerate $\mathcal{H}_n$ we introduce two related families that relax some of these requirements.
The first family $\mathcal{H}_n^{1}$ differs from $\mathcal{H}_n$ only by the requirement that the root edge is not on the spanning tree. 
The second family $\mathcal{H}_n^{2}$ differs from $\mathcal{H}_n$ by allowing the root edge, which must be in the spanning tree, to belong to a cycle of length two.
By convention, $\mathcal{H}_2^1 = \emptyset$ while $\mathcal{H}_2^{2}$ contains the degenerate map with a single edge, two vertices, and no faces.
The point of introducing these two families is that they admit a recursive decomposition.
To see this, it is convenient to represent them as triangulations of a 2-gon, by zipping open the root edge such the 2-gon lies on the right of the root edge.
In the case of $\mathcal{H}_n^2$ the root edge is dropped from the spanning tree, resulting in a pair of (possibly empty) trees based at the endpoints of the root edge that together span the vertices of the map.
If the map does not consist of a single edge, then the triangle adjacent to the root edge and its companion triangle together delimit a triple of triangulations of 2-gons. 
Each one of these is a single edge in the spanning tree or belongs to one of $\mathcal{H}_{n'}^1$ or $\mathcal{H}_{n'}^2$.
The possibilities are shown in Fig.~\ref{fig:hierarchicalrecurrence}.
Denoting by $H_1(z)$ and $H_2(z)$ the generating functions of $\mathcal{H}_n^1$ and $\mathcal{H}_n^2$ the decomposition is easily seen to lead to the relations
\begin{align*}
	H_1(z) &= z \left( H_2 + 4 H_1 H_2 + 3 H_1^2 H_2 \right)\\
	H_2(z) &= 1 + 2 z\, H_2^2 ( 1 + H_1 ) \\
	H(z) &= 2z\, H_2 (1+H_1).
\end{align*}
These are straightforwardly solved to give
\begin{align*}
	H(z) = \frac{1-\sqrt{1-12 z}}{3}.
\end{align*}
We thus conclude the following.

\begin{proposition}\label{prop:hierarchicalenum}
	The generating function of the tree-decorated hierarchical triangulations $\mathcal{H}_n$ is 
	\begin{equation}
		H(z) = \frac{1-\sqrt{1-12 z}}{3} = 2z + 6z^2 + 36 z^3+\cdots.
	\end{equation}
	Explicitly, the enumeration and its large-$n$ asymptotics are given by
	\begin{equation}
		|\mathcal{H}_n| = 2\cdot 3^{n-2}\operatorname{Cat}(n-2) \stackrel{n\to\infty}{\sim} \frac{1}{72\sqrt{\pi}} \frac{12^n}{n^{3/2}}.
	\end{equation}
\end{proposition}

\subsection{Enumeration of Apollonian triangulations}

Recall that $\mathcal{A}_n$ is the family of rooted tree-decorated Apollonian triangulations with $n+1$ vertices such that the root edge is in the spanning tree.
We denote its generating function by $A(z) = \sum_{n=2}^\infty |\mathcal{A}_n|z^n = 2 z^2 + \cdots$.
In order to enumerate $\mathcal{A}_n$, we use a similar strategy and introduce three new generating functions.
Let $A_1(z)$, $6A_2(z)$ and $3 A_3(z)$ be the generating functions of rooted tree-decorated Apollonian triangulations, with weight $z^{n-3}$ in case of $n$ vertices, such that respectively zero, one or two of the edges in the root face (the face on the right of the root) are part of the tree.
Note that in general these triangulation do not occur in $\mathcal{A}_n$, since in the latter we required the root edge to be part of the tree, but we can easily express the generating function of $\mathcal{A}_n$ in terms of these as
\begin{equation}
	A(z) = 2 z^2 A_2(z) + 2 z^2 A_3(z).
\end{equation}
To obtain a system of equations for $A_i(z)$ it is convenient to reinterpret them in terms of triangulations decorated by a forest as follows.
Let us denote the three vertices in the root face by $v_1$ (start of the root edge), $v_2$ (end of the root edge) and $v_3$.
We consider rooted Apollonian triangulations decorated by a forest (i.e.\ an acyclic subgraph) such that none of the edges on the root face is part of the forest and such that every vertex can be reached starting from at least one of $v_1$, $v_2$, $v_3$.
Then $A_1(z)$ enumerates such triangulations for which $v_1$, $v_2$, $v_3$ are in the same tree, $A_2(z)$ the triangulations for which $v_1$, $v_2$ belong to the same tree but not $v_3$, and $A_3(z)$ the triangulations for which $v_1$, $v_2$, $v_3$ are all in different trees.

\begin{figure}[h!]
	\centering
	\includegraphics[width=.7\linewidth]{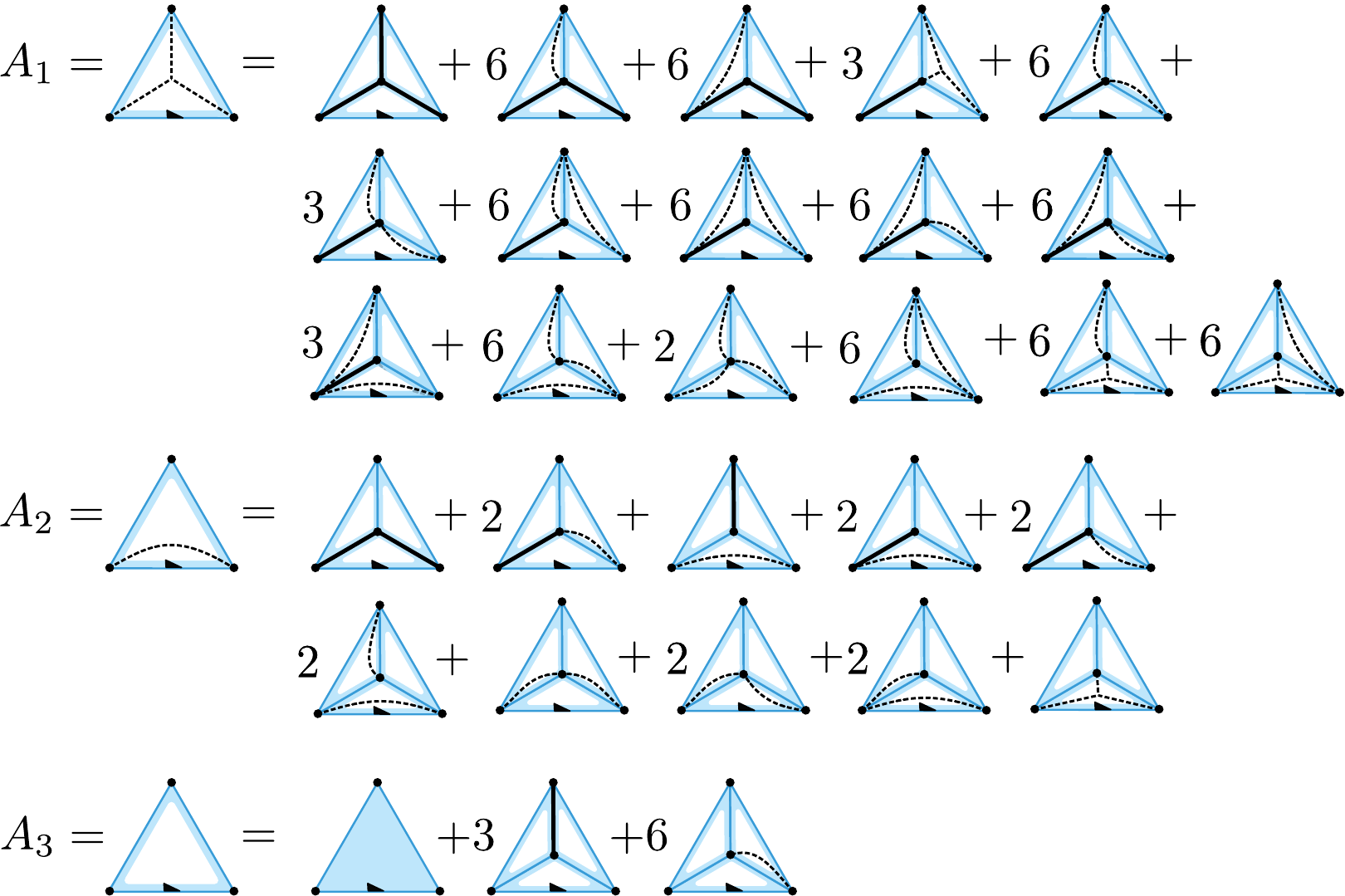}
	\caption{Pictorial representation of the system of equations satisfied by $A_1, A_2, A_3$. In each diagram the root face corresponds to the outside, the white regions represent arbitrary triangulations, and the dotted curves indicate the connectivity of the forest. The multiplicative factors account for the distinct rotations and mirror images of the displayed diagrams. \label{fig:appolonianrecurrence}}
\end{figure}
By examining the possibilities at the first star division (Figure \ref{fig:appolonianrecurrence}) one finds the system of equations
\begin{align*}
	A_1 &= z\left( A_3^3 + 12 A_2 A_3^2 + 3 A_1 A_3^2 + 36 A_2^2 A_3 + 14 A_2^3 + 12 A_1A_2A_3\right)\\
	A_2 &= z\left( A_3^3 + 7 A_2 A_3^2 + 7 A_2^2 A_3 + A_1 A_3^2\right)\\
	A_3 &= 1 + z\left( 3 A_3^3 + 6 A_2 A_3^2  \right)
\end{align*}
which uniquely determine them as formal power series in $z$. 
Eliminating $A_1$ and $A_2$ leads to the algebraic equation
\begin{equation*}
	81 z^2 A_3^6 - 54 z A_3^4 + 108 z A_3^3+4 A_3^3 + 9 A_3^2-48 A_3 + 35 = 0.
\end{equation*}
In particular we find the expansion
\begin{equation*}
	A(z) = 2 z^2 + 8 z^3 + 100 z^4 + 1680 z^5 + 32414 z^6 + 677810 z^7 + \cdots
\end{equation*}
Using singularity analysis one may easily deduce the asymptotics of the power series coefficients.
\begin{proposition}\label{prop:apollonianenum}
	The generating function $A(z)$ of the tree-decorated Apollonian triangulations $\mathcal{A}_n$ is algebraic. 
	The enumeration satisfies the asymptotics
	\begin{align*}
		|\mathcal{A}_n| \sim \frac{C^n}{c\,n^{3/2}},
	\end{align*}
	where $c= 1752.101\ldots$ and $C = 28.43330\ldots$ is the largest real solution to 
	\begin{equation*}
		1296 \,C^6+363232 \,C^5-16927248\, C^4+438097032 \,C^3-8977010004 \,C^2+28680825384 \,C-24111675 = 0.
	\end{equation*}
\end{proposition}

\subsection{A special class of triple trees}\label{sec:specialclass}

\begin{figure}[h!]
	\centering
	\includegraphics[width=.6\linewidth]{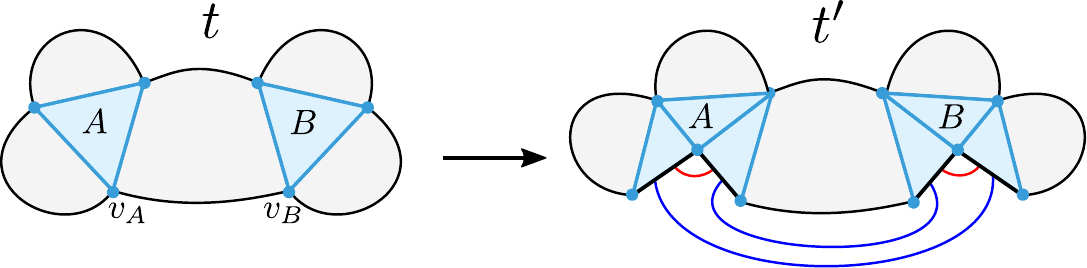}
	\caption{A schematic illustration of the outerplanar triangulation with two distinguished triangles $A$ and $B$. The grey regions represent arbitrary outerplanar triangulations or could be reduced a single edge. The blue and red arcs in the outer face of $t'$ represent the new pairs in ${\pi_{_\mathrm{H}}}'$ and ${\pi_{_\mathrm{A}}}'$ respectively. \label{fig:recursiveconstruction}}
\end{figure}

To obtain a first lower bound on $M_n(x) = \sum_{(t,\pi_{_\mathrm{H}},\pi_{_\mathrm{A}})\in\mathcal{M}_n} x^{N(\pi_{_\mathrm{H}},\pi_{_\mathrm{A}})}$, let us focus on an easy recursive construction of triple trees.
Suppose we have a triple tree $(t,\pi_{_\mathrm{H}},\pi_{_\mathrm{A}}) \in \mathcal{M}_n$ and let $A$, $B$ be a pair of triangles in $t$ that share three vertices in the hierarchical triangulation $h=\mathsf{Glue}(t,\pi_{_\mathrm{H}})$.
Let furthermore $v_A$ and $v_B$ be (not necessarily distinct) vertices of $t$ that are incident to $A$ and $B$ respectively, such that $v_A$ and $v_B$ are identified in $h$.
Then we can construct a new outerplanar triangulation $t' \in \mathcal{O}_{n+2}$ by cutting open $t$ from $v_A$ along each of the two sides of $A$ and inserting in each a new triangle with one side on the boundary and similarly for $B$, see Fig.~\ref{fig:recursiveconstruction}.
The non-crossing pairings $\pi_{_\mathrm{H}}$ and $\pi_{_\mathrm{A}}$ are extended to pairings ${\pi_{_\mathrm{H}}}'$ and ${\pi_{_\mathrm{A}}}'$ by adding the indicated pairs in blue and red respectively.
From the assumptions that $v_A$ and $v_B$ are identified in $h$ it follows that one can draw an arc in the outer face of $t$ from $v_A$ to $v_B$ that does not cross the arc system determined by $\pi_{_\mathrm{H}}$.
Hence ${\pi_{_\mathrm{H}}}'$ is non-crossing.
The two new triangles adjacent to $A$ in $t'$ each share three vertices in $\mathsf{Glue}(t',\pi_{_\mathrm{H}}')$ with their companion adjacent to $B$, implying that $\mathsf{Glue}(t',{\pi_{_\mathrm{H}}}') \in \mathcal{H}_{n+2}$ is hierarchical.
The other pairing ${\pi_{_\mathrm{A}}}'$ is obviously non-crossing and it should be clear that $\mathsf{Glue}(t',{\pi_{_\mathrm{A}}}')$ is obtained from $\mathsf{Glue}(t',{\pi_{_\mathrm{A}}}')$ by a star division of the triangles $A$ and $B$.
Hence $\mathsf{Glue}(t',{\pi_{_\mathrm{A}}}') \in \mathcal{A}_{n+2}$ and therefore we have found a new triple tree $(t',{\pi_{_\mathrm{H}}}',{\pi_{_\mathrm{A}}}')\in \mathcal{M}_{n+2}$.

From the point of view of the three-dimensional triangulation $T$ with spanning tree of tetrahedra $T_0$ and spanning tree of edges $E$, it corresponds to selecting a triangle that does not intersect $T_0$ and inserting two tetrahedra there that are glued along three of their sides, and extending $T_0$ and $E$ appropriately.

\begin{figure}[h!]
	\centering
	\includegraphics[width=.2\linewidth]{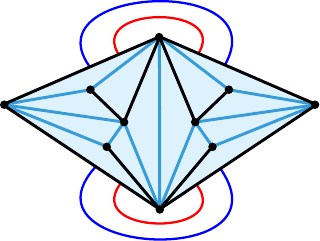}
	\caption{A possible result of applying the subdivision three times to a triple tree of size $n=2$. In this figure the outerplanar triangulation has been glued according to $\pi_{_A}$ except for the arcs indicated in red, while $\pi_{_H}$ pairs left-right mirror images. \label{fig:symmetrictripletree}}
\end{figure}

Since $t$ has $2n-2$ triangles, there are precisely $3n-3$ choices of $A$, $B$, $v_A$ and $v_B$ in the construction (modulo interchanging $A$ and $B$), which all lead to distinct triple trees.
Let us consider the family of triple trees obtained using this construction starting with one of the two triple trees in $\mathcal{M}_2$ shown in Fig.~\ref{fig:tripletreesn4n5}.
It should be clear that the mirror symmetry of the initial triple tree is preserved by the subdivision, as can be seen for the example of size $n=8$ shown in Fig.~\ref{fig:symmetrictripletree}.
The configuration is thus uniquely encoded by, say, its left half, which is a rooted Apollonian triangulation with $n/2$ triangles (excluding the outer face) in which at each star-division exactly one of the new edges is taken to be in the spanning tree. 
Since rooted Apollonian triangulations are in bijection with full ternary trees with $k=n/2-1$ nodes, their number is given by $\binom{3k+1}{k}/(3k+1)$. 
When allowing the special triple trees to be rooted arbitrarily, we thus find exactly 
\begin{equation*}
3^{k} \frac{2k+2}{3k+1}\binom{3k+1}{k}
\end{equation*}
triple trees of even size $n=2k+2$.
The number of loops of the meander system is easily seen to be $N(\pi_{_\mathrm{H}},\pi_{_\mathrm{A}}) = (n+2)/2 = k+2$, hence for $n$ even 
\begin{equation}
	M_n(x) \geq 3^{k} \frac{2k+2}{3k+1}\binom{3k+1}{k} x^{k+2}\quad\stackrel{n\to\infty}{\sim}\quad \frac{2}{27} \frac{(\tfrac92 \sqrt{x})^n}{\sqrt{\frac32 \pi n}}.
\end{equation} 

\end{document}